\documentclass{amsart}
\usepackage{amssymb,latexsym}
\usepackage{amsfonts}

\newcommand{\ba}{\begin{array}}
\newcommand{\ea}{\end{array}}

\def \qed{\cqfd}

\def \b{\beta}

\def\qed{\vbox{\hrule
\hbox{\vrule\hbox to 5pt{\vbox to 8pt{\vfil}\hfil}\vrule}\hrule}}

\newcommand{\beg}{\begin{eqnarray*}}
\newcommand{\begn}{\begin{eqnarray}}
\newcommand{\en}{\end{eqnarray*}}
\newcommand{\enn}{\end{eqnarray}}

\begin{document}
\title{The limit of the Yang-Mills-Higgs flow  on Higgs bundles}
\author{Jiayu Li}
\address{Key Laboratory of Wu Wen-Tsun Mathematics\\ Chinese Academy of Sciences\\School of Mathematical Sciences\\
University of Science and Technology of China\\
Hefei, 230026\\ and AMSS, CAS, Beijing, 100080, P.R. China\\} \email{jiayuli@ustc.edu.cn}
\author{Xi Zhang}
\address{Key Laboratory of Wu Wen-Tsun Mathematics\\ Chinese Academy of Sciences\\School of Mathematical Sciences\\
University of Science and Technology of China\\
Hefei, 230026,P.R. China\\ } \email{mathzx@ustc.edu.cn}
\thanks{The authors were supported in part by NSF in
China,  No.11131007.}

\begin{abstract} In this paper, we consider the gradient flow of the Yang-Mills-Higgs functional for Higgs pairs on a Hermitian vector bundle $(E, H_{0})$ over a compact K\"ahler manifold $(M, \omega )$. We study the asymptotic behavior of the Yang-Mills-Higgs flow for Higgs pairs at infinity, and show that the limiting Higgs sheaf is isomorphic to the double dual of the graded Higgs sheaves associated to the Harder-Narasimhan-Seshadri filtration of the initial Higgs bundle.
\end{abstract}

\maketitle

\section{Introduction}
\setcounter{equation}{0}

\hspace{0.4cm}

Let $(E, H_{0})$ be a Hermitian vector bundle
over a compact K\"ahler manifold $(M, \omega )$,
$ \textbf{A}_{H_{0}}$ be the space of connections of $E$
compatible with the metric $H_{0}$, and $ \textbf{A}^{1,1}_{H_{0}}$
be the space of unitary integrable connections of $E$ (i.e. those whose curvature is of type $(1, 1)$).
Given a unitary integrable connection
$A$ on $(E , H_{0})$ , then $D_{A}^{(0,1)}=\overline{\partial }_{A}$ defines a
holomorphic structure on $E$, and in fact, $A$ is the Chern connection on the holomorphic bundle $(E, \overline{\partial }_{A})$ with respect to $H_{0}$.

 A pair
$(A, \phi )\in  \textbf{A}^{1,1}_{H_{0}}\times \Omega^{1,0}(End(E))$
is called a Higgs pair if the relations
$\overline{\partial}_{A}\phi =0$ and $\phi \wedge \phi =0$ are
satisfied. Let $ \textbf{B}_{(E, H_{0})}$ denote the space of all
Higgs pairs on Hermitian vector bundle $(E , H_{0})$. Given a Higgs pair $(A, \phi )$, then $(E, \overline{\partial}_{A}, \phi )$ is a Higgs bundle, i.e. $(A, \phi )$ determines a Higgs structure on $E$.  Let us consider
 the Yang-Mills-Higgs functional which is defined on
$ \textbf{B}_{(E, H_{0})}$:
\begin{eqnarray}\label{YMHF}
YMH(A, \phi) =\int_{M}(|F_{A}+[\phi , \phi^{\ast}
]|^{2}+2|\partial_{A}\phi |^{2})\, dV_{g}.
\end{eqnarray}
We call $(A, \phi)$ a Yang-Mills Higgs pair if it is the critical points of
the above Yang-Mills-Higgs functional. Equivalently, the pair $(A,
\phi)$ satisfies the following Yang-Mills-Higgs equations:
\begin{eqnarray}\label{YMHE}
 \left \{\begin{array}{cll} &D_{A}^{\ast } F_{A}
+\sqrt{-1}(\partial_{A}\Lambda_{\omega}-\overline{\partial}_{A}\Lambda_{\omega})[\phi , \phi ^{\ast}]=0,\\
&[\sqrt{-1}\Lambda_{\omega}(F_{A}+[\phi , \phi^{\ast} ]), \phi ]=0,\\
\end{array}
\right.
\end{eqnarray}
where the operator $\Lambda_{\omega}$ is the contraction with
$\omega$, and $\phi^{\ast }$ denotes the dual of $\phi$ with respect to
the given metric $H_{0}$.

If $(A, \phi)$ satisfies the
following Hermitian-Einstein equation
\begin{eqnarray}\label{HEE}
\sqrt{-1}\Lambda_{\omega}(F_{A}+[\phi , \phi^{\ast} ])=\lambda
Id_{E},
\end{eqnarray}
then it must satisfy the above Euler-Lagrange equation (\ref{YMHE}). By Chern-Weil theory, in
fact it is the absolute minima of the above Yang-Mills-Higgs
functional.

On a Hermitian vector bundle $(E, H_{0})$, the Yang-Mills flow, as the gradient flow of the Yang-Mills functional, was first suggested by Atiyah-Bott in \cite{AB}.
Donaldson \cite{D1}  proved the global existence of the Yang-Mills flow in a
holomorphic bundle,  and proved the
convergence of the flow at infinity in the case that the
initial holomorphic structure is stable. He then used this fact to establish the correspondence between
existence of the Hermitian-Einstein metric and stability of the holomorphic structure  over complex algebraic surfaces. This correspondence was first shown by Narasimhan and Seshadri (\cite{NS}) in the case of compact Riemann surfaces, and is sometimes referred to as the Hitchin-Kobayashi
correspondence. The general K\"ahler manifold case  was proved by Uhlenbeck and Yau \cite{UY} by using the method of continuity.

Without the stability assumption,  Atiyah-Bott  \cite{AB} point out that there should be a correspondence between the Yang-Mills flow and the Harder-Narasimhan filtration over Riemann surfaces, and this is also conjectured by Bando and Siu \cite{BS} for the higher dimensional case. This correspondence was proved   by Daskalopoulos \cite{Da} in the case of Riemann surfaces,   by Daskalopoulos and Wentworth \cite{DW1} in the case of  K\"ahler surfaces. In the higher dimensional case, Hong and Tian \cite{HT} study the asymptotic behavior of the Yang-Mills flow, they proved that there is a subsequence along the Yang-Mills flow, modulo gauge transformations,  which converges smoothly to a limiting Yang-Mills connection away from the bubbling set $\Sigma_{an}$ of Hausdorff codimension $4$. Recently, Jacob in \cite{Ja3} and Sibley in \cite{Sib} studied the above correspondence for the higher dimension case.

Let us consider  the following  gradient flow  of the
Yang-Mills-Higgs functional of Higgs pairs, which is  called the Yang-Mills-Higgs flow. A regular solution is
given by a family of $(A(x,t), \phi (x,t))\in \textbf{B}_{(E, H_{0})}$
such that
\begin{eqnarray}\label{YMHH}
 \left \{\begin{array}{cll} &\frac{\partial A}{\partial t}=-D_{A}^{\ast } F_{A}
-\sqrt{-1}(\partial_{A}\Lambda_{\omega}-\overline{\partial}_{A}\Lambda_{\omega})[\phi , \phi ^{\ast}],\\
&\frac{\partial \phi }{\partial t}=-[\sqrt{-1}\Lambda_{\omega}(F_{A}+[\phi , \phi^{\ast} ]), \phi ].\\
\end{array}
\right.
\end{eqnarray}
It is interesting to study the Higgs pairs version of Atiyah-Bott's (or Bando-Siu's) conjecture, i.e.  there should be a correspondence between the Yang-Mills-Higgs flow for Higgs pairs and the Harder-Narasimhan-Seshadri filtration of Higgs bundle. In Riemann surface case, Wilkin \cite{Wi} develops the analytic results needed to construct a Morse theory for the Yang-Mills-Higgs functional on the space of all Higgs pairs, and proves the Higgs pairs version of Atiyah-Bott's conjecture. In \cite{LZ1}, the authors study the bubbling phenomena of the Yang-Mills-Higgs flow on K\"ahler surface case, and prove that the limit can be extended across the bubbling set $\Sigma_{an}$ (a finite collection of points) to a smooth Higgs bundle which isomorphic to the the double dual of the graded object of the Harder-Narasimhan-Seshadri filtration of the initial Higgs structure.

  A Higgs bundle $(E, \overline{\partial}_{E}, \phi )$ is
called  stable (semi-stable), if for every $\phi$-invariant coherent sub-sheaf
$E'\hookrightarrow (E, \overline{\partial }_{E})$ of lower rank, it holds:
\begin{eqnarray}
\mu (E')=\frac{deg (E')}{rank E'}< (\leq ) \mu (E)=\frac{deg (E)}{rank E},
\end{eqnarray}
where $\mu (E')$ is called the slope of $E'$. A Hermitian metric $H$ on $(E, \overline{\partial}_{E}, \phi )$ is called Hermitian-Einstein if it satisfies the  Hermitian-Einstein equation (\ref{HEE}), where $A$ is the Chern connection $A_{H}$ with respect to the metric $H$. Higgs bundles  first emerged twenty
years ago in  Hitchin's \cite{Hi} reduction of self-dual equation on $R^{4}$ to Riemann surface.
 Higgs bundles have a rich structure and
play a role in many different areas including gauge theory,
K\"ahler and hyperk\"ahler geometry, group representations and
nonabelian Hodge theory. In \cite{Si}, Simpson generalized it to the higher dimensional case  and proved that a Higgs bundle admits a Hermitian-Einstein metric if and only if it is Higgs poly-stable. This is  a
Higgs bundle version of the Donaldson-Uhlenbeck-Yau theorem.  This
correspondence has several interesting and important
generalizations and extensions where some extra structures are
added to the holomorphic bundles, see references: \cite{Hi}, \cite{Si},\cite{Br1},\cite{GP},
\cite{BG}, \cite{AG1},
\cite{Bi}, \cite{BT},  \cite{LN1}, \cite{LN2}, \cite{LY}, \cite{M}, \cite{Ta}.

\medskip

To a Higgs bundle $(E , \overline{\partial}_{E}, \phi )$ of rank $R$, as that for holomorphic bundles,  one can associate a
filtration by $\phi $-invariant holomorphic subsheaves, which is
called the Harder-Narasimhan filtration, whose successive
quotients are Higgs semi-stable. The topological type of the pieces in the associated graded objects is encoded into an $R$-tuple $\vec{\mu}=(\mu_{1}, \cdots , \mu_{R})$ of rational numbers called the Harder-Narasimhan type (abbr, HN-type ) of the  Higgs bundle $(E , \overline{\partial}_{E}, \phi )$.  For every semi-stable Higgs sub-sheaf, one can associate a Seshadri filtration, whose successive quotients are Higgs stable. Then, we have a double filtration which is called the Harder-Narasimhan-Seshadri filtration (abbr, HNS-filtration ) of the Higgs bundle, and we write $Gr^{HNS}(E , \overline{\partial}_{E}, \phi)$ for the associated graded object (i.e. the direct sum of the stable quotients) of the HNS filtration.

 Now, we consider  the  asymptotic behavior of the Yang-Mills-Higgs flow (\ref{YMHH}) of Higgs pairs for the higher dimensional case.
 In \cite{LZ1},  we  proved the global existence and
uniqueness of the solution for the Yang-Mills-Higgs flow, and obtained many basic
properties of the  flow, including the energy inequality,
Bochner-type inequality, monotonicities of certain quantities.  In this paper, we first give a
small action regularity estimate (Theorem 2.6.) for the higher dimensional case, which was proved  for the K\"ahler surface case in \cite{LZ1}.  Then,
following Hong-Tian's argument  in \cite{HT}, we conclude that there exists a sequence
of Higgs pairs along the solution of the Yang-Mills-Higgs flow (\ref{YMHH}) which converges, modulo gauge transformations, to a limiting Yang-Mills-Higgs pair  in $C^{\infty}_{loc}$ topology outside
the bubbling set $\Sigma_{an}\subset M$, where $\Sigma_{an}$ is a closed set of Hausdorff real codimension $4$. Furthermore, we show that the limiting $(E_{\infty}, A_{\infty}, \phi_{\infty})$ can be extended to the whole $M$ as a reflexive Higgs sheaf, and prove that this extended reflexive Higgs bundle is isomorphic to the the double dual of the graded object of the Harder-Narasimhan-Seshadri filtration of the initial Higgs bundle $(E, A_{0}, \phi_{0})$, i.e. we prove the Higgs version of  Atiyah-Bott's (or Bando-Siu's) conjecture for the higher dimension case. We obtain the following theorem.

\medskip

{\bf Theorem 1.1. } {\it Let $(E, H_{0})$  be a Hermitian  vector
bundle on a compact K\"ahler manifold $(M, \omega)$, and $(A(t),
 \phi (t) )$ be a global smooth solution of the Yang-Mills-Higgs flow (\ref{YMHH})  with smooth
initial Higgs pair $(A_{0}, \phi_{0})$.  Then:

 (1) There exists a sequence
$\{ t_{j}\}$ such that, as $t_{j}\rightarrow \infty$, $(A(t_{j}), \phi (t_{j}) )$ converges, modulo gauge transformations, to a Hermitian-Einstein
Higgs pair $(A_{\infty}, \phi_{\infty} )$ on Hermitian vector bundle $(E_{\infty}, H_{\infty})$ in $C_{loc}^{\infty}$ topology
outside a closed set $\Sigma^{an}\subset M$, where $\Sigma^{an}$ is
a closed set of Hausdorff codimension at least $4$.

(2)  The limiting $(E_{\infty}, \overline{\partial }_{A_{\infty}}, \phi_{\infty})$ can be extended to the whole $M$ as a reflexive Higgs sheaf  with a holomorphic orthogonal splitting
  \begin{eqnarray}
(E_{\infty}, H_{\infty}, A_{\infty},
\phi_{\infty})=\bigoplus_{i=1}^{l}(E_{\infty}^{i} , H_{\infty}^{i}
, A_{\infty}^{i} , \phi_{\infty} ^{i}),
\end{eqnarray}
where $H_{\infty}^{i}$ is an admissible Hermitian-Einstein metrics on the reflexive
Higgs sheaf $(E_{\infty}^{i} , A_{\infty}^{i} ,
\phi_{\infty}^{i})$.

(3) Moreover, the extended reflexive Higgs sheaf is isomorphic to the double dual of the graded Higgs sheaves associated to the HNS-filtration of the initial Higgs bundle, i.e. we have
\begin{eqnarray}
(E_{\infty}, \overline{\partial }_{A_{\infty}}, \phi_{\infty})\simeq Gr^{HNS}(E,
\overline{\partial }_{A_{0}},
 \phi_{0})^{\ast \ast }.
 \end{eqnarray}}

\medskip

We now give an overview of our proof. We (\cite{LZ1}) have already demonstrated  that there is an uniform $C^{0}$ bound on Higgs fields $\phi (t)$, so
the basic idea in \cite{DW1} for the Yang-Mills flow in the K\"ahler surface case can be used. But there are two points where we need new arguments for the higher dimensional case. The first one is to prove that the HN type of the limiting Higgs sheaf is in fact equal to the type of the initial Higgs bundle; and the second one is to construct a non-zero holomorphic map from any stable quotient Higgs sheaf in HNS filtration of initial Higgs bundle to the limiting Higgs sheaf.

The first idea is closely related to the existence of an $L^{p}$-approximate critical Hermitian metric (as defined in \cite{DW1}). Under the semi-stability assumption of the initial Higgs bundle $(E, A_{0}, \phi_{0})$ in \cite{LZ2}, we proved the existence of $L^{\infty}$-approximate metric by the heat flow method.  For the general case, we will use the cut-off argument by Daskalopoulos and Wentworth in \cite{DW1}. We use the resolution of singularities theorem of Hironaka \cite{Hio1}, and  write $\pi : \tilde{M}\rightarrow M$ as the composition of blow-ups involved in the resolution, then the pullback bundle $\pi^{\ast }E$ has a filtration by Higgs subbundles, which is precisely the HNS-filtration of the initial Higgs bundle away from the exceptional divisor $\tilde{\Sigma}$. The metric $\pi^{\ast}\omega $ is degenerated along the divisor $\tilde{\Sigma}$, and it can be approximated by a family of K\"ahler metrics $\omega_{\epsilon}$ on $\tilde{M}$ as that in \cite{BS}.  Since every pullback quotient bundle is stable with respect to K\"ahler metrics $\omega_{\epsilon}$ for small $\epsilon$, one can use Simpson's theorem (in \cite{Si}) to take the direct sum of the Hermitian-Einstein metrics on  quotient Higgs subbundles in the resolution. If one can get a flat Hermitian metric in the neighborhood of singularities, then one may use Daskalopoulos and Wentworth's cut-off argument (where the singularity set is a collection of finite points).  After getting some uniform estimates, and modifying this metric, one can show that its Hermitian-Einstein tensor is close to the HN type in the $L^{p}$ norm.  By pushing this metric down, one can obtain a smooth $L^{p}$-approximate critical Hermitian metric on the Higgs bundle $(E, A_{0}, \phi_{0})$. This idea was used by Sibley in \cite{Sib}. However, since the singularity set for the filtration is  complex codimension $2$ which is not necessary the collection of finite points, in general we can not get a flat Hermitian metric on a neighborhood of singularities. But we should point out that this is not a crucial issue. Using Sibley's good observation ( see Lemma 5.8) and by choosing any fixed Hermitian metric on a neighborhood of singularities, we can also obtain an uniform estimate, this is enough to obtain a smooth $L^{p}$-approximate critical Hermitian metric, see Proposition 5.11. for details.

For the second one, we use Donaldson's idea to construct a nonzero holomorphic map to the limiting bundle as the limit of the sequence of gauge transformations (by rescaled) defined by the flow. The difficulty is to prove that the limiting map is in fact non-zero, because we have no uniform $L^{\infty}$ bound on the mean curvature (i.e. $|\sqrt{-1}\wedge_{\omega}F_{A}|$) for subsheaves. If the singularities are finite points,  i.e. the K\"ahler surface case, we can follow the argument in \cite{DW1} by a complex analytic argument to get uniform $C^{0}$
estimate, and then prove the limiting holomorphic map is non-zero. This argument is not suitable for higher dimensional case, since we do not know whether the complement of the singular set has a strictly pseudo-concave boundary.
 Using the resolution of singularities, we consider the  pullback bundle which has a filtration by subbundles. Evolving the Hermtian metric by the Donaldson heat flow with respect to K\"ahler metric $\omega_{\epsilon}$, by the result in \cite{BS}, we get uniform $L^{\infty}$ bounds on the mean curvature of $H(t)$ for positive $t$. By uniform local $C^{0}$ estimate of the evolved Hermitian metrics and using the standard elliptic estimates, we can construct a nonzero holomorphic map which we need (see Proposition 4.1. for detail). In \cite{Ja3} and \cite{Sib}, the authors studied the same question for holomorphic vector bundles, they have good observations there. We should point out that our argument is different from the ones they used.

This paper is organized as follows. In Section 2, we recall some basic estimates for Donaldson's heat flow and the Yang-Mills-Higgs flow, and prove the first and second part of Theorem 1.1., see Theorem 2.7. and Proposition 2.10..  In section 3, we consider the resolution of the HNS filtration of Higgs bundle.  In section 4, we construct a non-zero holomorphic map between Higgs sheaves, where Proposition 4.1. is the key technical part in the proof of Theorem 1.1.. In section 5, we use the cut-off argument to obtain $L^{p}$-approximate critical Hermitian metric, and prove that the Harder-Narasimhan type of the limiting Higgs sheaf is in fact equal to the type of the initial Higgs bundle.  In section 6, we complete the proof of Theorem 1.1. using an inductive argument.

\hspace{0.3cm}

\section{Analytic preliminaries and basic estimates }
\setcounter{equation}{0}

Let Let $(M, \omega )$ be a compact K\"ahler manifold of complex dimension $n$, and $(E, \overline{\partial}_{E}, \phi )$  be a Higgs
bundle on $M$. Suppose $H(t)$ is a solution of the  following Donaldson's heat flow with initial metric $H_{0}$,
\begin{eqnarray}\label{D1}
H^{-1}\frac{\partial H}{\partial
t}=-2(\sqrt{-1}\Lambda_{\omega}(F_{H}+[\phi , \phi ^{\ast H}
])-\lambda Id_{E}).
\end{eqnarray}
Let
$h(t)=H_{0}^{-1}H(t)$, using the identities
\begin{eqnarray}\label{id1}
\begin{array}{lll}
&& \partial _{H}-\partial_{H_{0}} =h^{-1}\partial_{H_{0}}h ;\\
&& F_{H}-F_{H_{0}}=\overline{\partial }_{E} (h^{-1}\partial_{H_{0}} h) ;\\
&& \phi^{\ast H}=h^{-1} \phi^{\ast H_{0}} h ,\\
\end{array}
\end{eqnarray}  then we can rewrite (\ref{D1}) as
\begin{eqnarray}
\frac{\partial h}{\partial
t}=-2\sqrt{-1}h\Lambda_{\omega}(F_{H_{0}}+\overline{\partial
}_{E}(h^{-1}\partial_{H_{0}}h)+[\phi, h^{-1}
\phi^{\ast H_{0}} h])+2\lambda h.
\end{eqnarray}

In \cite{Si}, Simpson proved the existence of long time solution of the heat flow (\ref{D1}).
 The following lemma is essentially proved by Simpson (\cite{Si} Lemma 6.1).

\medskip

{\bf Lemma 2.1. } {\it Let $H(t)$ be a solution of the  heat flow (\ref{D1}) with initial metric $H_{0}$, then we have:
\begin{eqnarray}\label{F1}
(\frac{\partial }{\partial t}-\triangle )tr (\Lambda_{\omega}(F_{H}+[\phi , \phi^{\ast H}
]))=0
\end{eqnarray}
and
\begin{eqnarray}\label{F2}
(\frac{\partial }{\partial t}-\triangle )|\Lambda_{\omega}(F_{H}+[\phi , \phi^{\ast H}
])|_{H}^{2}=-4|D''_{ \phi} (\Lambda_{\omega}(F_{H}+[\phi , \phi^{\ast H}
]))|^{2}_{H},
\end{eqnarray}
where $D''_{ \phi}=\overline{\partial}_{E} +\phi $.}

  \medskip

Denote the complex gauge group (unitary gauge group) of the Hermitian
vector bundle $(E, H_{0} )$ by $\textbf{G}^{C}$ ($\textbf{G}$, where
$\textbf{G}=\{\sigma \in \textbf{G}^{C}| \sigma^{\ast H_{0}}\sigma
=Id\}$). $\textbf{G}^{C}$ acts on the space of Higgs pairs $\textbf{B} _{(E, H_{0})}$ as follows: let
$\sigma \in \textbf{G}^{C}$
\begin{eqnarray}\label{id2}
\overline{\partial }_{\sigma(A)}=\sigma \circ \overline{\partial
}_{A}\circ \sigma^{-1}, \quad \partial _{\sigma (A)}=(\sigma^{\ast
H_{0}})^{-1} \circ \partial _{A}\circ \sigma^{\ast H_{0}};
\end{eqnarray}
\begin{eqnarray}\label{id3}
\sigma (\phi )=\sigma \circ \phi \circ \sigma^{-1}.
\end{eqnarray}

\medskip

In \cite{LZ1}, we  proved the following proposition.

\medskip

{\bf Proposition 2.2. (Theorem 2.1 in \cite{LZ1}) } {\it
Given any Higgs pair $(A_{0} , \phi_{0})$, the Yang-Mills-Higgs flow (\ref{YMHH}) has a unique solution $(A(t), \phi (t))$ in the complex gauge orbit of $(A_{0} , \Phi_{0})$. In fact, $(A(t), \phi (t))=g(t)(A_{0}, \phi_{0})$, where  $g(t) \in \textbf{G}^{C}$ satisfies $g^{\ast H_{0}}(t)g(t)=H_{0}^{-1}H(t)$, and $H(t)$ is the solution of Donaldson's flow (\ref{D1}) on Higgs bundle $(E, \overline{\partial }_{A_{0}}, \phi_{0})$ with initial metric $H_{0}$.}

 \medskip

Furthermore, we have the following Bochner type inequality (\cite[p.~1384]{LZ1})
\begin{eqnarray}\label{8}
\Big(\triangle-\frac{\partial}{\partial t}\Big)| \phi |^2\geq2| \nabla_A \phi|^2+\tilde{C}_1(|\phi|^2+1)^2-\tilde{C}_2(|\phi|^2+1),
\end{eqnarray}
where constants $\tilde{C}_{1}$ and $\tilde{C}_{2}$ depend only on the geometry of $(X,\omega)$ and the initial data $(A_0,\phi_0)$. By the maximum principle, we have the following uniform $C^{0}$ estimate for $\phi (t)$.

\medskip

{\bf Lemma 2.3. (Lemma 2.3. in \cite{LZ1}) } {\it Let $( A(t) , \phi (t) )$ be a solution of the Yang-Mills-Higgs
flow (\ref{YMHH}) with initial Higgs pair $( A_{0} , \phi_{0} )$, then we
have
\begin{eqnarray}\label{9}
|\phi (x, t)|_{H_{0}}^{2}\leq C ,
\end{eqnarray}
where $C$ is a constant depending only on
$\phi_{0}$ and the geometry of $(M , \omega )$. }

\medskip

For simplicity, we set
\begin{eqnarray}
\theta (A, \phi ) =\frac{1}{2\pi}\Lambda_{\omega}(F_{A}+[\phi  , \phi  ^{\ast H_{0}}]),
\end{eqnarray}
and
\begin{eqnarray}
I(t)=\int_{M}|D_{A(t)}\theta (A(t), \phi (t))|_{H_{0}}^{2}+2|[\theta (A(t), \phi (t)), \phi (t) ]|_{H_{0}}^{2}\frac{\omega ^{n}}{n!}.
\end{eqnarray}

Let $( A(t) , \phi (t) )$ be the solution of the heat
flow (\ref{YMHH}) on the Hermitian bundle $(E, H_{0})$ with initial Higgs pair $( A_{0} , \phi_{0} )$, and $H(t)$ be the solution of the Donaldson's flow (\ref{D1}) on Higgs bundle $(E, \overline{\partial }_{A_{0}}, \phi_{0})$ with initial metric $H_{0}$. As above, we know that $( A(t) , \phi (t) )=g(t)(A_{0}, \phi_{0})$, where $g(t)\in \textbf{G}^{C}$ and satisfies $g(t)^{\ast H_{0}} g(t)=h(t)=H_{0}^{-1}H(t)$. By (\ref{id1}), (\ref{id2}), (\ref{id3}), it is easy to check that
\begin{eqnarray}\label{id4}
F_{A(t)}=g(t)\circ F_{H(t)} \circ g(t)^{-1},
\end{eqnarray}
\begin{eqnarray}\label{id5}
[\phi (t) , \phi (t)^{\ast H_{0}}]=g(t)\circ [\phi_{0}, \phi_{0}^{\ast H(t)}]\circ g(t)^{-1},
\end{eqnarray}
and
\begin{eqnarray}\label{id6}
|\theta (A(t), \phi (t) )|_{H_{0}}=\frac{1}{2\pi}|\Lambda_{\omega}(F_{H(t)}+[\phi_{0} , \phi_{0}^{\ast H(t)}
])|_{H(t)},
\end{eqnarray}
where $F_{H(t)}$ is the curvature of the Chern connection on $(E, \overline{A_{0}})$ with respect to the metric $H(t)$.
Direct calculation shows that
\begin{eqnarray}
I(t)\rightarrow 0, \quad (t\rightarrow \infty),
\end{eqnarray}
the  proof can be found in \cite{LZ1} (page 1384-1386).

Furthermore, we recall the monotonicity inequality (Theorem 2.6 in \cite{LZ1} ) for the solution $(A(t), \phi (t))$ of (\ref{YMHH}).  For any point $x_0\in X$, there exists complex normal coordinates $\{z_1, \ldots, z_n\}$ in the geodesic ball $B_r(x_0)$ with center at $x_0$ and radius $r\leq i_X$ ($i_X$ is the infimum of the injectivity radius), such that
$$|g_{i\bar{j}}(z)-\delta_{ij}|\leq C|z|^2,\ \ \ \ \Big| \frac{\partial g_{i\bar{j}}}{\partial z_k}\Big|\leq C|z|$$
 for any $z\in B_{r}(x_0)$, where $(g_{i\bar{j}})$ is given by $g_{i\bar{j}}=g(\frac{\partial}{\partial z_i},\frac{\partial}{\partial z_{\bar{j}}})$ and $C$ is a constant which only depends on $x_{0}$.
  For a fixed point $u_0=(x_0,t_0)\in X\times \mathbb{R}_+$, we denote
\begin{equation}
\begin{split}
 &T_r(x_0,t_0)=\left\{u=(x,t):t_0-4r^2<t<t_0-r^2, x\in X\right\},
 \\& P_r(u_0)=B_r(x_0)\times[t_0-r^2,t_0+r^2].
\end{split}
\end{equation}
\par The fundamental solution of (backward) heat equation with singularity at $(z_0,t_0)$ is
\begin{equation}
G_{(z_0,t_0)}(z,t)=\frac{1}{(4\pi(t_0-t))^n}\exp{\Big(-\frac{|z-z_0|^2}{4(t_0-t)}\Big)},\ \ \ (t<t_0).
\end{equation}

Assume $(A(t),\phi(t))$ is a solution of the heat flow (\ref{YMHH}) with initial value $(A_0,\phi_0)$. Let $\varphi\in C_0^\infty(B_{i_X}(x_{0}))$ be a smooth cut-off function such that $\varphi\equiv1$ on $B_{i_X/2}(x_0)$, $ |\varphi|\leq1$, and $|\nabla \varphi| \leq 4/i_X$ in $B_{i_X}(x_0)\setminus B_{i_X/2}(x_0)$. Then we set
\begin{eqnarray}
e(A,\phi)(x,t)=|F_A+[\phi,\phi^{*H_0}]|_{H_{0}}^2+2|\partial_A\phi|_{H_{0}}^2
\end{eqnarray}
and
\begin{eqnarray}
\Phi(r;A,\phi)=r^2\int_{T_r(x_0,t_0)}e(A,\phi)\varphi^2G_{u_0}dV_g\,dt.
\end{eqnarray}
 We have the following monotonicity inequality.

\medskip

{\bf Proposition 2.4. (Theorem 2.6 in \cite{LZ1}} {\it
Let $(A(t),\phi(t))$ be a solution of the heat flow (\ref{YMHH}) with initial value $(A_0,\phi_0)$. Then, for any $(x_0,t_0)\in X\times [0,T]$ and  $0<r_1\leq r_2\leq\min\{i_X,\sqrt{t_0}/2\}$, we have
\begin{equation}\label{mo}
\Phi(r_1;A,\phi)\leq C\exp(C(r_2-r_1))\Phi(r_2;A,\phi)+C(r_2^2-r_1^2)\mbox{\rm{YMH}}(A_0,\phi_0)
\end{equation}
where $C$ is a positive constant depending only on the geometry of $(X, \omega )$.}

\medskip

For further consideration, we give an estimate for the integral of $|\nabla_A\phi|^2$ over $P_R(u_0)$, where $R\leq\min\{\sqrt{t_0},i_X/2\}$.

\medskip

{\bf Lemma 2.5. } {\it
 Let $(A(t),\phi(t))$ is a smooth solution of (\ref{YMHH}) with initial data $(A_{0},\phi_{0})$. Then, for any point  $u_0=(x_0,t_0) \in X\times\mathbb{ R}_+$ and $R\in\left(0,\min\{\sqrt{t_0},i_X/2\}\right]$, we have
\begin{eqnarray}\label{phi}
\int_{P_R(u_0)}|\nabla_A\phi|^2 dv_{g} dt\leq C R^{2n},
\end{eqnarray}
where  $C$ is a constant  depending only on the geometry of $(X,\omega)$ and the initial data $(A_0,\phi_0)$.}

\begin{proof}
Let $f\in C_0^\infty(B_{2R}(x_0))$ be a smooth cut-off function such that $f\equiv 1$ on $B_{R}(x_0)$, $| f|\leq1$ and $|\nabla f|\leq 4/R,\,|\triangle f|\leq 16/R^2$ in $B_{2R}(x_0)\setminus B_{R}(x_0)$.
By the Bochner type inequality (\ref{8}) and  the $C^{0}$-estimate (\ref{9}), we have
\begin{equation}\label{20}
\begin{split}
(\triangle-\frac{\partial}{\partial t})f^2|\phi|^2=&4f\langle\nabla f,\nabla| \phi|^2\rangle +2| \nabla f|^2|\phi|^2+2f\triangle f| \phi|^2+f^2(\triangle-\frac{\partial}{\partial t})| \phi|^2
\\ \geq&-8| f| |\nabla f| |\nabla_A\phi| |\phi|+2|\nabla f|^2|\phi|^2+2f\triangle f|\phi|^2
\\&+2f^2|\nabla_A\phi|^2-C_{2}(|\phi |^{2}+1)f^2
\\ \geq &-14|\nabla f|^2|\phi|^2+f^2|\nabla_A\phi|^2
\\&+2f\triangle f |\phi|^2-C_{2}(|\phi |^{2}+1)f^2
\\ \geq &f^2|\nabla_A\phi|^2 -C_{4}-C_{5}(|\nabla f|^{2}+|\triangle f|),
\end{split}
\end{equation}
where  $C_{4}$ and $C_{5}$ are constants depending only on the geometry of $(X,\omega)$ and the initial data $(A_0,\phi_0)$.

Integrating both sides of the inequality (\ref{20}) over $X\times [t_0-R^2,t_0+R^2]$, we have
\begin{equation*}
\begin{split}
&\int_{t_0-R^2}^{t_0+R^2}\int_Xf^2|\nabla_A\phi|^2
\\ \leq &\int_{t_0-R^2}^{t_0+R^2}\int_X\Big\{\Big(\triangle-\frac{\partial}{\partial t}\Big)f^2|\phi|^2+C_{4}+C_{5}(|\nabla f|^{2}+|\triangle f|)\Big\}dV_g\,dt
\\ \leq &\int_{B_{2R}(x_0)}|\phi(\cdot , t_0-R^2)|^2dV_g+C_{6}R^{2n}+C_{6}R^{2n+2}\leq C_7R^{2n},
\end{split}
\end{equation*}
where constant $C_7$ depends only on the geometry of $(X,\omega)$ and the initial data $(A_0,\phi_0)$.
\end{proof}

\medskip

Using the above monotonicity inequality (\ref{mo}) and the estimate (\ref{phi}) of $|\nabla_{A}\phi |^{2}$, we can obtain the following $\epsilon $-regularity theorem for the Yang-Mills-Higgs flow (\ref{YMHH}). The argument is very similar with the K\"ahler surface case (Theorem 3.1. in \cite{LZ1}), but  there is a difference in handling the term $|\nabla_{A}\phi |^{2}$. So we  give a proof in detail for reader's convenience.

\medskip

{\bf Theorem 2.6. ($\epsilon $-regularity theorem)} {\it
 Let $(A,\phi)$ be a smooth solution of (\ref{YMHH}) over an $n$-dimensional compact K\"ahler manifold $(X,\omega)$ with initial value $(A_0,\phi_0)$. There exist positive constants $\epsilon_0,\delta_0<1/4$, such that, if for some $0<R<\min\{i_{X}, \frac{\sqrt{t_0}}{2}\}$, the inequality
\begin{equation}
  R^{2-2n}\int_{P_{R}(x_{0},t_{0})}e(A(t),\phi(t))dV_g\,dt\leq \epsilon_0
\end{equation}
holds, then for any $ \delta\in(0,\, \delta_0)$, we have
\begin{equation}\label{k1}
\sup\limits_{P_{\delta R }(x_{0},t_{0})}e(A(t),\phi(t)) \leq 16(\delta R)^{-4},
\end{equation}
and
\begin{equation}\label{k2}
\sup\limits_{P_{\delta R}(x_{0},t_{0})}|\triangledown_{A(t)}\phi(t)|^{2}\leq C,
\end{equation}
where $C$ is a constant depending only on the geometry of $(X,\omega)$, the initial data $(A_0,\phi_0)$, $\delta_0$ and $R$.
}

\medskip

\begin{proof}
For any $ \delta\in (0,\delta_0]$, we construct a function \textit{f}: $[0,2\delta R]\longrightarrow \mathbb{R}$
\begin{eqnarray}
f(r)=(2\delta R-r)^4\sup\limits_{P_r(x_0,t_0)}e(A,\phi).
\end{eqnarray}
Since $f(r)$ is continuous and $f(2\delta R)=0$, we can assume that $f(r)$ attains its maximum at a certain $r_0\in[0,2\delta R)$. There exists a point $(x_1,t_1)\in \overline{P}_{r_0}(x_0,t_0)$, such that $$e(A,\phi)(x_1,t_1)=\sup\limits_{P_{r_0}(x_0,t_0)}e(A,\phi).$$

We claim that when $\epsilon_0,\,\delta_0$ are small enough, $f(r_0)\leq16$. Otherwise, set
$$\rho_0:=\frac{(2\delta R-r_0)}{\sqrt[4]{f(r_0)}}=e(A,\phi)(x_1,t_1)^{-1/4}<\frac{2\delta R-r_0}{2}.$$
Rescaling the Riemannian metric $\widetilde{g}=\rho_0^{-2}g$ and $t=t_1+\rho_0^2 \widetilde{t}$, we have
\begin{equation*}
\begin{split}
&|F_A+[\phi,\phi^{*H_0}]|^2_{\widetilde{g}}=\rho_0^4|F_A+[\phi,\phi^{*H_0}]|^2_g,
\\ &|\partial_A\phi|_{\widetilde{g}}^2=\rho_0^4|\partial_A\phi|_g^2,\ \ \
|\nabla_A\phi|_{\widetilde{g}}=\rho_0^4|\nabla_A\phi|^2_g.
\end{split}
\end{equation*}
Setting
\begin{equation*}
\begin{split}
e_{\rho_0}(x,\tilde{t}):=&|F_A+[\phi,\phi^{*H_0}]|^2_{\widetilde{g}}+2|\partial_A\phi|_{\widetilde{g}}^2=
 \rho_0^4e(A,\phi)(x,t_1+\rho_0^2\tilde{t}),
\\b_{\rho_0}(x,\tilde{t}):=&|\nabla_A\phi|_{\widetilde{g}}^2=\rho^4|\nabla_A\phi|^2_g(x,t_1+\rho_0^2\tilde{t}),
\\ \tilde{P}_{\tilde{r}}(x_1,0):=&B_{\rho_0\tilde{r}}(x_1)\times [-\tilde{r}^2,\tilde{r}^2],
\end{split}
\end{equation*}
we have $e_{\rho_0}(x_1,0)=\rho_0^4e(a,\phi)(x_1,t_1)=1$, and
\begin{equation*}
\begin{split}
\sup\limits_{\tilde {P}_1(x_1,0)}e_{\rho_0}&=\rho_0^4\sup\limits_{P_{\rho_0}(x_1,t_1)}e(A,\phi)
\leq\rho_0^4\sup\limits_{P_{2\delta R+r_0/2}(x_0,t_0)}e(A,\phi)
\\&\leq\rho_0^4f(r_0)(\frac{2\delta R-r_0}{2})^{-4}
=16.
\end{split}
\end{equation*}

This implies that
\begin{eqnarray}\label{FF}
|F_A+[\phi,\phi^*]|_{\widetilde{g}}+2|\partial_A\phi|_{\widetilde{g}}^2\leq 16,\ \ \ \mbox{on}\  \tilde{P}_1(x_1,0).
\end{eqnarray}

 Using the Bochner type inequality $(2.11)$ in \cite{LZ1}
\begin{equation}\label{12}
\begin{split}
&\Big(\triangle_g-\frac{\partial}{\partial t}\Big)|\nabla_A\phi|_g^2-2|\nabla_A\nabla_A\phi|_g^2
\\ \geq&-C(n)(|F_A|_g+|Rm|_g+|Ric_X|_g+|\phi|_g^2)|\nabla_A\phi|_g^2\\&-C(n)|\phi|_g|\nabla Ric_X|_g|\nabla_A \phi|_g,
\end{split}
\end{equation}
and inequality (\ref{9}), we have
\begin{equation}
\begin{split}
&\Big(\frac{\partial}{\partial\tilde{t}}-\triangle_{\tilde{g}}\Big)(b_{\rho_0}+\rho_0^4)
=\rho_0^6(\frac{\partial}{\partial t}-\triangle_g)(|\nabla_A \phi|^2_g+1)
\\ \leq &C(n)\rho_0^6(|F_A|_g+|Rm|_g+|Ric_X|_g+|\phi|_g^2)|\nabla_A \phi|_g^2\\&+C(n)\rho_0^6|\phi|_g|\nabla Ric_X|_g|\nabla_A \phi|_g
\\ \leq &C_8(b_{\rho_0}+\rho_0^4),
\end{split}
\end{equation}
where $C_8$ is a positive constant depending only on the initial data $(A_0,\phi_0)$ and the geometry of $(X,\omega)$. Then using the  parabolic mean value inequality (Theorem 14.7 in \cite{PL}) and the Lemma 2.5., we have
\begin{equation}\label{36}
\begin{split}
\sup\limits_{\tilde{P}_{1/2}(x_1,0)}(b_{\rho_0}+\rho_0^4)\leq &C^*\int_{\tilde{P}_1(x_1,0)}(b_{\rho_0}+\rho_0^4)dV_{\tilde{g}}d\tilde{t}
\\=&C^*\rho_0^{2-2n}\int_{P_{\rho_0}(x_1,t_1)}(|\nabla_A \phi|_g^2+1)dV_{g}dt
\\ \leq& C\rho_0^2\leq C_9
\end{split}
\end{equation}
where $C_9$ is a positive constant depending only on the initial data $(A_0,\phi_0)$ and the geometry of $(X,\omega)$.
By using the Bochner type inequality (2.12) in \cite{LZ1}
\begin{equation}\label{11}
\begin{split}
&\Big(\triangle_g-\frac{\partial}{\partial t}\Big)(|F_{A(t)}+[\phi(t),\phi(t)^{*H_0}]|_g^{2}+2|\partial_{A(t)}\phi(t)|_g^{2})
\\&-2|\nabla_{A(t)}(|F_{A(t)}+[\phi(t),\phi(t)^{*H_0}]|_g^{2})|-4|\nabla_{A(t)}(\partial_{A(t)}\phi(t))|_g^{2}
\\ \geq&-C(n)(|F_{A(t)}+[\phi(t),\phi(t)^{*H_0}]|_g+|\nabla_{A(t)}\phi(t)|_g+|\phi(t)|_g^{2}+|Rm|_g)
\\&(|F_{A(t)}+[\phi(t),\phi(t)^{*H_0}]|_g^{2}+2|\partial_{A(t)}\phi(t)|_g^{2}),
\end{split}
\end{equation}
 (\ref{FF}), (\ref{9}) and (\ref{36}), we have
\begin{equation*}
\begin{split}
\Big(\frac{\partial}{\partial \tilde{t}}-\triangle_{\tilde{g}}\Big)e_{\rho_0}=&\rho_0^6(\frac{\partial}{\partial t}-\triangle_g)e(A,\phi)
\\ \leq &\rho_0^6C(n)(|F_{A}+[\phi,\phi^{*H_0}]|_g+|\nabla_{A}\phi|_g+|\phi|_g^{2}+|Rm|_g)
\\&(|F_{A}+[\phi,\phi^{*H_0}]|_g^{2}+2|\partial_{A}\phi|_g^{2})
\\ \leq&C_{10} e_{\rho_0},
\end{split}
\end{equation*}
where $C_{10}$ is a positive constant depending only on $i_X$, the initial data $(A_0,\phi_0)$ and the geometry of $(X,\omega)$. Using the parabolic mean value inequality again, we have
\begin{equation}\label{T1}
\begin{split}
1=e_{\rho_0}(x_1,0)&\leq\sup\limits_{\tilde{P}_{1/4}(x_1,0)}e_{\rho_0}(x,\tilde{t})
\leq C\int_{\tilde{P}_{1/2}(x_1,0)}e_{\rho_0}dV_{\tilde{g}}\,d\tilde{t}
\\&\leq C_{11}\rho_0^{2-2n}\int_{P_{\rho_0}(x_1,t_1)}e(A,\phi)dV_g\,dt,
\end{split}
\end{equation}
where $C_{11}$ is a positive constant depending only on the initial data $(A_0,\phi_0)$ and the geometry of $(X,\omega)$.
 \par
We choose normal complex coordinates centred at $x_1$, and a smooth cut-off function $\varphi\in C_0^\infty(B_{R/2}(x_1))$, such that $\varphi\equiv1$ on $B_{R/2}(x_1)$, $|\varphi|\leq 1$ and $|\nabla\varphi|\leq 8/R$ on $B_R(x_1)\setminus B_{R/2}(x_1)$. Let $r_1=\rho_0$ and $r_2=\delta_0R$, using the monotonicity inequality (\ref{mo}) and the properties of the fundamental solution $G$,  we have
\begin{equation}\label{T2}
\begin{split}
&\rho_0^{2-2n}\int_{P_{\rho_0}(x_1,t_1)}e(A,\phi)dV_g\,dt
\\ \leq &C\rho_0^2\int_{P_{\rho_0}(x_1,t_1)}e(A,\phi)G_{(x_1,t_1+2\rho_0^2)}\varphi^2dV_g\,dt
\\ \leq&C\rho_0^2\int_{T_{\rho_0}(x_1,t_1+2\rho_0^2)}e(A,\phi)G_{(x_1,t_1+2\rho_0^2)}\varphi^2dV_g\,dt
\\ \leq&Cr_2^2\int_{T_{r_2}(x_1,t_1+2\rho_0^2)}e(A,\phi)\varphi^2G_{(x_1,t_1+2\rho_0^2)}dV_g\,dt
+C\delta_0^2R^2\mbox{\rm{YMH}}(A_0,\phi_0)
\\ \leq&C\delta_0^{2-2n}R^{2-2n}\int_{P_R(x_0,t_0)}e(A,\phi)dV_g\,dt
    +C\delta_0^2R^2\mbox{\rm{YMH}}(A_0,\phi_0)
\\ \leq&C_{12}(\epsilon_0\delta_0^{2-2n}+\delta_0),
\end{split}
\end{equation}
where $C_{12}$ is a positive constant depending only on the geometry of $(X,\omega)$ and the initial data $(A_0,\phi_0)$.

\par
From (\ref{T1}) and (\ref{T2}), we have the inequality $1\leq C_{11}C_{12}(\epsilon_0\delta_0^{2-2n}+\delta_0)$. Choosing $\epsilon_0,\,\delta_0$ small enough, we get a contradiction. Thus we have proved the claim. So, we have
\begin{eqnarray}\label{ee}
(2\delta R-r)^4\sup\limits_{P_r(x_0,t_0)}e(A,\phi)\leq 16
\end{eqnarray}
for any $\delta\in (0, \delta_0  ]$ and $r \in[0,2\delta R)$. It is easy to see that (\ref{ee}) implies (\ref{k1}), and
\begin{equation}
\sup\limits_{P_{\frac{3}{2}\delta_0 R }(x_0,t_0)}(| F_A+[\phi,\phi^{*{H_0}}]|_g^2+2| \partial_A \phi |_g^2)\leq f(r_0)(\delta_0R/2)^{-4}\leq 256(\delta_0R)^{-4}.
\end{equation}
In $P_{\frac{3}{2}\delta_0 R}(x_0,t_0)$, using the Bochner type inequality (\ref{12}), the parabolic mean value inequality and the Lemma 2.5., we have
\begin{equation}
\begin{split}
\sup\limits_{P_{\delta R}(x_0,t_0)}(|\nabla_A\phi|_g^2+1)&\leq C_{13} \int_{P_{3\delta_0 R/2}(x_0,t_0)}(|\nabla_A\phi|_g^2+1) dV_g\,dt
\\ &\leq C_{14}
\end{split}
\end{equation}
for any $\delta \in (0, \delta_{0})$, where $C_{13}$ and $C_{14}$ are positive constants depending only on the geometry of $(X,\omega)$, $R$, $\delta_0$ and the initial data $(A_0,\phi_0)$. This complete the proof.
\end{proof}

\medskip

             Using the above $\epsilon$-regularity theorem, and following the argument of Hong and Tian \cite{HT} for the Yang-Mills flow case, we can analyze the limiting behavior of the  heat flow (\ref{YMHH}). We have the following theorem.

\medskip

{\bf Theorem 2.7.} {\it Let $(A(t), \phi (t) )$ be a  smooth
solution of the gradient heat flow (\ref{YMHH}) on a Hermitian vector bundle $(E, H_{0})$ over K\"ahler manifold $(M
, \omega )$ with smooth initial data. Then there exists a sequence
$\{ t_{i}\}$ such that, as $t_{i}\rightarrow \infty$, $(A(t_{i}), \phi (t_{i})
)$ converges, modulo gauge transformations, to a
solution $(A_{\infty}, \phi _{\infty } )$ of the Yang-Mills-Higgs
equation (\ref{YMHE}) on Hermitian vector bundle $(E_{\infty} , H_{\infty})$  in $C^{\infty}_{loc}$ topology outside
$\Sigma^{an}\subset M$, where $\Sigma^{an}$ is a closed set of Hausdorff codimension at least $4$, and there exists an isometry between $(E, H_{0})$ and $(E_{\infty}, H_{\infty})$ outsides $\Sigma_{an}$.}

\medskip

{\bf Corollary 2.8.} {\it Let $(A(t_{i}) , \phi (t_{i}))$ be a sequence
of Higgs pairs along the gradient heat flow (\ref{YMHH}) with Uhlenbeck
limit $(A_{\infty } , \phi_{\infty})$, then:

(1) $\theta (A(t_{i}) , \phi (t_{i}))\rightarrow \theta (A_{\infty} ,
\phi_{\infty})$ strongly in $L^{p}$ for all $1\leq p<\infty$,  and consequently,
$lim_{t\rightarrow \infty }\int_{M}|\theta (A_{t} ,
\phi_{t})|^{2}=\int_{M}|\theta (A_{\infty} , \phi_{\infty})|^{2}$;

(2) $\|\theta (A_{\infty} , \phi_{\infty})\|_{L^{\infty}}\leq
\|\theta (A(t_{j}) , \phi(t_{j}))\|_{L^{\infty}}\leq \|\theta (A_{t_{0}} ,
\phi_{t_{0}})\|_{L^{\infty}}$ for $0\leq t_{0}\leq t_{j}$.

}

\medskip

{\bf Remark 2.9. } {\it Since we already have the  $\epsilon$-regularity theorem (Theorem 2.6.) for the higher dimensional case, the above theorem and corollary can be proved in the same way as that in the K\"ahler surface case (Theorem 3.2. and Corollary 3.12. in \cite{LZ1}). So, we omit the proof. }

\medskip

From the equation (\ref{YMHE}), we know that
\begin{eqnarray}
D_{A_{\infty}}\theta_{\infty} =0, \quad
[\theta_{\infty} , \phi_{\infty} ]=0,
\end{eqnarray}
where $\theta _{\infty} = \Lambda_{\omega }(F_{A_{\infty}}+[\phi_{\infty} , \phi_{\infty}^{\ast }]) $.
Since $\theta_{\infty} $ is parallel and $(\sqrt{-1}\theta_{\infty})^{\ast
}=\sqrt{-1}\theta_{\infty} $, we can decompose $E_{\infty}$ according to the
eigenvalues of $\sqrt{-1}\theta_{\infty}$. We obtain a holomorphic
orthogonal decomposition
\begin{eqnarray}
E_{\infty}=\bigoplus_{i=1}^{l}E_{\infty}^{i},
\end{eqnarray}
and
\begin{eqnarray}
\phi_{\infty} : E_{\infty}^{i} \rightarrow E_{\infty}^{i}
\end{eqnarray}
on $M\setminus \Sigma_{an}$. Let $H_{\infty}^{i}$ be the restrictions of $H_{\infty}$ to $E_{\infty}^{i}$, $\phi_{\infty}^{i}$ be
the restriction of $\phi_{\infty} $ to $E^{i}$, and $A_{\infty}^{i}=A_{\infty}|_{E^{i}}$.
Then $(A_{\infty}^{i} , \phi^{i})$ is a Higgs pair on $(E_{\infty}^{i} , H_{\infty}^{i})$ and
satisfies
\begin{eqnarray}
\sqrt{-1}\Lambda_{\omega }(F_{A_{\infty}^{i}}+[\phi_{\infty}^{i} , (\phi_{\infty} ^{i})^{\ast
}])=2\pi \lambda_{i}Id_{E_{\infty}^{i}}.
\end{eqnarray}
So $(A_{\infty}^{i} , \phi_{\infty}^{i})$ is a Hermitian-Einstein Higgs pair on
$(E_{\infty}^{i} , H_{\infty}^{i})$, i.e. $(E_{\infty}^{i} , H_{\infty}^{i} , A_{\infty}^{i} , \phi_{\infty}^{i})$ is a
Hermitian-Einstein Higgs bundle on $M\setminus \Sigma_{an}$.

The Yang-Mills-Higgs functional is decreasing along the gradient flow (\ref{YMHH}), and $\phi (t)$ is uniformly $C^{0}$ bounded (by Lemma 2.3.), then we have
\begin{eqnarray}
\int_{M\setminus \Sigma_{an}} |F_{A_{\infty}}|^{2}_{H_{\infty}} \frac{\omega^{n}}{n!}\leq C <\infty .
\end{eqnarray}
  Since the singularity set $\Sigma_{an}$ is of Hausdorff codimension $4$,  $ \phi_{\infty}$ is holomorphic and $C^{0}$ bounded,  and every metric $H_{\infty}^{i}$ (or the connection $A_{\infty}^{i}$) satisfies the Hermitian-Einstein equation (\ref{HEE}), a similar argument as that in the proof of Theorem 2 in Bando and Siu 's paper \cite{BS} can show that, every $(E_{\infty}^{i}, \overline{\partial }_{A_{\infty}^{i}})$ can be extended to the whole $M$ as a reflexive sheaf (which is also denoted by $(E_{\infty}^{i}, \overline{\partial }_{A_{\infty}^{i}})$ for simplicity), $\phi_{\infty}^{i}$ and $H_{\infty}^{i}$ can be smoothly extended over the place where the sheaf $(E_{\infty}^{i}, \overline{\partial }_{A_{\infty}^{i}})$ is locally free. So we have the following proposition.

  \medskip

  {\bf Proposition 2.10. } {\it The limiting $(E_{\infty}, \overline{\partial }_{A_{\infty}}, \phi_{\infty})$ can be extended to the whole $M$ as a reflexive Higgs sheaf  with a holomorphic orthogonal splitting
  \begin{eqnarray}
(E_{\infty}, H_{\infty}, A_{\infty},
\phi_{\infty})=\bigoplus_{i=1}^{l}(E_{\infty}^{i} , H_{\infty}^{i}
, A_{\infty}^{i} , \phi_{\infty} ^{i}),
\end{eqnarray}
where $H_{\infty}^{i}$ is an admissible Hermitian-Einstein metrics on the reflexive
Higgs sheaf $(E_{\infty}^{i} , A_{\infty}^{i} ,
\phi_{\infty}^{i})$.
  }

\medskip

Let $S$ be a $\phi$-invariant torsion free sub-sheaf of $(E , \overline{\partial}_{E},  \phi )$ with a Hermitian metric $H$. Since we can view $S$ as a holomorphic vector sub-bundle off the singular set $\Sigma$ where $S$ fails to be locally free, away from $\Sigma$ we have a corresponding orthogonal projection
$\pi : E\rightarrow E$ with $\pi (E)=S$. Since $S$ is $\phi$-invariant and holomorphic, on  almost everywhere of $M$, we have  $\pi^{2}=\pi=\pi^{\ast H}$;
$(Id -\pi) \overline{\partial }\pi =0$; and $(Id -\pi) [\phi , \pi]=0$. Furthermore, it can be shown that $\pi$ extends to an $L_{1}^{2}$ section of $End{E}$.
Conversely such $\pi$ determines a coherent $\phi$-invariant subsheaf.

\medskip

{\bf Definition 2.11. } {\it An element $\pi $ is called a weakly $\phi $ holomorphic subbundle if $\pi \in L_{1}^{2}(End(E))$ and
\begin{eqnarray}\label{WHC}
\begin{array}{lll}
&&(Id -\pi) \overline{\partial }\pi =0;\\
&& \pi^{2}=\pi=\pi^{\ast H};\\
&& (Id -\pi) [\phi , \pi]=0\\
\end{array}
\end{eqnarray}
hold almost everywhere.}

\medskip

In \cite{UY}, Uhlenbeck and Yau prove that the above $\pi$ determines a coherent subsheaf $S$, the last term in conditions (\ref{WHC}) implies that $S$ is $\phi$-invariant. So, we have the following proposition.

\medskip

{\bf Proposition 2.12. } {\it A weakly $\phi $ holomorphic subbundle $\pi $ determines a $\phi$-invariant coherent subsheaf of the Higgs bundle $(E , \overline{\partial}_{E},  \phi )$.  }

\section{Resolution of the HNS filtration }
\setcounter{equation}{0}

 Given a Higgs bundle $(E, A, \phi)$ on a
K\"ahler manifold $(M, \omega)$. A Higgs sub-sheaf of $(E , A , \phi
)$ is a $\phi
$-invariant coherent analytic  subsheaf $S\subset (E , A)$ .  The $\omega$-slope $\mu
(S)$ of a torsion-free sheaf $S\rightarrow M$ is defined by:
\begin{eqnarray}
\mu _{\omega }(S)=\frac{deg_{\omega
}(S)}{rank(S)}=\frac{1}{rank(S)}\int_{M}C_{1}(S)\wedge \frac{\omega
^{n-1} }{(n-1)!}.
\end{eqnarray}

For any subsheaf $S$, its singular set $\Sigma_{S}$ is the set of points where  $S$ fails to be locally free.   If $S$ is a saturated subsheaf then the singular set $\Sigma_{S}$ is a closed complex analytic subset of $M$ of complex codimension at least 2.
A torsion-free Higgs subsheaf $S$ is said $\omega$-stable
(resp. $\omega$-semistable) if for all proper $\phi $-invariant
saturated subsheaves $F\subset S$, $\mu_{\omega
}(F)<\mu_{\omega}(S)$ ($\mu_{\omega }(F)\leq \mu_{\omega}(S)$).

In the following, we will give a description of the appropriate
Higgs bundle versions of the Harder-Narasimhan filtration and the
Harder-Narasimhan-Seshadri filtration, the proof is almost the
same as the one used in the holomorphic bundles case (\cite{Ko2},
v.7.15, 7.17, 7.18; \cite{HT}, section 7 ), the only difference is that we always
consider $\phi$-invariant subsheaves instead of general
subsheaves. We omit the details here.

\medskip

{\bf Proposition 3.1. }  {\it Let $(E, A, \phi )$ be a Higgs bundle on K\"ahler manifold $(M ,
\omega )$. Then there is a filtration of $E$ by $\phi$-invariant coherent
 sub-sheaves
$$
0=E_{0}\subset E_{1}\subset \cdots \subset E_{l}=E ,$$ called the
Harder-Narasimhan filtration of Higgs bundle $(E, A, \phi )$ (abbr,
HN-filtration ), such that $Q_{i}=E_{i}/E_{i-1}$ is torsion-free and
Higgs semistable. Moreover, $\mu (Q_{i})>\mu (Q_{i+1})$, and the
associated graded object $Gr^{hn}(E, A, \phi
)=\oplus_{i=1}^{l}Q_{i}$ is uniquely determined by the isomorphism
class of $(E, A, \phi)$.}

\hspace{0.3cm}

{\bf Proposition 3.2. } {\it Let $(V , \phi )$ be a semistable Higgs
sheaf on K\"ahler manifold $(M , \omega )$, then there is a
filtration of $V$ by $\phi$-invariant subsheaf $$ 0=V_{0}\subset
V_{1}\subset \cdots \subset V_{l}=V ,$$ called the Seshadri
filtration of  $(V, \phi )$, such that $V_{i}/V_{i-1}$ is
torsion-free and Higgs stable. Moreover, $\mu (V_{i}/V_{i-1})=\mu
(V)$ for each $i$, and the associated graded object $Gr^{s}(V, \phi
)=\oplus_{i=1}^{l}V_{i}/V_{i-1}$ is uniquely determined by the
isomorphism class of $(V, \phi)$.}

\medskip

Combining the two previous proposition, fora given Higgs bundle $(E, A, \phi )$,  there exists a double filtration , called the
Harder-Narasimhan-Seshadri filtration of  $(E, A, \phi
)$, which we abbreviate by HNS-filtration.

\medskip

{\bf Proposition 3.3. } {Let $(E, A, \phi )$ be a Higgs bundle on K\"ahler manifold $(M ,
\omega )$. Then there is a double filtration $\{E_{i, j}\}$ with the following properties: if
$\{E_{i}\}_{i=1}^{l}$ is the HN filtration of $(E, A, \phi )$, then
$$
E_{i-1}=E_{i , 0}\subset E_{i , 1}\subset \cdots \subset E_{i ,
l_{i}}=E_{i}
$$
and the successive quotient $Q_{i , j}=E_{i , j}/E_{i , j-1}$ are
Higgs stable torsion-free sheaves. Moreover, $\mu (Q_{i , j})=\mu
(Q_{i , j+1})$ and $\mu (Q_{i , j})>\mu (Q_{i+1 , j})$, the
associated graded object: $$ Gr^{hns}(E, A, \phi
)=\oplus_{i=1}^{l}\oplus_{j=1}^{l_{i}}Q_{i , j}
$$is uniquely determined by the isomorphism class of $(E, A, \phi )$.}

\hspace{0.3cm}

The following lemma can be proved by an argument similar  the one used in
Chapter 5, V.7.11; 7.12 in \cite{Ko2} for the case of holomorphic bundles.

\hspace{0.3cm}

{\bf Lemma 3.4. (Lemma 6.3. in \cite{LZ1}) } {\it Let $(E_{1}, \overline{\partial}_{A_{1}} , \phi_{1})$ and
$(E_{2} , \overline{\partial}_{A_{2}}, \phi_{2})$ be  two semistable Higgs sheaves with
same rank and degree. If $(E_{1}, \overline{\partial}_{A_{1}} , \phi_{1})$ is Higgs stable, and let  $f: E_{1}\rightarrow
E_{2}$ be a  sheaf homomorphism satisfying $f\circ
\phi_{1}=\phi_{2}\circ f$, then either $f=0$ or $f$
 is injective.}

\medskip

{\bf Definition 3.5. (Harder-Narasimhan type ) } {\it For a Higgs
bundle $(E, A, \phi )$ of rank $R$, construct a nonincreasing
$R$-tuple of numbers
\begin{eqnarray}
\vec{\mu }(E, A, \phi )=(\mu_{1} , \cdots , \mu_{R})
\end{eqnarray}
from the HN filtration by setting: $\mu_{i}=\mu (Q_{j})$, for
$rk(E_{j-1})+1\leq i\leq rk(E_{j})$. We call $\vec{\mu }(E, A,
\phi )$ the Harder-Narasimhan type of $(E, A, \phi )$.}

\hspace{0.3cm}

{\bf Remark: } {\it For a pair $\vec{\mu }$, $\vec{\lambda  }$ of
$R$-tuple's satisfying
$\sum_{i=1}^{R}\mu_{i}=\sum_{i=1}^{R}\lambda_{i}$, we define:
\begin{eqnarray}
\vec{\mu }\leq \vec{\lambda  } \quad \Leftrightarrow \quad
\sum_{i\leq k}\mu_{i}\leq \sum_{i\leq k}\lambda_{i}
\end{eqnarray}
for all $k=1, \cdots , R$.

 }

\medskip

It will be convenient to denote the $\phi $-invariant subsheaf
$E_{i}$ in the HN-filtration by $\texttt{F}_{i}^{hn}(E , A, \phi )$
 or by $\texttt{F}_{i , \omega }^{hn}(E , A, \phi )$  when we wish
to emphasize the role of the K\"ahler structure. Let $\{E_{i,j}\}$ be the HNS-filtration of the Higgs bundle $(E, A, \phi )$, we set
\begin{eqnarray}
\Sigma_{alg}=\cup_{i,j}Sing(E_{i,j})\cup Sing(Q_{i,j}),
\end{eqnarray}
this is a complex analytic subset of complex codimension at least two. We will refer to it as the singular set of the HNS-filtration.
Since the HNS-filtration fails to be given by subbundles on the singular set $\Sigma_{alg}$, it makes difficult to do analysis.  We can use the singularities theorem of Hironaka (\cite{Hio1}) to resolve the singularities $\Sigma_{alg}$, and obtain a filtration by subbundles. This idea had been used by Bando and Siu \cite{BS} to obtain admissible Hermitian-Einstein metric on reflexive stable sheaf. The following proposition has been proved by Sibley in \cite{Sib}.

\medskip

{\bf Proposition 3.6. (Proposition 4.3. in \cite{Sib} )} {\it Let $0=E_{0}\subset E_{1}\subset \cdots \subset E_{l-1}\subset E_{l}=E$ be a filtration of a holomorphic vector bundle $E$ on a complex manifold $M$ by saturated subsheaves and let $Q_{i}=E_{i}/E_{i-1}$. Then there is a finite sequence of blowups along complex submanifolds of $M$ whose composition $\pi : \tilde{M} \rightarrow M$ enjoys the following properties. There is a filtration
\begin{eqnarray*}
0=\tilde{E}_{0}\subset \tilde{E}_{1}\subset \cdots \subset \tilde{E}_{l-1}\subset \tilde{E}_{l}=\tilde{E}
\end{eqnarray*}
by subbundles such that $\tilde{E}_{i}$ is the saturation of $\pi^{\ast } E_{i}$. If $\tilde{Q}_{i}=\tilde{E}_{i}/\tilde{E}_{i-1}$, then we have exact sequences:
\begin{eqnarray*}
0\rightarrow E_{i} \rightarrow \pi_{\ast }\tilde{E}_{i}\rightarrow T_{i} \rightarrow 0
\end{eqnarray*}
and
\begin{eqnarray*}
0\rightarrow Q_{i} \rightarrow \pi_{\ast }\tilde{Q}_{i}\rightarrow T_{i}' \rightarrow 0
\end{eqnarray*}
where $T_{i}$ and $T_{i}'$ are torsion sheaves supported on the singular sets of $E_{i}$ and $Q_{i}$ respectively, and furthermore $\pi_{\ast }\tilde{E}_{i}=E_{i}$ and $Q_{i}^{\ast \ast }=(\pi_{\ast }\tilde{Q}_{i})^{\ast \ast}$.
}

\medskip

Let $\phi $ be a Higgs field on holomorphic bundle $(E, \overline{\partial}_{A})$ and $\tilde{\phi }=\pi^{\ast }\phi $
 be the pullback Higgs field on $\tilde{E}$. If the filtration $\{E_{i}\}_{i=1}^{l}$ is by $\phi $-invariant subsheaves, then the pullback filtration $\{\tilde{E}_{i}\}_{i=1}^{l}$ in the above proposition is by $\tilde{\phi }$-invariant subbundles. So, we have the following proposition.

\medskip

{\bf Proposition 3.7. } {\it Let $\{E_{i,j}\}$ be the HNS-filtration of a Higgs bundle $(E, A, \phi )$ on complex manifold $M$ and let $Q_{i, j}=E_{i, j}/E_{i, j-1}$. Then there is a finite sequence of blowups along complex submanifolds of $M$ whose composition $\pi : \tilde{M} \rightarrow M$ enjoys the following properties. There is a filtration  $\{\tilde{E}_{i,j}\}$
by $\tilde{\phi }$-subbundles such that $\tilde{E}_{i,j}$ is the saturation of $\pi^{\ast } E_{i,j}$, and $\pi_{\ast }\tilde{E}_{i, j}=E_{i, j}$ and $Q_{i, j}^{\ast \ast }=(\pi_{\ast }\tilde{Q}_{i, j})^{\ast \ast}$, where $\tilde{\phi }=\pi^{\ast }\phi $.
}

 \medskip

 The following proposition is well known, the proof can be found in, for example \cite{GH}.

 \medskip

 {\bf Proposition 3.8. } {\it Let $(M, \omega )$ be a compact K\"ahler manifold and $\pi : \tilde{M} \rightarrow M$ be a blow up with non-singular center. Then $\tilde{M}$ is also K\"ahler, moreover $\tilde{M}$ possesses a one family of K\"ahler metrics given by $\omega_{\epsilon }=\pi^{\ast}\omega +\epsilon \eta $ where $\epsilon >0$ and $\eta $ is K\"ahler metric on $\tilde{M}$.}

 \medskip

 Furthermore, Bando and Siu (\cite{BS}) proved the uniform boundedness of heat kernel for $\omega_{\epsilon}$.

 \medskip

 {\bf Proposition 3.9. (Proposition 2 in \cite{BS})} {\it Let $(M, \omega )$ be an $n$-dimensional compact K\"ahler manifold and $\pi : \tilde{M} \rightarrow M$ be a blow up with non-singular center of co-dimensional at least two. Fix a K\"ahler metric $\eta $ on $\tilde{M}$ and set $\omega_{\epsilon }=\pi^{\ast}\omega +\epsilon \eta $ for $0<\epsilon \leq 1$. Let $K_{\epsilon}(t, x, y)$ be the heat kernel with respect to the metric $\omega_{\epsilon}$, then we have a uniform estimate $0\leq K_{\epsilon}(t, \cdot , \cdot)\leq C (1+t^{-n})$ with a positive constant $C$.}

\medskip

In the following, we consider the $\omega_{\epsilon}$ slope of an arbitrary coherent subsheaf of a holomorphic vector bundle $\tilde{E}$ on the blow up $\tilde{M}$. One can show that the  $\omega_{\epsilon}$ slope converges to the $\omega $ slope of the push forward sheaf on the base $M$, and the stability will be preserved for small $\epsilon $. These properties should be well known, see Bando and Siu's paper \cite{BS}, more details can be found in Sibley's paper \cite{Sib}.

\medskip

{\bf Proposition 3.10.} {\it Let $\tilde{S}$ be a coherent subsheaf (with torsion free quotient $\tilde{Q}$) of a holomorphic vector bundle $\tilde{E}$ on $\tilde{M}$, and $\pi : \tilde{M} \rightarrow M$ be a blow up with non-singular center of co-dimensional at least two. Then $\mu _{\omega_{\epsilon }}(\tilde{S}) \rightarrow \mu_{\omega }(\pi_{\ast } \tilde{S})$ and $\mu _{\omega_{\epsilon }}(\tilde{Q}) \rightarrow \mu_{\omega }(\pi_{\ast } \tilde{Q})$ as $\epsilon \rightarrow 0$. Furthermore there is a uniform constant $B$ independent of $S$ such that $\mu _{\omega_{\epsilon }}(\tilde{S}) \leq  \mu_{\omega }(\pi_{\ast } \tilde{S})+\epsilon B $ and $\mu _{\omega_{\epsilon }}(\tilde{Q}) \geq  \mu_{\omega }(\pi_{\ast } \tilde{Q})-\epsilon B $. }

\medskip

{\bf Proof. }
By the definition, we have
\begin{eqnarray}\label{3101}
\begin{array}{lll}
deg (\tilde{S}, \omega_{\epsilon })&=&\int_{\tilde{M}}c_{1} (det (\tilde{S})) \wedge \frac{\omega_{\epsilon}^{n-1}}{(n-1)!}\\
&=&\int_{\tilde{M}}c_{1} (det (\tilde{S})) \wedge \frac{(\pi^{\ast}\omega)^{n-1}}{(n-1)!}\\
&&+\frac{1}{(n-1)!}\sum_{i=1}^{n-1}\epsilon ^{i}C_{n}^{i}\int_{\tilde{M}}c_{1} (det (\tilde{S})) \wedge \eta ^{i}\wedge (\pi^{\ast}\omega )^{n-1-i}.\\
\end{array}
\end{eqnarray}
Since the blow up set is co-dimensional at least two, so
\begin{eqnarray}\label{3102}
\begin{array}{lll}
&&\int_{\tilde{M}}c_{1} (det (\tilde{S})) \wedge \frac{\omega_{\epsilon}^{n-1}}{(n-1)!}\\
&=& \int_{\tilde{M}}c_{1} (\pi_{\ast}det (\tilde{S})) \wedge \frac{\omega_{\epsilon}^{n-1}}{(n-1)!}\\
&=& \int_{\tilde{M}}c_{1} (det (\pi_{\ast}\tilde{S})) \wedge \frac{\omega_{\epsilon}^{n-1}}{(n-1)!}\\
&=& deg (\pi_{\ast } S, \omega ),
\end{array}
\end{eqnarray}
where we used the isomorphism $det (\pi_{\ast}\tilde{S})=\pi_{\ast}det (\tilde{S})$. Then (\ref{3101}) and (\ref{3102}) imply $\mu _{\omega_{\epsilon }}(\tilde{S}) \rightarrow \mu_{\omega }(\pi_{\ast } \tilde{S})$ as $\epsilon \rightarrow 0$.

Let $H$ be a Hermitian metric on $\tilde{E}$,  we can view $\tilde{S}$ as a holomorphic sub-bundle off the singular set $\Sigma$ which is co-dimensional at least two, away from $\Sigma$ we has a corresponding orthogonal projection
$\pi_{\tilde{S}} : \tilde{E}\rightarrow \tilde{E}$ with $\pi_{\tilde{S}} (E)=\tilde{S}$. It is well known that $\pi_{\tilde{S}}$ can be extended to an
$L_{1}^{2}$ section of $End{E}$. Using the Gauss-Codazzi equation, we have
\begin{eqnarray}\label{3103}
\begin{array}{lll}
&&\int_{\tilde{M}}c_{1} (det (\tilde{S})) \wedge \eta ^{i}\wedge (\pi^{\ast}\omega )^{n-1-i}\\
&=&\frac{\sqrt{-1}}{2\pi}\int_{\tilde{M}}tr(\pi_{\tilde{S}}\circ F_{A_{H}}\circ \pi_{\tilde{S}} +\overline{\partial }_{\tilde{E}}\pi_{\tilde{S}} \wedge \partial _{A_{H}}\pi_{\tilde{S}}) \wedge \eta ^{i}\wedge (\pi^{\ast}\omega )^{n-1-i}\\
\end{array}
\end{eqnarray}
where $A_{H}$ is the Chern connection with respect to $H$. Since $\pi^{\ast}\omega $ is nonnegative,  the second term in right hand side of the above equality is non-positive. Since the first term is uniform bound, we see that $\int_{\tilde{M}}c_{1} (det (\tilde{S})) \wedge \eta ^{i}\wedge (\pi^{\ast}\omega )^{n-1-i}$ has an uniform upper bound  independent of $\tilde{S}$. By (\ref{3101}), (\ref{3102}) and (\ref{3103}), we know there is a constant $B$ such that $\mu _{\omega_{\epsilon }}(\tilde{S}) \leq  \mu_{\omega }(\pi_{\ast } \tilde{S})+\epsilon B $. For the quotient sheaf $\tilde{Q}$, we can consider the dualised sequence $0\rightarrow \tilde{Q}^{\ast}\rightarrow \tilde{E}^{\ast } \rightarrow \tilde{S}^{\ast} \rightarrow 0$, a similar argument as above implies the statements for $Q$ in the proposition.

\hfill $\Box$ \\

\medskip

{\bf Remark. 3.11. } {\it  If there is a sequence of blow-ups:
$$\pi_{i}: \overline{M}_{i}\rightarrow \overline{M}_{i-1}, \quad i=1, \cdots , r
$$
and  $\pi =\pi_{r}\circ \cdots \circ \pi_{1} $, where $\overline{M}_{0}=M$ and $\overline{M}_{r}=\tilde{M}$ and every $\pi_{i}$ is blow up along a smooth complex submanifold of co-dimensional at least two.  On each blow-up $\overline{M}_{i}$, we have a family of K\"ahler metrics defined iteratively by $\omega_{\epsilon_{1}\cdots \epsilon_{i}}=\pi_{i}^{\ast}\omega_{\epsilon_{1}\cdots \epsilon_{i-1}}+\epsilon_{i}\eta_{i}$, where $\eta_{i}$ is a K\"ahler metric on $\overline{M}_{i}$ and $\epsilon_{i} >0$. For simplicity, in the following, we will denote $\epsilon =(\epsilon_{1}, \cdots , \epsilon _{r})$,  $\omega_{\epsilon }=\omega_{\epsilon_{1}\cdots \epsilon_{r}}$, and $\|\epsilon \|=\max_{i}\epsilon_{i}$. It is easy to see that the above proposition is also valid for such $\pi$. }

\medskip

{\bf Proposition 3.12.} {\it Let $\pi: \tilde{M} \rightarrow M$ be a composition of finitely many blowups along complex submanifolds of co-dimensional at least two, $(\tilde{E} ,
\tilde{\phi })$ be a Higgs bundle
 over $\tilde{M}$,
and $(E , \phi )$ be a Higgs sheaf  over $M$ with
$\pi_{\ast }\tilde{E}=E$ ,  $\phi (\pi_{\ast
}X)=\pi_{\ast }\tilde{\phi}(X)$ for any $X\in
\tilde{E}$. If the Higgs sheaf $(E , \phi )$
is $\omega $-stable, then there is a number $\epsilon_{0}>0$,
such that the Higgs sheaf $(\tilde{E} ,
\tilde{\phi })$ is $\omega_{\epsilon }$-stable for all
$0<\epsilon \leq \epsilon_{0}$.}

\medskip

{\bf Proof. } By the assumption that E is $\omega$-stable, then there is constant $\delta >0$ such that $\mu_{\omega } (E)-\mu_{\omega } (S)>\delta $ for any proper Higgs subsheaf $S$. By proposition 3.10, for any proper Higgs subsheaf $\tilde{S}\subset \tilde{E}$ we have
\begin{eqnarray}
\begin{array}{lll}
&&\mu_{\omega_{\epsilon }} (\tilde{S})-\mu_{\omega_{\epsilon }} (\tilde{E})\\
&\leq & \mu_{\omega } (\pi_{\ast} (\tilde{S})) -\mu_{\omega } (E) +2\|\epsilon \| B \\
&<& -\delta +2\|\epsilon \| B <0
\end{array}
\end{eqnarray}
for $\|\epsilon\| <\frac{\delta }{2B}$, where we used that $B$ is independent on $\tilde{S}$. This completes the proof.

\hfill $\Box$ \\

\medskip

Let $(S, \phi )$ be a Higgs sheaf on $M$, we  define $\mu_{\max , \omega }(S)$ to be the maximum $\omega $-slope of $\phi $-invariant subsheaves of $S$, and
$\mu_{\min , \omega }(S)$ to be the minimal $\omega $-slope of $\phi $-invariant torsion free quotient sheaves of $S$. It is easy to check that $\mu_{\min , \omega }(S)=-\mu_{\max , \omega }(S^{\ast})$. By Proposition 3.10, we have the following corollary.

\medskip

{\bf Corollary 3.13.} {\it Let $\pi: \tilde{M} \rightarrow M$ be a composition of finitely many blowups along complex submanifolds of co-dimensional at least two, $(\tilde{E} ,
\tilde{\phi })$ be a Higgs bundle
 over $\tilde{M}$,
and $(E , \phi )$ be a Higgs sheaf  over $M$ with
$\pi_{\ast }\tilde{E}=E$ ,  $\phi (\pi_{\ast
}X)=\pi_{\ast }\tilde{\phi}(X)$ for any $X\in
\tilde{E}$. There is a constant $B>0$ such that:

(1): $\mu_{\max , \omega_{\epsilon } }(\tilde{E})\leq \mu_{\max , \omega }(E) +\|\epsilon \|B$;

(2): $\mu_{\min , \omega_{\epsilon } }(\tilde{E})\geq \mu_{\min , \omega }(E) -\|\epsilon \|B$.
}

\medskip

Let the filtration
\begin{eqnarray*}
0=\tilde{E}_{0}\subset \tilde{E}_{1}\subset \cdots \subset \tilde{E}_{l-1}\subset \tilde{E}_{l}=\tilde{E}=\pi^{\ast}E
\end{eqnarray*}
be the resolution of the HN filtration of the Higgs bundle $(E, \phi )$. By proposition, we have
 \begin{eqnarray}
 \mu_{\omega_{\epsilon}} (\tilde{Q}_{i})\rightarrow \mu_{\omega}(Q_{i})
 \end{eqnarray}
for all $i$ as $\epsilon \rightarrow 0$.
 Using the properties $\mu_{\omega }(Q_{i})>\mu_{\omega}(Q_{i+1})$, $\mu_{\min ,\omega}(E_{i})=\mu_{\omega }(Q_{i})$, $\mu_{\max , \omega }(E/E_{i})=\mu_{\omega }(Q_{i+1})$, and Corollary 3.13., we have
\begin{eqnarray}
\mu_{\min , \omega_{\epsilon }}(\tilde{E}_{i})>\mu_{\max , \omega_{\epsilon }}(\tilde{E}/\tilde{E}_{i}).
\end{eqnarray}
By the above inequality, it is easy to see that the resolution appears within the HN filtration of the Higgs bundle $(\tilde{E}, \tilde{\phi })$ with respect to $\omega_{\epsilon}$, and two successive  Higgs bundles in the resolution are seperated by the HN filtration of the larger Higgs bundle. By an inductive argument, repeatedly using proposition 3.10 one can show the convergence of the HN type, so we have the following proposition.

\medskip

{\bf Proposition 3.14. }  {\it Let $(E, \phi )$ be a Higgs bundle
 over $M$, $\pi : \tilde{M} \rightarrow M$ be a sequence of blowups resolving the HNS filtration
and $(\tilde{E} ,
\overline{\phi })=\pi^{\ast }(M, \phi )$ be the pull back Higgs bundle over $\tilde{M}$. Let $\vec{\mu }_{\epsilon }$ denote the HN
type of $(\tilde{E} ,
\overline{\phi })$ with respect
to $\omega_{\epsilon }$ and $\vec{\mu}$ the HN type of
$(E , \phi )$ with respect to $\omega $, then
$\vec{\mu}_{\epsilon }\rightarrow \vec{\mu }$ as $\epsilon
\rightarrow 0$.}

\medskip

\medskip

\section{Existence of non-zero holomorphic map }
\setcounter{equation}{0}

In order to prove theorem 1.1, we need construct non-zero holomorphic maps from subsheaves in the HNS filtration of the original Higgs bundle to the limiting  reflexive sheaf. For bundle case, we can follow Donaldson's argument in \cite{D0}  to construct non-zero holomorphic maps. But in general the HNS filtration is given by subsheaves, so Donaldson's argument can not be applied directly in our case.
The following proposition is the key to construct an isomorphism between $(E_{\infty}, \overline{\partial}_{A_{\infty}}, \phi_{\infty})$ and the double dual of the stable quotients of the HNS filtration $Gr^{HNS}(E, \overline{\partial}_{A_{0}}, \phi_{0})$.

\medskip

{\bf Proposition 4.1. } {\it Let $(M, \omega )$ be a K\"ahler manifold,  $(E, A_{0}, \phi_{0})$ be a Higgs sheaf on $M$ with Hermitian metric $H_{0}$, $S$ be a Higgs sub-sheaf of $(E, A_{0}, \phi_{0})$, and $(A_{j}, \phi_{j})=g_{j}(A_{0}, \phi_{0})$ be a sequence of Higgs pairs on $E$, where $g_{j}$ is a sequence of complex gauge transformations.  Suppose that there exits a sequence of blow-ups:
$
\pi_{i}: \overline{M}_{i}\rightarrow \overline{M}_{i-1},$ $ i=1, \cdots , r
$
(where $\overline{M}_{0}=M$,  every $\pi_{i}$ is a blow up with non-singular center;  denoting $\pi =\pi_{r}\circ \cdots \circ \pi_{1}$ ); such that $\pi^{\ast }E$ and $\pi^{\ast }S$ are bundles, the pulling back geometric objects $\pi^{\ast}(A_{0}, \phi_{0})$, $\pi^{\ast }g_{j}$ and $\pi^{\ast }H_{0}$ can be extended smoothly on the whole $M_{r}$.  Assume that $(A_{j}, \phi_{j})$ converges  to $(A_{\infty }, \phi_{\infty })$ outside a closed subset $\Sigma_{An}$ of Hausdorff complex codimension 2, and $|\Lambda _{\omega}(F_{A_{j}})|_{H_{0}}$ is bounded uniformly in $j$ in $L^{1}(\omega_{0})$. Let $i_{0}: (S, \bar{\partial}_{A_{0}})\rightarrow (E, \bar{\partial}_{A_{0}})$ be the holomorphic inclusion, then there is a subsequence of $g_{j} \circ i_{0}$, up to rescale,  converges to a non-zero holomorphic map $f_{\infty }: (S, \overline{\partial}_{A_{0}})\rightarrow (E_{\infty }, \overline{\partial}_{A_{\infty}})$ in $C^{\infty}_{loc}$ off $\Sigma \cup \Sigma_{An}$, and $f_{\infty} \circ \phi_{0}=\phi_{\infty} \circ f_{\infty}$, where $\Sigma$ is the singular set of $S$ and $E$.}

\medskip

{\bf Proof. }
 On each blow-up $\overline{M}_{i}$, we have a family of K\"ahler metrics defined iteratively by $\omega_{\epsilon_{1}\cdots \epsilon_{i}}=\pi_{i}^{\ast}\omega_{\epsilon_{1}\cdots \epsilon_{i-1}}+\epsilon_{i}\eta_{i}$, where $\eta_{i}$ is a K\"ahler metric on $\overline{M}_{i}$.
For simplicity, we write $\omega_{\epsilon} =\omega_{\epsilon _{1}, \cdots , \epsilon _{r}}$, $\tilde{E}=\pi ^{\ast }E$. In the following, we denote geometric objects and their pulling back by the same notation, and $\tilde{H}_{0}=\pi^{\ast }H_{0}$.

Define the map $\tilde{\eta}_{j}: (\tilde{S}, \bar{\partial}_{A_{0}})\rightarrow (\tilde{E}, \bar{\partial}_{A_{j}})$ by $\tilde{\eta}_{j}=g_{j} \circ i_{0}$. It is easy to check that
\begin{eqnarray}
\overline{\partial }_{A_{0} , A_{j}} \tilde{\eta }_{j}=0 , \quad \tilde{\eta}_{j}\circ
\phi_{0} =\phi_{j} \circ \tilde{\eta}_{j},
\end{eqnarray}
i.e. $\tilde{\eta }_{j}$ is a $\phi $-invariant holomorphic map. For simplicity,
we will denote $\overline{\partial }_{A_{0} , A_{j}}$ by
$\overline{\partial }_{0 , j}$, and the trace Laplacian operator on the section of $S^{\ast}\otimes E$ with respect to the connection $A_{0}\otimes A_{j}$ by $\triangle _{0,j}$.

Let $H_{j, \epsilon } (t)$ and $H^{S}_{\epsilon }(t)$ be the solutions of Donaldson's flow on holomorphic bundles $(\tilde{E}, \overline{\partial}_{A_{j}})$ and $(\tilde{S}, \overline{\partial}_{A_{0}})$ with the fixed initial metrics $\tilde{H}_{0}$ and $H^{S}_{0}$  and with respect to the metric $\omega_{\epsilon}$, i.e. it satisfies the following heat equation
\begin{eqnarray}\label{DDD1}
H^{-1}\frac{\partial H}{\partial
t}=-2\sqrt{-1}\Lambda_{\omega_{\epsilon}}F_{H}.
\end{eqnarray}

 By the heat flow, we have
\begin{eqnarray}\label{H00}
(\triangle_{\epsilon } -\frac{\partial }{\partial t}) |\Lambda _{\omega_{\epsilon }}(F_{H_{j, \epsilon}(t)}) |_{H_{j, \epsilon}(t)}\geq  0,
\end{eqnarray}
 The maximum principal implies that, for $t>0$,
 \begin{eqnarray}\label{H000}
 |\Lambda _{\omega_{\epsilon }}(F_{H_{j, \epsilon}(t)}) |_{H_{j, \epsilon}(t)}(x) \leq \int_{\tilde{M}}K_{\epsilon} (t, x, y)|\Lambda _{\omega_{\epsilon }}(F_{\tilde{A}_{j}}) |_{\tilde{H}_{0}}\frac{\omega_{\epsilon }^{n}}{n!} ,
\end{eqnarray}
where $K_{\epsilon} (t, x, y)$ is the heat kernel of the Laplacian with respect to $\omega_{\epsilon }$. On the other hand, we have
\begin{eqnarray}
\begin{array}{lll}
&&\triangle_{\epsilon } |\tilde{\eta}_{j}|_{H^{S}_{\epsilon}(t), H_{j, \epsilon } (t)}^{2}=2|\partial_{H^{S}_{\epsilon}(t), H_{j, \epsilon } (t)}\tilde{\eta}_{j}|^{2}\\
&& -2 <\sqrt{-1}\Lambda_{\omega_{\epsilon }}F_{H_{j, \epsilon } (t)}\circ \tilde{\eta}_{j} , \eta_{j}> +2<\tilde{\eta}_{j}\circ \sqrt{-1}\Lambda_{\omega_{\epsilon }}F_{H^{S}_{\epsilon}(t)} , \eta_{j}>,
\end{array}
\end{eqnarray}
\begin{eqnarray}\label{kkk1}
\begin{array}{lll}
&&\frac{\partial }{\partial t} |\tilde{\eta}_{j}|_{H^{S}_{\epsilon}(t), H_{j, \epsilon } (t)}^{2} =\\
&& -2 <\sqrt{-1}\Lambda_{\omega_{\epsilon }}F_{H_{j, \epsilon } (t)}\circ \tilde{\eta}_{j} , \eta_{j}> +2<\tilde{\eta}_{j}\circ \sqrt{-1}\Lambda_{\omega_{\epsilon }}F_{H^{S}_{\epsilon}(t)} , \eta_{j}>,
\end{array}
\end{eqnarray}
and then
\begin{eqnarray}
(\triangle_{\epsilon }-\frac{\partial }{\partial t}) |\tilde{\eta}_{j}|_{H^{S}_{\epsilon}(t), H_{j, \epsilon } (t)}^{2}\geq 0.
\end{eqnarray}
Using the Maximum principle again, we have
 \begin{eqnarray}\label{H0001}
 |\tilde{\eta}_{j}|^{2}_{H^{S}_{\epsilon}(t_{0}+t), H_{j, \epsilon } (t_{0}+t)}(x) \leq \int_{\tilde{M}}K_{\epsilon} (t, x, y)|\tilde{\eta}_{j}|^{2}_{H^{S}_{\epsilon}(t_{0}), H_{j, \epsilon } (t_{0})}\frac{\omega_{\epsilon }^{n}}{n!} ,
\end{eqnarray}
for any $t_{0}\geq 0$ and $t>0$.

By \cite{BS} (Lemma 4), for fixed $\epsilon ' =(\epsilon_{1}, \cdots , \epsilon_{r-1})$ the heat kernel $K_{\epsilon} (t, x, y) $ has a uniform bound for $0<\epsilon_{r}\leq 1$.   It is easy to see that $K_{\epsilon }$ converges to the heat kernel $K_{\epsilon '}$ on $M_{r-1}$ outside the exceptional divisor as $\epsilon_{r} \rightarrow 0$. Bando and Siu (\cite{BS}) had shown that we could choose a subsequence of $H_{j, \epsilon } (t)$ (and the same for $H^{S}_{\epsilon }(t)$) which converges to a solution of the Donaldson's heat flow (\ref{DDD1}) on $M_{r-1}$ as $\epsilon_{r}$ tends to $0$. Then we have
\begin{eqnarray}\label{H0002}
 |\Lambda _{\omega_{\epsilon' }}(F_{H_{j, \epsilon'}(t)}) |_{H_{j, \epsilon'}(t)}(x) \leq \int_{M_{r-1}}K_{\epsilon'} (t, x, y)|\Lambda _{\omega_{\epsilon' }}(F_{\tilde{A}_{j}}) |_{\tilde{H}_{0}}\frac{\omega_{\epsilon' }^{n}}{n!} ,
\end{eqnarray}
and
 \begin{eqnarray}\label{H0003}
 |\tilde{\eta}_{j}|^{2}_{H^{S}_{\epsilon'}(t_{0}+t), H_{j, \epsilon' } (t_{0}+t)}(x) \leq \int_{M_{r-1}}K_{\epsilon'} (t, x, y)|\tilde{\eta}_{j}|^{2}_{H^{S}_{\epsilon'}(t_{0}), H_{j, \epsilon' } (t_{0})}\frac{\omega_{\epsilon' }^{n}}{n!} ,
\end{eqnarray}
for all $x$ outside the exceptional set.

          Taking the limit $\epsilon_{r-1}\rightarrow 0$, and repeating the argument, we have a solution of the heat flow (\ref{DDD1}) $H_{j}(t)$ (and $H^{S}(t)$) on $M$. We also have
\begin{eqnarray}\label{H0004}
 |\Lambda _{\omega}(F_{H_{j}(t)}) |_{H_{j}(t)}(x) \leq \int_{M}K (t, x, y)|\Lambda _{\omega}(F_{A_{j}}) |_{H_{0}}\frac{\omega ^{n}}{n!} ,
\end{eqnarray}
and
 \begin{eqnarray}\label{H0005}
 |\tilde{\eta}_{j}|^{2}_{H^{S}(t_{0}+t), H_{j } (t_{0}+t)}(x) \leq \int_{M}K (t, x, y)|\tilde{\eta}_{j}|_{H^{S}(t_{0}), H_{j} (t_{0})}^{2}\frac{\omega ^{n}}{n!} ,
\end{eqnarray}
for all $x$ outside $\Sigma $, where $K(t, x, y)$ is the heat kernel of $(M, \omega )$. Using  $K(t, x, y)\leq C_{K} (1+ \frac{1}{t^{n}})$, and the uniform $L^{1}$ bound in the assumptions, we have a uniform constant $C_{F}$ which is independent on $j$ such that
\begin{eqnarray}\label{FFF1}
2( |\Lambda _{\omega}(F_{H_{j}(t)}) |_{H_{j}(t)}+  |\Lambda _{\omega}(F_{H^{S}(t)}) |_{H^{S}(t)}) (x)\leq C_{F}
\end{eqnarray}
for all $x\in M\setminus \Sigma $ and $t\geq t_{0}>0$.

By (\ref{kkk1}) and (\ref{FFF1}), we have
\begin{eqnarray}
-C_{F}\leq \frac{\partial }{\partial t} \ln  |\tilde{\eta}_{j}|^{2}_{H^{S}(t), H_{j } (t)}(x)\leq C_{F} ,
\end{eqnarray}
for all $x\in M\setminus \Sigma $ and $t\geq t_{0}>0$. Then
\begin{eqnarray}
e^{-C_{F}\delta }\leq \frac{|\tilde{\eta}_{j}|^{2}_{H^{S}(t_{0}+\delta ), H_{j } (t_{0}+\delta)}}{|\tilde{\eta}_{j}|^{2}_{H^{S}(t_{0}), H_{j } (t_{0})}}(x)\leq e^{C_{F}\delta },
\end{eqnarray}
and
\begin{eqnarray}\label{key11}
\begin{array}{lll}
&&|\tilde{\eta}_{j}|^{2}_{H^{S}(t_{0}), H_{j } (t_{0})}(x)\leq e^{C_{F}\delta } |\tilde{\eta}_{j}|^{2}_{H^{S}(t_{0}+\delta ), H_{j } (t_{0}+\delta)}(x)\\
&& \leq e^{C_{F}\delta } \int_{M}K (t, x, y)|\tilde{\eta}_{j}|_{H^{S}(t_{0}), H_{j} (t_{0})}^{2}\frac{\omega ^{n}}{n!}\\
&& \leq C_{K}e^{C_{F}\delta } (1+\delta^{-n})\int_{M}|\tilde{\eta}_{j}|_{H^{S}(t_{0}), H_{j} (t_{0})}^{2}\frac{\omega ^{n}}{n!}.\\
\end{array}
\end{eqnarray}

 Denote $h_{j, \epsilon }(t)=\tilde{H}_{0}^{-1}H_{j, \epsilon }(t)$, it is easy to check that
\begin{eqnarray}\label{H01}
\begin{array}{lll}
&& (\triangle_{\epsilon } -\frac{\partial }{\partial t}) \ln (tr (h_{j, \epsilon }(t))+ tr (h_{j, \epsilon }(t)^{-1}))\\
&\geq & -2 |\Lambda _{\omega_{\epsilon }}(F_{A_{j}}) |_{\tilde{H}_{0}}.\\
\end{array}
\end{eqnarray}
From the above inequality, we have
\begin{eqnarray}\label{H0002}
\begin{array}{lll}
&& \int_{\tilde{M}} \ln (tr (h_{j, \epsilon }(t))+ tr (h_{j, \epsilon }(t)^{-1}))\frac{\omega_{\epsilon }^{n}}{n!}-\ln 2rank(\tilde{E})Vol (\omega_{\epsilon })\\
&\leq & 2t \int_{\tilde{M}}|\Lambda _{\omega_{\epsilon }}(F_{A_{j}})|_{\tilde{H}_{0}}\frac{\omega_{\epsilon }^{n}}{n!} .\\
\end{array}
\end{eqnarray}
Recall the result of Bando and Siu in \cite{BS}, by choosing a subsequence, we know that $H_{j, \epsilon }$ converges to a solution $H_{j}$ of the heat flow (\ref{DDD1}) on $M\setminus \Sigma $ as $\epsilon \rightarrow 0$. Then
\begin{eqnarray}\label{H02}
\begin{array}{lll}
&& \int_{M} \ln (tr (h_{j }(t))+ tr (h_{j }(t)^{-1}))\frac{\omega^{n}}{n!}-\ln 2rank(\tilde{E})Vol (\omega)\\
&\leq & 2t \int_{M}|\Lambda _{\omega }(F_{A_{j}})|_{H_{0}}\frac{\omega^{n}}{n!} .\\
\end{array}
\end{eqnarray}

On the other hand, we have
\begin{eqnarray}\label{H03}
\begin{array}{lll}
&& \triangle \ln (tr (h_{j }(t))+ tr (h_{j }(t)^{-1}))\\
&\geq &-2 |\Lambda _{\omega }(F_{H_{j }(t)}) |_{H_{j }(t)} -2 |\Lambda _{\omega }(F_{A_{j}}) |_{H_{0}}\\
\end{array}
\end{eqnarray}
on $M\setminus \Sigma $, for all $t>0$.

For any compact subset $\Omega \subset M \setminus (\Sigma\cup \Sigma_{An})$, setting $d(\Omega)=\inf \{\rho (x, y) | x\in \Omega , y\in \Sigma\cup \Sigma_{An}\}>0$, where $\rho $ is the distance function on $(M, \omega )$. Let $B=\cup_{y\in \Sigma \cup \Sigma_{An} } B_{y}(\frac{1}{2} d(\Omega))$ and $\Omega '=\tilde{M}-B$, then we choose the cut-off function $\varphi $ such that $\varphi \equiv 1$ on $\Omega $, $\varphi \equiv 0$ on $B$, and $|d \varphi |_{\omega }\leq \frac{4}{d(\Omega )}$. By the assumption, we known that $A_{j}$ are locally bounded in $C^{\infty}$ outside $\Sigma _{An}$, so we have
\begin{eqnarray}\label{C2}
|\Lambda _{\omega}(F_{A_{j}}) |_{H_{0}}\leq C_{c}
\end{eqnarray}
on $\Omega '$, where $C_{c}$ is a constant independent of  $j$. Using  (\ref{H02}), (\ref{H03}), (\ref{FFF1}), (\ref{C2}), the cut-off function $\varphi $ and  the Moser's iteration, we have
\begin{eqnarray}\label{mc0}
\begin{array}{lll}
&& \sup_{x\in \Omega } \ln (tr (h_{j }(1))+ tr (h_{j }(1)^{-1}))\\
&\leq & C_{d} \int_{M} \ln (tr (h_{j }(1))+ tr (h_{j }(1)^{-1}))\frac{\omega^{n}}{n!}\\
&\leq & C_{e},\\
\end{array}
\end{eqnarray}
where $C_{d}$, $C_{e}$ are constants independent of $j$. In a similar way, we have locally $C^{0}$ bound on metrics $H^{S}(1)$, i.e. for any compact subset $\Omega $, there exists a constant $C_{f}$  such that
 \begin{eqnarray}\label{mc0s}
 \sup_{x\in \Omega } \ln (tr ((H^{S}_{0})^{-1}(H^{S} (1)))+ tr ((H^{S} (1))^{-1}H^{S}_{0}))\leq C_{f}.
\end{eqnarray}

\medskip

By (\ref{key11}) and rescale $\eta_{j}$, we have a sequence $\phi $-invariant $\overline{\partial }_{0 , j}$-holomorphic map $f_{j}$ such that
\begin{eqnarray}\label{C011}
|f_{j}|^{2}_{H^{S}(1), H_{j } (1) } \leq C_{a}, \quad  \int_{M}|f_{j}|^{2}_{H^{S}(1), H_{j } (1)}\frac{\omega^{n}}{n!}=1.
\end{eqnarray}

 For any compact subset $\Omega \subset M \setminus (\Sigma \cup \Sigma _{An})$, by (\ref{C011}) and the above uniform locally $C^{0}$ bound on $H_{j } (1)$ and $H^{S} (1)$ (i.e. (\ref{mc0}) and (\ref{mc0s})), we have
 \begin{eqnarray}\label{C0111}
 \sup_{x\in \Omega }|f_{j}|^{2}_{H_{0}^{S}, H_{0 } }(x)\leq C(\Omega),
 \end{eqnarray}
where $C(\Omega)$ is a constant independent of $j$.

 Since $f_{j}$ is $\overline{\partial }_{0 , j}$-holomorphic, we have
\begin{eqnarray}
\begin{array}{lll}
\triangle _{0, j} f_{j}&=&\sqrt{-1}\Lambda_{\omega }(\partial_{0 ,
j}\overline{\partial_{0 , j}}-\overline{\partial_{0 , j}}\partial_{0
, j})f_{j}\\
&=&-\sqrt{-1}\Lambda_{\omega }(\partial_{0 , j}\overline{\partial_{0
, j}}+\overline{\partial_{0 , j}}\partial_{0
, j})f_{j}\\
&=& -\sqrt{-1}\Lambda_{\omega }(F_{A_{j}}f_{j}-f_{j}F_{A_{0}}),\\
\end{array}
\end{eqnarray}
and
\begin{eqnarray}\label{6.1}
\overline{\partial }_{0 , \infty} f_{j}=\overline{\partial }_{
A_{\infty }}\circ f_{j} -f_{j}\circ \overline{\partial }_{A_{0} } =
-\beta_{j}^{0, 1}\circ f_{j},
\end{eqnarray}
where $\beta
_{j}=A_{j}-A_{\infty}$.

By the above uniform locally $C^{0}$ bound of $f_{j}$ (i.e. (\ref{C0111})) and  the assumption that $A_{j}\rightarrow A_{\infty}$ in $C^{\infty}_{loc}$ topology outside $\Sigma_{An}$, the elliptic theory implies that there exists a subsequence of $f_{j}$ (for simplicity, also denoted by $f_{j}$) such that $f_{j} \rightarrow f_{\infty}$ in $C^{\infty}_{loc}$ topology outside $\Sigma \cup \Sigma_{An}$, and
\begin{eqnarray}
\overline{\partial }_{A_{0} , A_{\infty}} f_{\infty}=0 , \quad f_{\infty}\circ
\phi_{0} =\phi_{\infty} \circ f_{\infty}.
\end{eqnarray}

Since  $\Sigma \cup \Sigma_{An}$ is Hausdorff codimension at least $2$, for any small $\delta >0$, we can choose a compact subset $\Omega_{\delta} \subset M \setminus (\Sigma \cup \Sigma_{An})$ such that
\begin{eqnarray}\label{dd}
\int_{M\setminus \Omega_{\delta}} 1 \frac{\omega^{n}}{n!}\leq \delta .
\end{eqnarray}

From (\ref{C011}) and (\ref{dd}), we have
\begin{eqnarray}\label{n0}
      \int_{\Omega_{\delta }}|f_{j}|^{2}_{H^{S}(1), H_{j } (1) }\frac{\omega ^{n}}{n!}\geq 1-\delta C_{a}.
\end{eqnarray}

Using the above uniform locally $C^{0}$ bound on $H_{j } (1)$ and $H^{S}(1)$ ((\ref{mc0}) and (\ref{mc0s})) again, we have a positive constant $\tilde{C}(\Omega_{\delta })$ such that,
\begin{eqnarray}
      \int_{\Omega_{\delta }}|f_{j}|^{2}_{H_{0}^{S}, H_{0 }}\frac{\omega^{n}}{n!}\geq \tilde{C}(\Omega_{\delta })(1-\delta C_{a})>0 ,
\end{eqnarray}
for every $j$. Then, we get $
      \int_{\Omega_{\delta }}|f_{\infty}|^{2}_{H_{0}^{S}, H_{0 } }\frac{\omega^{n}}{n!}>0
$, and
so $f_{\infty}$ is a non-zero holomorphic map.

\hfill $\Box$ \\

\medskip

\section{ The HN type of the Uhlenbeck limit}
\setcounter{equation}{0}

Let $(A_{t} , \phi_{t})$ be a smooth solution of the gradient heat
flow (\ref{YMHH}) over a K\"ahler manifold $(M, \omega )$ with initial data $(A_{0} ,
\phi_{0})$, and let $(A_{\infty} , \phi_{\infty})$ be an Uhlenbeck
limit. From Theorem 2.7., we know that $(A_{\infty} ,
\phi_{\infty})$ is a smooth Yang-Mills Higgs pair on Hermitian
bundle $(E_{\infty} , H_{\infty})$ over $M\setminus \Sigma_{an}$, and $\theta (A_{\infty} ,
\phi_{\infty})$ is parallel, then the constant eigenvalues vector
$\vec{\lambda}_{\infty}=(\lambda_{1} , \cdots , \lambda_{R})$ of
$\sqrt{-1}\theta (A_{\infty} , \phi_{\infty})$ is just the HN type
of the extended Uhlenbeck limit Higgs sheaf $(E_{\infty} , A_{\infty} ,
\phi_{\infty})$. Let $\vec{\mu}$ be the HN type  of the initial
Higgs bundle $(E , A_{0} , \phi_{0})$. In this section, we will
prove that $\vec{\lambda }_{\infty}=\vec{\mu}$ .

Let $\verb"u"(R)$ denote the Lie algebra of the unitary group
$U(R)$. Fix a real number $\alpha \geq 1$,  for any $\verb"a"\in
\verb"u"(R)$, let $\varphi_{\alpha
}(\verb"a")=\sum_{j=1}^{R}|\lambda_{j}|^{\alpha }$, where
$\sqrt{-1}\lambda_{j}$ are eigenvalues of $\verb"a"$. It is easy
to see that we can find a family $\varphi_{\alpha , \rho }$ of smooth convex ad-invariant functions
$0<\rho \leq 1$,  such that
$\varphi_{\alpha , \rho}\rightarrow \varphi_{\alpha }$ uniformly
on compact subsets of $\verb"u"(R)$ as $\rho \rightarrow 0$.
Hence, from \cite{AB} ( Prop.12.16) it follows that $\varphi_{\alpha }$
is a convex function on $\verb"u"(R)$. For a given real number
$N$, define the Hermitian-Yang-Mills type functionals as follows:
\begin{eqnarray}
HYM_{\alpha , N}(A, \phi )=\int_{M}\varphi_{\alpha }(\theta (A,
\phi ) -\sqrt{-1}N Id_{E}) \frac{\omega^{n}}{n!}.
\end{eqnarray}
In the following we assume that $Vol(M, \omega )=1$, and set $HYM_{\alpha , N}(\vec{\mu })=HYM_{\alpha
,}(\vec{\mu }+N)= \varphi_{\alpha }(\sqrt{-1}((\vec{\mu
}+N)))$, where $\vec{\mu }+N=diag (\mu_{1}+N , \cdots ,
\mu_{R}+N)$. We need the following two lemmas, the proofs can
be found in \cite{DW1} ( Lemma 2.23 ; Prop.2.24).

\hspace{0.3cm}

{\bf Lemma 5.1.} {\it The functional $\verb"a" \mapsto (\int_{M}
\varphi_{\alpha }(\verb"a") dvol )^{\frac{1}{\alpha }}$, defines a
norm on $L^{\alpha }(\verb"u"(E))$ which is equivalent to the
$L^{\alpha }$ norm.}

\hspace{0.3cm}

{\bf Lemma 5.2.} {\it (1) If $\vec{\mu}\leq \vec{\lambda }$, then
$\varphi_{\alpha }(\sqrt{-1}\vec{\mu})\leq \varphi_{\alpha
}(\sqrt{-1}\vec{\lambda })$ for all $\alpha \geq 1$.

(2) Assume $\mu_{R}\geq 0$ and $\lambda_{R}\geq 0$. If
$\varphi_{\alpha }(\sqrt{-1}\vec{\mu})= \varphi_{\alpha
}(\sqrt{-1}\vec{\lambda })$ for all $\alpha $ in some set $S\subset
[1 , \infty )$ possessing a limit point, then $\vec{\mu}=
\vec{\lambda }$.}

\medskip

For any smooth convex ad-invariant functions $\varphi $, we have
\begin{eqnarray}\label{2.33}
(\triangle -\frac{\partial }{\partial t})\varphi(\theta (A_{t},
\phi_{t} ) -\sqrt{-1}N Id_{E})\geq 0,
\end{eqnarray}
whose proof can be found in \cite{LZ1} (Section two). Since we can approximate $\varphi_{\alpha }$ by
smooth convex ad-invariant functions $\varphi_{\alpha , \rho}
\rightarrow \varphi_{\alpha }$, by (\ref{2.33}) we know that $t\mapsto HYM_{\alpha , N}(A_{t} , \phi_{t})$ is
nonincreasing along the flow.
By Corollary 2.8., we can choose a sequence
$t_{j}\rightarrow \infty$, such that
\begin{eqnarray}
HYM_{\alpha , N}(A_{t_{j}}
, \phi_{t_{j}})\rightarrow HYM_{\alpha , N}(A_{\infty} ,
\phi_{\infty}).
\end{eqnarray}
 Then we have the following proposition.

 \medskip

{\bf Proposition 5.3. } {\it Let $(A_{t} , \phi_{t})$ be a
solution of the gradient flow (\ref{YMHH}) and $(A_{\infty} ,
\phi_{\infty })$ be a subsequential Uhlenbeck limit of $(A_{t} ,
\phi_{t})$. Then for any $\alpha \geq 1 $ and any $N$, $t\mapsto
HYM_{\alpha , N}(A_{t} , \phi_{t})$ is nonincreasing, and
$lim_{t\rightarrow \infty }HYM_{\alpha , N}(A_{t} ,
\phi_{t})=HYM_{\alpha , N}(A_{\infty} , \phi_{\infty})$.}

\medskip

{\bf Lemma 5.4 } {\it Let $(A_{j} , \phi_{j})=g_{j}(A_{0} ,
\phi_{0})$ be a sequence of complex gauge equivalent Higgs pairs on
a complex vector bundle $E$ of rank $R$ with Hermitian metric
$H_{0}$. Let $S$ be a coherent $\phi_{0}$-invariant subsheaf of $(E
, A_{0})$. Suppose that $\sqrt{-1}\Lambda_{\omega
}(F_{A_{j}}+[\phi_{j} , \phi_{j}^{\ast }])\rightarrow \verb"a"$ in
$L^{1}$ as $j\rightarrow \infty $, where $\verb"a"\in
L^{1}(\sqrt{-1}\verb"u"(E))$, and that eigenvalues $\lambda_{1 }\geq
\cdots \geq \lambda_{R}$ of $\frac{1}{2\pi}\verb"a"$ are constant almost
everywhere. Then: $deg (S)\leq \sum_{i\leq rank(S)}\lambda_{i}$.}

\hspace{0.3cm}

{\bf Proof.} Since $deg (S)\leq deg(Sat_{E}(S))$, we may assume that
$S$ is saturated. Let $\pi_{j}$ denote the orthogonal projection
onto $g_{j}(S)$ with respect to the Hermitian metric $H_{0}$. It is
well known that $\pi_{j}$ are $L_{1}^{2}$ sections of the smooth
endomorphism bundle of $E$, and satisfy
$\pi_{j}^{2}=\pi_{j}=\pi_{j}^{\ast}$,
$(Id-\pi_{j})\overline{\partial }_{A_{j}}\pi_{j}=0$ and
$(Id-\pi_{j})\phi_{i}\pi_{j}=0$ (since $g_{j}(S)$ are
$\phi_{j}$-invariant). By the usual degree formula (see
Lemma 3.2 in \cite{Si}), we have
\begin{eqnarray}\label{5.4}
\begin{array}{lll}
deg(S)&=&\frac{1}{2\pi}\int_{M}(Tr(\sqrt{-1}\Lambda_{\omega}(F_{A_{j}}+[\phi_{j} ,
\phi_{j}^{\ast}])\pi_{j})-|\overline{\partial }_{A_{j}+\phi_{j}}\pi_{j}|^{2})\\
&\leq & \frac{1}{2\pi}\int_{M}(Tr(\sqrt{-1}\Lambda_{\omega}(F_{A_{j}}+[\phi_{j} ,
\phi_{j}^{\ast}])\pi_{j})\\
&= & \frac{1}{2\pi}\int_{M}(Tr(\verb"a"\pi_{j})+\frac{1}{2\pi}\int_{M}(Tr(\sqrt{-1}\Lambda_{\omega}((F_{A_{j}}+[\phi_{j} ,
\phi_{j}^{\ast}])-\verb"a")\pi_{j}).
\end{array}
\end{eqnarray}
By a result from linear algebra (Lemma 2.20 in \cite{DW1}), we have $\frac{1}{2\pi}Tr(\verb"a"\pi_{j})\leq \sum_{i\leq rank(S)}\lambda_{i}$. So, we have
$deg (S)\leq \sum_{i\leq rank(S)}\lambda_{i} +\frac{1}{2\pi}\|\sqrt{-1}\Lambda_{\omega}(F_{A_{j}}+[\phi_{j} ,
\phi_{j}^{\ast}])-\verb"a"\|_{L^{1}}$. Let $j\rightarrow \infty $, this completes the proof of the lemma.

\hfill $\Box$ \\

\medskip

Let $\vec{\mu}_{0}=(\mu_{1} ,
\cdots , \mu_{R})$ be the HN type of Higgs bundle $(E , A_{0} ,
\phi_{0})$, by (\ref{F1}) in Lemma 2.1, we have
\begin{eqnarray}
\sum_{\alpha =1}^{R}\mu_{\alpha} =deg (E, \overline{\partial}_{A_{0}})=deg (E_{\infty}, \overline{\partial}_{A_{\infty}})=\sum_{\alpha =1}^{R}\lambda_{\alpha}.
\end{eqnarray}
Let $\{E_{i}\}_{i=1}^{l}$ be the HN filtration of the Higgs bundle $(E, A_{0}, \phi_{0})$.
Using Corollary 2.8. and Lemma 5.4.,  we have:
\begin{eqnarray}
\sum_{\alpha \leq rank E_{i}} \mu_{\alpha }=deg(E_{i})\leq \sum_{\alpha \leq rank E_{i}} \lambda_{\alpha }
\end{eqnarray}
for all $i$. By Lemma 2.3 in \cite{DW1}, we have the following proposition.

\medskip

{\bf Proposition 5.5. } {\it Let $(A_{i} , \phi_{i})$ be a sequence
of Higgs pairs along the gradient heat flow (1.6) with Uhlenbeck
limit $(A_{\infty } , \phi_{\infty})$. Let $\vec{\mu}_{0}=(\mu_{1} ,
\cdots , \mu_{R})$ be the HN type of Higgs bundle $(E , A_{0} ,
\phi_{0})$, and let $\vec{\lambda}_{\infty}=(\lambda_{1} , \cdots ,
\lambda_{R})$ be the type of Higgs bundle $(E_{\infty } , A_{\infty}
, \phi _{\infty})$. Then $\vec{\mu }_{0}\leq
\vec{\lambda}_{\infty}$.}

\hspace{0.3cm}

\hspace{0.3cm}

Let $H$ be a smooth Hermitian metric on the holomorphic bundle
$\textbf{E}=(E, \overline{\partial}_{E})$, and let $\textit{F}=\{F_{i}\}_{i=1}^{l}$ be a
filtration of $E$ by saturated subsheaves: $0=F_{0}\subset
F_{1}\subset \cdots \subset F_{l-1}\subset F_{l}=\textbf{E}$.
Associated to each $F_{i}$ and the metric $H$ we have the unitary
projection $\pi_{i}^{H}$ onto $F_{i}$. It is well known that
$\pi_{i}^{H}$ are bounded $L_{1}^{2}$ Hermitian endomorphisms. For
convenience, we set $\pi_{0}^{H}=0$. Given real numbers $\mu_{1}
,\cdots , \mu_{l}$ and a filtration $\textit{F}$, we define a
bounded $L_{1}^{2}$ Hermitian endomorphism of $\textbf{E} $ by
$\Psi (\textit{F }, (\mu_{1} ,\cdots , \mu_{l}), H
)=\Sigma_{i=1}^{l}\mu_{i}(\pi_{i}^{H}-\pi_{i-1}^{H})$. Given a
Hermitian metric on a Higgs bundle $(\textbf{E} , \phi  )$, the
Harder-Narasimhan projection, $\Psi_{\omega }^{HN}(\textbf{E},
\phi , H)$ is the bounded $L_{1}^{2}$ Hermitian endomorphism
defined above in the particular case where $\textit{F}$ is the HN
filtration $F_{i}=\texttt{F}_{i}^{hn}(E)$ and $\mu_{i}=\mu
(F_{i}/F_{i-1})$.

\hspace{0.3cm}

{\bf Definition 5.6. } {\it Fix $\delta >0$ and $1\leq p\leq
\infty$. An $L^{p}$-$\delta $-approximate critical Hermitian
metric on a Higgs bundle $(\textbf{E} , \phi )$ is a smooth $H$
such that
$$\|\frac{\sqrt{-1}}{2\pi}\Lambda_{\omega }(F_{A_{H}}+[\phi , \phi^{\ast H
}])-\Psi ^{HN}(\textbf{E}, \phi , H)\|_{L^{p}(\omega )}\leq
\delta ,
$$where $A_{H}$ is the Chern connection determined by
$(\overline{\partial }_{E} , H)$.}

\medskip

{\bf Proposition 5.7. } {\it Let $(\textbf{E} , \phi )$ be a
Higgs bundle on a smooth K\"ahler manifold $(M , \omega )$, and
and let $\textit{F}=\{F_{i}\}_{i=1}^{l}$ be a
filtration of $E$ by saturated subsheaves, where every $F_{i}$ is $\phi $ invariant. Let  $\pi :
\overline{M}\rightarrow M$ be a blow-up along a smooth complex manifold $\Sigma$ of complex co-dimensional at least $2$, $\overline{\textbf{E}}=\pi^{\ast }\textbf{E}$ be the pull-back bundle and $\overline{\phi }=\pi ^{\ast } \phi $. Let $\eta $ be a K\"ahler metric on $\overline{M}$, and  set a family of
K\"ahler metrics $\omega_{\epsilon }=\pi^{\ast
}\omega +\epsilon \eta $. Suppose that the filtration $\overline{\textit{F}}=\{\overline{F}_{i}\}_{i=1}^{l}=\{Sat_{\overline{E}}(\pi^{\ast }F_{i})\}_{i=1}^{l}$  of
$\overline{\textbf{E}}$ is given by subbundles, and every quotient $\overline Q_{i}=\overline{F}_{i}/\overline{F}_{i-1}$ is $\overline{\phi}$-Higgs $\omega_{\epsilon }$-stable for $0<\epsilon < \epsilon ^{\ast}$.  Then for any $\tilde{\delta } >0$ and any $0<\epsilon < \epsilon ^{\ast}$, there is
 a smooth Hermitian metric $\overline{H}$ on
$\overline{\textbf{E}}$ such that
 \begin{eqnarray}\label{a1}
\|\frac{\sqrt{-1}}{2\pi}\Lambda_{\omega_{\epsilon}
}(F_{(\overline{\partial }_{\overline{E}} ,
\overline{H})}+[\overline{\phi} , \overline{\phi}^{\ast \overline{H}
}])-\Psi(\overline{\textit{F}}, (\mu_{\epsilon , 1} , \cdots , \mu_{\epsilon , l}) ,
\overline{H})\|_{L^{\infty}}\leq \tilde{\delta } ,
\end{eqnarray}
where $(\overline{\partial }_{\overline{E}} ,
\overline{H})$ denotes the Chern connection with respect to holomorphic structure $\overline{\partial }_{\overline{E}}$ and metric $\overline{H}$ , and $\mu_{\epsilon , i}$ is the slope of quotient $\overline{Q}_{i}$ with respect to the metric $\omega_{\epsilon }$.}

 \hspace{0.3cm}

{\bf Proof. } Let $\overline{\phi}_{i}$ be the induced Higgs field on the quotient $\overline{Q}_{i}$. Since  Higgs bundles
$(\overline{Q_{i}}, \overline{\phi}_{i})$ are
$\omega_{\epsilon}$-stable for all $0<\epsilon \leq \epsilon^{\ast}$,
by Simpson's result (Theorem 1 in \cite{Si}), we have a Hermitian-Einstein metric
$\overline{H}_{i , \epsilon}$ on the Higgs bundle $(\overline{Q_{i}},
\overline{\phi}_{i})$ with respect to
$\omega_{\epsilon}$. In particular:
\begin{eqnarray}
\frac{\sqrt{-1}}{2\pi}\Lambda_{\omega_{\epsilon}
}(F_{(\overline{\partial }_{\overline{Q}_{i}}, \overline{H}_{i,\epsilon})}+[\overline{\phi}_{i} ,
\overline{\phi}_{i}^{\ast }])-\mu_{\epsilon
}(\overline{Q_{i}})Id_{\overline{Q}_{i}}=0.
\end{eqnarray}

We will use Donaldson's argument in \cite{D1}. Recall that $\overline{E}$ and $\oplus _{i} \overline{Q}_{i}$ are isomorphic vector bundle,  and take the direct sum $\overline{H}_{\epsilon }=\oplus _{i} \overline{H}_{i , \epsilon}$. By the equivalence of holomorphic structures and integrable unitary connections, we see that it suffices to show that for a fixed Hermitian metric $\overline{H}$ there is  a smooth complex gauge transformation $\sigma $ preserving the filtration $\overline{\textit{F}}$ such that
\begin{eqnarray}
\begin{array}{lll}
&&\|\frac{\sqrt{-1}}{2\pi}\Lambda_{\omega_{\epsilon}
}(F_{(\sigma (\overline{\partial }_{\overline{E}}) ,
\overline{H})}+[\overline{\sigma(\phi )} , \sigma (\overline{\phi})
^{\ast \overline{H}
}])-\Psi(\overline{\textit{F}}, (\mu_{\epsilon , 1} , \cdots , \mu_{\epsilon , l}) ,
\overline{H})\|_{L^{\infty}}\\& &\leq  \delta ,\\
\end{array}
\end{eqnarray}

We only  consider the case $l=2$, the other case $l>2$ can be solved by inductive argument. We can write the holomorphic structure $\overline{\partial }_{\overline{E}}$ and the Higgs field $\overline{\phi }$ as
\begin{eqnarray}
\overline{\partial }_{\overline{E}}=\left (
\begin{matrix}
\overline{\partial}_{\overline{Q}_{1}} &  B \\
0   & \overline{\partial}_{\overline{Q}_{2}}\\
\end{matrix}
\right ),
\quad
 \overline{\phi } = \left (
\begin{matrix}
\overline{\phi}_{1} &  \zeta \\
0   & \overline{\phi }_{2}\\
\end{matrix}
\right ),
\end{eqnarray}
where $B$ is the second fundamental form. Define the complex gauge transformation $\sigma_{t}$ to be the following block diagonal matrix
\begin{eqnarray}
\sigma_{t}=\left (
\begin{matrix}
t Id_{\overline{Q}_{1}} &  0 \\
0   & t^{-1} Id_{\overline{Q}_{2}}\\
\end{matrix}
\right ).
\end{eqnarray}
Then, we have
\begin{eqnarray}
\sigma_{t}(\overline{\partial }_{\overline{E}})=\left (
\begin{matrix}
\overline{\partial}_{\overline{Q}_{1}} &  t^{2}B \\
0   & \overline{\partial}_{\overline{Q}_{2}}\\
\end{matrix}
\right ),
\quad
\sigma_{t}(\overline{\phi })= \left (
\begin{matrix}
\overline{\phi}_{1} &  t^{2}\zeta \\
0   & \overline{\phi }_{2}\\
\end{matrix}
\right ),
\end{eqnarray}
\begin{eqnarray}
F_{(\sigma_{t} (\overline{\partial }_{\overline{E}}) ,
\overline{H}_{\epsilon})}=\left (
\begin{matrix}
F_{( \overline{\partial }_{\overline{Q}_{1}} ,
\overline{H}_{1, \epsilon
})}+t^{4}B\wedge B^{\ast} &  t^{2}\partial _{H}B \\
t^{2}\overline{\partial} B^{\ast}   & F_{(\overline{\partial }_{\overline{Q}_{2}} ,
\overline{H}_{2, \epsilon })}+t^{4}B^{\ast}\wedge B\\
\end{matrix}
\right ),
\end{eqnarray}
and
\begin{eqnarray}
[\sigma_{t}(\overline{\phi }), \sigma_{t}(\overline{\phi })^{\ast}]= \left (
\begin{matrix}
[\overline{\phi}_{1}, \overline{\phi}_{1}^{\ast}]+t^{4}\zeta \wedge \zeta^{\ast} &  t^{2}(\zeta\wedge \overline{\phi}_{2}^{\ast}+\overline{\phi}_{1}^{\ast}\wedge \zeta ) \\
 t^{2}(\zeta^{\ast}\wedge \overline{\phi}_{1}+\overline{\phi}_{2}\wedge \zeta ^{\ast})   & [\overline{\phi}_{2}, \overline{\phi}_{2}^{\ast}]+t^{4}\zeta^{\ast} \wedge \zeta \\
\end{matrix}
\right )
.
\end{eqnarray}

We can also write $\Psi(\overline{\textit{F}}, (\mu_{\epsilon , 1} , \mu_{\epsilon , 2}) ,
\overline{H}_{\epsilon})$ as follows
\begin{eqnarray}
\left (
\begin{matrix}
\mu_{\epsilon , 1} Id_{\overline{Q}_{1}}&  0 \\
0   & \mu_{\epsilon , 2}Id_{\overline{Q}_{2}}\\
\end{matrix}
\right ).
\end{eqnarray}
Then, it is easy to see that
\begin{eqnarray}
\begin{array}{lll}
&&|\frac{\sqrt{-1}}{2\pi}\Lambda_{\omega_{\epsilon}
}(F_{(\sigma_{t} (\overline{\partial }_{\overline{E}}) ,
\overline{H}_{\epsilon})}+[\sigma_{t}(\overline{\phi }) , \sigma_{t} (\overline{\phi})
^{\ast \overline{H}_{\epsilon}
}])-\Psi(\overline{\textit{F}}, (\mu_{\epsilon , 1} , \mu_{\epsilon , 2}) ,
\overline{H}_{\epsilon})|_{\overline{H}_{\epsilon}}\\
&&\leq  \sum_{i=1}^{2}|\frac{\sqrt{-1}}{2\pi}\Lambda_{\omega_{\epsilon}
}(F_{(\overline{\partial }_{\overline{Q}_{i}}, \overline{H}_{i,\epsilon})}+[\overline{\phi}_{i} ,
\overline{\phi}_{i}^{\ast }])-\mu_{\epsilon
}(\overline{Q_{i}})Id_{\overline{Q}_{i}}|_{\overline{H}_{i, \epsilon}}\\
&& + f(t, B),
\end{array}
\end{eqnarray}
where $f(t , B) \rightarrow 0$ as $t\rightarrow 0$. Let $\overline{H}_{t}=t^{2}\overline{H}_{1, \epsilon}\oplus t^{-2}\overline{H}_{2, \epsilon }$,  then $\sigma_{t}^{\ast \overline{H}_{\epsilon}}\circ \sigma_{t}=\overline{H}_{\epsilon}^{-1}\overline{H}_{t}$. On the other hand, we have
\begin{eqnarray}
\begin{array}{lll}
&&|\frac{\sqrt{-1}}{2\pi}\Lambda_{\omega_{\epsilon}
}(F_{(\overline{\partial }_{\overline{E}} ,
\overline{H}_{t})}+[\overline{\phi } , \overline{\phi}
^{\ast \overline{H}_{t}
}])-\Psi(\overline{\textit{F}}, (\mu_{\epsilon , 1} , \mu_{\epsilon , 2}) ,
\overline{H}_{t})|_{\overline{H}_{t}}\\
&=& |\frac{\sqrt{-1}}{2\pi}\Lambda_{\omega_{\epsilon}
}(F_{(\sigma_{t} (\overline{\partial }_{\overline{E}}) ,
\overline{H}_{\epsilon})}+[\sigma_{t}(\overline{\phi }) , \sigma_{t} (\overline{\phi})
^{\ast \overline{H}_{\epsilon}
}])-\Psi(\overline{\textit{F}}, (\mu_{\epsilon , 1} , \mu_{\epsilon , 2}) ,
\overline{H}_{\epsilon})|_{\overline{H}_{\epsilon}}\\
\end{array}
\end{eqnarray}
Choosing $t$ small enough, we obtain a metric $\overline{H}$ which satisfies (\ref{a1}).

\hfill $\Box$ \\

\medskip

The following lemma was proved by Sibley in \cite{Sib} (Lemma 5.3.),  we give a proof for reader's convenience.

\medskip

{\bf Lemma 5.8. } {\it Let $(M, \omega )$ be a compact K\"ahler manifold of complex dimension $n$, and  $\pi :
\overline{M}\rightarrow M$ be a blow-up along a smooth complex sub-manifold $\Sigma$ of complex co-dimension $k$ where $k\geq 2$. Let $\eta $ be a K\"ahler metric on $\overline{M}$, and consider the family of
K\"ahler metrics $\omega_{\epsilon }=\pi^{\ast
}\omega +\epsilon \eta $. Then for any $0\leq \gamma < \frac{1}{k-1}$, we have $\frac{\eta ^{n}}{ \omega_{\epsilon}^{n}}\in L^{\gamma}(\overline{M}, \eta )$, and the $L^{\gamma}(\overline{M}, \eta )$-norm of $\frac{\eta ^{n}}{ \omega_{\epsilon}^{n}}$ is uniformly bounded in $\epsilon$, i.e. there is a positive constant $C^{\ast}$ such that
\begin{eqnarray}\label{Int1}
\int_{\overline{M}}( \frac{\eta ^{n}}{ \omega_{\epsilon}^{n}})^{\gamma } \frac{\eta ^{n}}{n!}\leq C^{\ast}
\end{eqnarray}
for all $\epsilon$.
}

\medskip

{\bf Proof. } Since $\pi^{\ast} \omega $ is only degenerated along the exceptional divisor $\pi^{-1}(\Sigma )$, on the complement of a neighborhood of $\pi^{-1}(\Sigma )$ there is a constant $C$ such that $C^{-1}\eta \leq \pi^{\ast}\omega \leq C \eta $. So, it is suffices to prove the result in a neighborhood of $\pi^{-1}(\Sigma )$. One can choose a local coordinate chart $U$ with coordinates $(z_{1}, \cdots , z_{n})$, such that locally $\Sigma $ is given by the slice $\{z_{1}=\cdots =z_{k}=0\}$. On the blow-up $\overline{M}$ we have local coordinate charts $\overline{U}_{i} \subset \pi^{-1}(U)$ where $\overline{U}_{i}=\{z\in U\setminus \Sigma | z_{i}\neq 0\}\cup \{(z, [v])\in \textbf{CP}(\Sigma)|_{U\cap \Sigma }|v_{i}\neq 0\}$, where $\textbf{CP}(\Sigma)$ is the projective bundle of the normal bundle of $\Sigma $ and $1\leq i\leq k$. Let $(w_{1}, \cdots , w_{n})$ denote local coordinates on $\overline{U}_{i}$, then the map $\pi : \overline{M} \rightarrow M$ is given by:
\begin{eqnarray*}
(w_{1}, \cdots , w_{n}) \rightarrow (w_{1}w_{i},\cdots , w_{i-1}w_{i}, w_{i}, w_{i+1}w_{i},\cdots , w_{k}w_{i}, w_{k+1} \cdots , w_{n}).
\end{eqnarray*}
Set $\omega =\sqrt{-1} g_{i \bar{j}} dz_{i}\wedge d\bar{z}_{j}$, then
\begin{eqnarray}
\pi^{\ast} \omega ^{n}=\pi^{\ast}(det g_{i \bar{j}}) |w_{i}|^{2k-2} (\sqrt{-1})^{n} d w_{1} \wedge d \bar{w}_{1} \wedge \cdots \wedge d w_{n} \wedge d \bar{w}_{n} .
\end{eqnarray}
Note that $\pi^{\ast}(det g_{i \bar{j}})$
is positive and $\pi^{\ast }\omega^{n}< \omega_{\epsilon }$, then we have
\begin{eqnarray}
\begin{array}{lll}
&&\int_{\overline{U}_{i}}( \frac{\eta ^{n}}{ \omega_{\epsilon}^{n}})^{\gamma } \frac{\eta ^{n}}{n!}\\
&\leq & \int_{\overline{U}_{i}}( \frac{\eta ^{n}}{ (\pi^{\ast} \omega )^{n}})^{\gamma } \frac{\eta ^{n}}{n!}\\
&\leq & C \int_{\overline{U}_{i}} |w_{i}|^{-(2k-2)\gamma } (\sqrt{-1})^{n} d w_{1} \wedge d \bar{w}_{1} \wedge \cdots \wedge d w_{n} \wedge d \bar{w}_{n},\\
\end{array}
\end{eqnarray}
where $C$ is a positive constant independent of $\epsilon$. By the condition we have $-(2k-2)\gamma >-2$, and then we see the result follows.

\hfill $\Box$ \\

\medskip

\medskip

{\bf Lemma 5.9. } {\it Let $\pi : \overline{M} \rightarrow M $, the complex submanifold $\Sigma $ and the family of metrics $\omega_{\epsilon }$ be the same as in the previous lemma. Then, for any $1<\alpha < 1+\frac{1}{2k-1}$, $\frac{\alpha }{1-(k-1)(\alpha -1)}<\tilde{\alpha }$ , and any neighborhood $U$ of the exceptional divisor $\pi^{-1}(\Sigma )$, there exists a positive constant $C $ independent of $\epsilon $, $\epsilon_{1}$ and $U$, and
a positive constant $C(U)$ depending only on $U$, $\alpha $, metric $\eta $, and $\pi^{\ast} \omega$, such that for any $End(E)$-valued $(1, 1)$ form $F$
\begin{eqnarray}\label{Int2}
\begin{array}{lll}
&&\|\Lambda_{\omega_{\epsilon }}F\|_{L^{\alpha }(\overline{M}, \omega_{\epsilon})} \leq  |\epsilon_{1}-\epsilon |C(U)\|F\|_{L^{\alpha }(\overline{M}, \omega_{\epsilon_{1}})}\\
&+& C(\|\Lambda_{\omega_{\epsilon_1 }}F\|_{L^{\tilde{\alpha } }(\overline{M}, \omega_{\epsilon_1})} +(Vol(U, \omega_{\epsilon_1}))^{C(\alpha , k)}\|F\|_{L^{2 }(\overline{M}, \omega_{\epsilon_{1}})})\\
\end{array}
\end{eqnarray}
for all $ 0<\epsilon \leq \epsilon_{1}\leq 1$, where $C(\alpha , k)$ is a positive constant depending only on $\alpha $ and $k$.
}

\medskip

{\bf Proof. } By the definition, we know that
\begin{eqnarray}
\Lambda_{\omega_{\epsilon }}F =\frac{\omega_{\epsilon_{1}}^{n}}{\omega_{\epsilon}^{n}}(\Lambda_{\omega_{\epsilon_{1} }}F +n \frac{F \wedge (\omega_{\epsilon}^{n-1}- \omega_{\epsilon_{1} }^{n-1})}{\omega_{\epsilon_{1}}^{n}}),
\end{eqnarray}
and then
\begin{eqnarray}\label{LLL1}
\begin{array}{lll}
\|\Lambda_{\omega_{\epsilon }}F\|_{L^{\alpha }(\overline{M}, \omega_{\epsilon})} &\leq & \|\frac{\omega_{\epsilon_{1}}^{n}}{\omega_{\epsilon}^{n}}(\Lambda_{\omega_{\epsilon_{1} }}F)\|_{L^{\alpha }(\overline{M}, \omega_{\epsilon})}\\
&+& \|n \frac{F \wedge (\omega_{\epsilon}^{n-1}- \omega_{\epsilon_{1} }^{n-1})}{\omega_{\epsilon_{1}}^{n}}\|_{L^{\alpha }(\overline{M}, \omega_{\epsilon})}.\\
\end{array}
\end{eqnarray}
Firstly, we have
\begin{eqnarray}\label{LLL2}
\begin{array}{lll}
&&\|\frac{\omega_{\epsilon_{1}}^{n}}{\omega_{\epsilon}^{n}}(\Lambda_{\omega_{\epsilon_{1} }}F )\|_{L^{\alpha }(\overline{M}, \omega_{\epsilon})}^{\alpha }
= \int_{\overline{M}} |\Lambda_{\omega_{\epsilon_{1} }}F|^{\alpha } (\frac{\omega_{\epsilon_{1}}^{n}}{\omega_{\epsilon}^{n}})^{\alpha -1} \frac{\omega_{\epsilon_{1}}^{n}}{n!}\\
&\leq & (\int_{\overline{M}} |\Lambda_{\omega_{\epsilon_{1} }}F|^{\alpha \cdot p } \frac{\omega_{\epsilon_{1}}^{n}}{n!})^{\frac{1}{p}}
(\int_{\overline{M}}  (\frac{\omega_{\epsilon_{1}}^{n}}{\omega_{\epsilon}^{n}})^{(\alpha -1)\cdot q} \frac{\omega_{\epsilon_{1}}^{n}}{n!})^{\frac{1}{q}},
\end{array}
\end{eqnarray}
where $\alpha \cdot p=\tilde{\alpha }$, since $\frac{\alpha }{1-(k-1)(\alpha -1)}<\tilde{\alpha }$ and $\frac{1}{p}+\frac{1}{q}=1$, we have $(\alpha -1)q <\frac{1}{k-1}$.

Since $\overline{M}$ is compact, there is a constant $C_{M}$ such that $\pi^{\ast}\omega \leq C_{M}\eta$ on $\overline{M}$. Let $U$ be a neighborhood of the exceptional divisor $\pi^{-1}(\Sigma ) $, since $\pi^{
\ast}\omega $ is degenerated only along $\pi^{-1}(\Sigma ) $, we can suppose that $\pi^{\ast}\omega \geq C_{u}\eta $ on $\overline{M}\setminus U$ for some positive constant $C_{u}$. Then, we have
\begin{eqnarray}\label{LLL3}
\begin{array}{lll}
&&\int_{\overline{M} \setminus U}|\frac{F \wedge (\omega_{\epsilon}^{n-1}- \omega_{\epsilon_{1} }^{n-1})}{\omega_{\epsilon}^{n}}|^{
\alpha } \frac{\omega_{\epsilon}^{n}}{n!}\\&=&\int_{\overline{M} \setminus U}|\epsilon_{1}-\epsilon |^{\alpha }|\frac{F \wedge \eta  \wedge (\sum_{i=0}^{n-2}\omega_{\epsilon}^{n-i-2}\wedge \omega_{\epsilon_{1} }^{i})}{\omega_{\epsilon_{1}}^{n}}|^{
\alpha } (\frac{\omega_{\epsilon_{1}}^{n}}{\omega_{\epsilon}^{n}})^{(\alpha -1)} \frac{\omega_{\epsilon_{1}}^{n}}{n!}\\
&\leq &C(n)C_{u}^{-(n+1)\alpha } |\epsilon_{1}-\epsilon |^{\alpha }\int_{\overline{M} \setminus U} |F|_{\omega_{\epsilon_1}}^{\alpha }(\frac{\omega_{\epsilon_{1}}^{n}}{\omega_{\epsilon}^{n}})^{(\alpha -1)} \frac{\omega_{\epsilon_{n}}^{n}}{n!}\\
&\leq &C(n)C_{u}^{-(n+1)\alpha +n }(C_{M}+\epsilon_{1})^{n(\alpha -1)} |\epsilon_{1}-\epsilon |^{\alpha }\int_{\overline{M} \setminus U} |F|_{\omega_{\epsilon_1}}^{\alpha } \frac{\omega_{\epsilon_{1}}^{n}}{n!}.\\
\end{array}
\end{eqnarray}
On the other hand, we have
\begin{eqnarray}\label{LLL4}
\begin{array}{lll}
&&\int_{ U}|\frac{F \wedge (\omega_{\epsilon}^{n-1}- \omega_{\epsilon_{1} }^{n-1})}{\omega_{\epsilon}^{n}}|^{
\alpha } \frac{\omega_{\epsilon}^{2}}{n!}\\&= &\int_{ U}|\frac{F \wedge (\omega_{\epsilon}^{n-1}- \omega_{\epsilon_{1} }^{n-1})}{\omega_{\epsilon_{1}}^{n}}|^{
\alpha }(\frac{\omega_{\epsilon_{1}}^{n}}{\omega_{\epsilon}^{n}})^{(\alpha -1)} \frac{\omega_{\epsilon_{1}}^{n}}{n!}\\
&\leq &C(n)\int_{ U} |F|_{\omega_{1}}^{\alpha }(\frac{\omega_{\epsilon_{1}}^{n}}{\omega_{\epsilon}^{n}})^{(\alpha -1)} \frac{\omega_{\epsilon}^{n}}{n!}\\
&\leq &C(n) (\int_{U} |F|_{\omega_{\epsilon_{1} }}^{2 } \frac{\omega_{\epsilon_{1}}^{n}}{n!})^{\frac{\alpha }{2}}
(\int_{U}  (\frac{\omega_{\epsilon_{1}}^{n}}{\omega_{\epsilon}^{n}})^{\frac{2\alpha -2}{2-\alpha }} \frac{\omega_{\epsilon_{1}}^{n}}{n!})^{\frac{2-\alpha }{2}}\\
&\leq & C(n) (Vol(U, \omega_{1} ))^{\frac{2-\alpha }{2}(1-\frac{1}{q})}(\int_{U} |F|_{\omega_{\epsilon_{1} }}^{2 } \frac{\omega_{\epsilon_{1}}^{n}}{n!})^{\frac{\alpha}{2}}
(\int_{U}  (\frac{\omega_{\epsilon_{1}}^{n}}{\omega_{\epsilon}^{n}})^{\frac{2\alpha -2}{2-\alpha }\cdot q} \frac{\omega_{\epsilon_{1}}^{n}}{n!})^{\frac{2-\alpha }{2q}},\\
\end{array}
\end{eqnarray}
where $q= \frac{1}{2}(\frac{2\alpha -2}{2-\alpha } +\frac{1}{k-1})\cdot (\frac{2\alpha -2}{2-\alpha })^{-1}$, and note that by the condition on $\alpha $ we have $\frac{2\alpha -2}{2-\alpha }\cdot q <\frac{1}{k-1}$.

Using (\ref{Int1}) in the previous lemma, we see that the result follows from (\ref{LLL1}), (\ref{LLL2}), (\ref{LLL3}) and (\ref{LLL4}).

\hfill $\Box$ \\

\medskip

{\bf Proposition 5.10. }  {\it Let $(\textbf{E} , \phi )$ be a
Higgs bundle on a smooth K\"ahler manifold $(M , \omega )$, and let $\textit{F}=\{F_{i}\}_{i=1}^{l}$ be a
filtration of $E$ by saturated subsheaves, where every $F_{i}$ is $\phi $ invariant. Let  $\pi :
\overline{M}\rightarrow M$ be a blow-up along a smooth complex manifold $\Sigma$ of complex co-dimension $k \geq 2$, $\overline{\textbf{E}}=\pi^{\ast }\textbf{E}$ be the pull-back bundle and $\overline{\phi }=\pi ^{\ast } \phi $. Let $\eta $ be a K\"ahler metric on $\overline{M}$, and  set a family of
K\"ahler metrics $\omega_{\epsilon }=\pi^{\ast
}\omega +\epsilon \eta $, and  $\overline{\textit{F}}=\{\overline{F}_{i}\}_{i=1}^{l}=\{Sat_{\overline{E}}(\pi^{\ast }F_{i})\}_{i=1}^{l}$  is a filtration of
$\overline{\textbf{E}}$ (not necessary given by subbundles).
 Suppose that for any $\tilde{\delta } >0$ and any $0<\epsilon < \epsilon ^{\ast}$, there is
 a smooth Hermitian metric $\overline{H}$ on
$\overline{\textbf{E}}$ such that
 \begin{eqnarray}\label{ca1}
\|\frac{\sqrt{-1}}{2\pi}\Lambda_{\omega_{\epsilon}
}(F_{(\overline{\partial }_{\overline{E}} ,
\overline{H})}+[\overline{\phi} , \overline{\phi}^{\ast \overline{H}
}])-\Psi(\overline{\textit{F}}, (\mu_{\epsilon , 1} , \cdot , \mu_{\epsilon , l}) ,
\overline{H})\|_{L^{2}(\overline{M}, \omega_{\epsilon})}\leq \tilde{\delta }.
\end{eqnarray}
Then for any $\delta >0$ and any $1\leq p < 1+\frac{1}{2k-1}$ there are $\epsilon_{1}>0$ and a smooth Hermitian metric $\overline{H}_{1}$ on $\overline{E}$ such that
  \begin{eqnarray}\label{a2}
\|\frac{\sqrt{-1}}{2\pi}\Lambda_{\omega_{\epsilon}
}(F_{(\overline{\partial }_{\overline{E}} ,
\overline{H}_{1})}+[\overline{\phi} , \overline{\phi}^{\ast \overline{H}_{1}
}])-\Psi(\overline{\textit{F}}, (\mu_{ 1} , \cdot , \mu_{l}) ,
\overline{H}_{1})\|_{L^{p}(\overline{M}, \omega_{\epsilon})}\leq \delta ,
\end{eqnarray}
for all $0<\epsilon \leq \epsilon_{1}$, where $ \mu_{i}$ is the $\omega$-slope of sheaf $F_{i}$.}

\medskip

{\bf Proof. } Let $\epsilon_{1} \in (0, \epsilon ^{\ast}) $, by the condition, we can choose a smooth metric $\overline{H}_{1}$ satisfies (\ref{ca1}) for $\epsilon_{1}$ and $\tilde{\delta}$ which will be chosen small enough later. For simplicity, we denote $\Theta_{1} =\frac{\sqrt{-1}}{2\pi}(F_{(\overline{\partial }_{\overline{E}} ,
\overline{H}_{1})}+[\overline{\phi} , \overline{\phi}^{\ast \overline{H}_{1}
}])$.  Then
 \begin{eqnarray}\label{ZZ1}
 \begin{array}{lll}
&&\|\frac{\sqrt{-1}}{2\pi}\Lambda_{\omega_{\epsilon}
}(F_{(\overline{\partial }_{\overline{E}} ,
\overline{H}_{1})}+[\overline{\phi} , \overline{\phi}^{\ast \overline{H}_{1}
}])-\Psi(\overline{\textit{F}}, (\mu_{ 1} , \cdots , \mu_{ l}) ,
\overline{H}_{1})\|_{L^{p}(\omega_{\epsilon})}\\
&\leq &\|\Lambda_{\omega_{\epsilon}
}\{\Theta_{1}-\frac{\omega_{\epsilon_{1}}}{n}\Psi(\overline{\textit{F}}, (\mu_{\epsilon_{1} , 1} , \cdots , \mu_{\epsilon_{1} , l}) ,
\overline{H}_{1})\}\|_{L^{p}(\omega_{\epsilon})}\\
&& +\| \frac{1}{n}\Lambda_{\omega_{\epsilon}
}(\omega_{\epsilon _{1}}-\omega_{\epsilon })\Psi(\overline{\textit{F}}, (\mu_{\epsilon_{1} , 1} , \cdots , \mu_{\epsilon_{1} , l}) ,
\overline{H}_{1}) \|_{L^{p}(\omega_{\epsilon})}\\
&& +\|\Psi(\overline{\textit{F}}, (\mu_{ 1} , \cdots , \mu_{ l}) ,
\overline{H}_{1})- \Psi(\overline{\textit{F}}, (\mu_{\epsilon_{1} , 1} , \cdots , \mu_{\epsilon_{1} , l}) ,
\overline{H}_{1}) \|_{L^{p}(\omega_{\epsilon})}.
\end{array}
\end{eqnarray}

For simplicity, we set $\Theta_{2}=\Theta_{1}-\frac{\omega_{\epsilon_{1}}}{n}\Psi(\overline{\textit{F}}, (\mu_{\epsilon_{1} , 1} , \cdots , \mu_{\epsilon_{1} , l}) ,
\overline{H}_{1})$. From the equality
\begin{eqnarray}
\begin{array}{lll}
&&\int_{\overline{M}} (2C_{2}(\overline{E})-\frac{r-1}{r}C_{1}(\overline{E})\wedge C_{1}(\overline{E}))\frac{\omega_{\epsilon_{1}}^{n-2}}{(n-2)!}\\
&=&\int_{M}|\Theta_{1}|_{\overline{H}_{1}}^{2}-|\Lambda_{\omega_{\epsilon_{1}}} \Theta_{1}|_{\overline{H}}^{2} \frac{\omega_{\epsilon_{1}}^{n}}{n!},\\
\end{array}
\end{eqnarray}
we know that $\|\Theta_{2}\|_{L^{2 }(\overline{M}, \omega_{\epsilon_{1}})}$ is bounded uniformly.
By (\ref{Int2}), we have
\begin{eqnarray}
\begin{array}{lll}
&&\|\Lambda_{\omega_{\epsilon }}\Theta_{2}\|_{L^{p }(\overline{M}, \omega_{\epsilon})} \leq  |\epsilon_{1}-\epsilon |C(U)\|\Theta_{2}\|_{L^{p }(\overline{M}, \omega_{\epsilon_{1}})}\\
&+& C(\|\Lambda_{\omega_{\epsilon_1 }}\Theta_{2}\|_{L^{2 }(\overline{M}, \omega_{\epsilon_1})} +(Vol(U, \omega_{1}))^{C(\alpha , k)}\|\Theta_{2}\|_{L^{2 }(\overline{M}, \omega_{\epsilon_{1}})}).\\
\end{array}
\end{eqnarray}
We may choose $U$ such that $Vol(U, \omega_{1})$ small enough first, and then $\tilde{\delta }$ and $\epsilon_{1}$ both sufficiently small so that
\begin{eqnarray}
\|\Lambda_{\omega_{\epsilon }}\Theta_{2}\|_{L^{p }(\overline{M}, \omega_{\epsilon})} \leq \frac{\delta }{3}.
\end{eqnarray}

 By (\ref{Int1}) in Lemma 5.8., it is not difficult to see that $\|\Lambda_{\omega_{\epsilon}}\eta \|_{L^{p}(\overline{M}, \omega_{\epsilon})}$ is uniformly bounded.
On the other hand, since $\mu_{\epsilon , i} \rightarrow \mu_{i}$ as $\epsilon \rightarrow 0$,  we may choose $\epsilon_{1}$ small enough so that the second and third terms in (\ref{ZZ1}) are both smaller than $\frac{\delta }{3}$, so (\ref{a2}) follows.

\hfill $\Box$ \\

\medskip

\medskip

{\bf Proposition 5.11. } {\it Let $(\textbf{E} , \phi )$ be a
Higgs bundle on a smooth K\"ahler manifold $(M , \omega )$, and
 let $\textit{F}=\{F_{i}\}_{i=1}^{l}$ be a
filtration of $E$ by saturated subsheaves, where every $F_{i}$ is $\phi $ invariant. Let  $\pi :
\overline{M}\rightarrow M$ be a blow-up along a smooth complex manifold $\Sigma$ of complex co-dimension $k \geq 2$, $\overline{\textbf{E}}=\pi^{\ast }\textbf{E}$ be the pull-back bundle and $\overline{\phi }=\pi ^{\ast } \phi $. Let $\eta $ be a K\"ahler metric on $\overline{M}$, and  set a family of
K\"ahler metrics $\omega_{\epsilon }=\pi^{\ast
}\omega +\epsilon \eta $, and  $\overline{\textit{F}}=\{\overline{F}_{i}\}_{i=1}^{l}=\{Sat_{\overline{E}}(\pi^{\ast }F_{i})\}_{i=1}^{l}$  is a filtration of
$\overline{\textbf{E}}$ (not necessary given by subbundles). Suppose that for any $\delta >0$ and any $1\leq p < 1+\frac{1}{2k-1}$ there is a smooth metric $\overline{H}_{1}$ on $\overline{E}$ and $\epsilon_{1}>0$ such that (\ref{a2}) hold for all $0<\epsilon \leq \epsilon_{1}$. Then for any $\delta' >0$ and any $1\leq p < 1+\frac{1}{2k-1}$ there is a smooth metric $H$ on $E$ such that
  \begin{eqnarray}\label{a3}
\|\frac{\sqrt{-1}}{2\pi }\Lambda_{\omega
}(F_{(\overline{\partial }_{E} ,
H)}+[\phi , \phi ^{\ast H
}])-\Psi(\textit{F}, (\mu_{ 1} , \cdots , \mu_{l}) ,
H)\|_{L^{p}(M, \omega)}\leq \delta' ,
\end{eqnarray}
where $ \mu_{i}$ is the $\omega$-slope of sheaf $F_{i}$.}

\medskip

{\bf Proof } We use a cut-off argument to get the smooth metric on bundle $E$. Since $\Sigma $ is a smooth complex submanifold, the open set $\{(x, \nu ) \in \textbf{N}_{\Sigma } | |\nu | <R \}$ in the normal bundle $\textbf{N}_{\Sigma }$ of $\Sigma $, is diffeomorphic to an open neighborhood $U_{R}$ of $\Sigma $ for $R$ sufficiently small. For any small $R$, we may choose a smooth cut-off function $\psi_{R}$ which supported in $U_{R}$ and identically $1$ on $U_{\frac{R}{2}}$, $0\leq \psi_{R} \leq 1$, and furthermore  $|\partial \psi_{R}|_{\omega}^{2}+|\partial \overline{\partial } \psi_{R}|_{\omega }\leq C R^{-2}$, where $C$ is a positive constant
independent of $R$. Let $H_{D}$ be a smooth Hermitian metric on bundle $E$, and $\overline{H}_{1}$ be the metric on $\overline{E}$ such that (\ref{a2}) holds for all $0<\epsilon \leq \epsilon_{1}$  where $\delta \leq \frac{\delta '}{4} $. Note that $E$ is isomorphic to $\overline{E}$ outsides $\Sigma $, we can define
\begin{eqnarray}
H_{R}=(1-\psi_{R})\overline{H}_{1}+\psi_{R} H_{D}
\end{eqnarray}
on bundle $E$, and $\overline{H}_{R}=\pi^{\ast}H_{R}$ on bundle $\overline{E}$.

As above, we denote $\Theta (\overline{H}_{R})=\frac{\sqrt{-1}}{2\pi}(F_{(\overline{\partial }_{\overline{E}} ,
\overline{H}_{R})}+[\overline{\phi} , \overline{\phi}^{\ast \overline{H}_{R}
}]) $. We have
\begin{eqnarray}\label{ZZZ1}
\begin{array}{lll}
&&\int_{\overline{M}} |\Lambda_{\omega_{\epsilon}
}\Theta (\overline{H}_{R}) -\Psi(\overline{\textit{F}}, (\mu_{ 1} , \cdots , \mu_{l}) ,
\overline{H}_{R})|_{\overline{H}_{R}}^{p} \frac{\omega_{\epsilon}^{n}}{n!} \\
&\leq & \int_{\pi^{-1}(U_{\frac{R}{2}})} |\Lambda_{\omega_{\epsilon}
}\Theta (\pi^{\ast}{H}_{D}) -\Psi(\overline{\textit{F}}, (\mu_{ 1} , \cdots , \mu_{l}) ,
\pi^{\ast}{H}_{D})|_{\pi^{\ast}{H}_{D}}^{p} \frac{\omega_{\epsilon}^{n}}{n!} \\
&&+\int_{M\setminus \pi^{-1}(U_{R})} |\Lambda_{\omega_{\epsilon}
}\Theta (\overline{H}_{1}) -\Psi(\overline{\textit{F}}, (\mu_{ 1} , \cdots , \mu_{l}) ,
\overline{H}_{1})|_{\overline{H}_{1}}^{p} \frac{\omega_{\epsilon}^{n}}{n!} \\
&&+C(p)\int_{\pi^{-1}(U_{R}\setminus U_{\frac{R}{2}})} |\Lambda_{\omega_{\epsilon}
}(\Theta (\overline{H}_{R}) -\Theta (\overline{H}_{1}))|_{\overline{H}_{R}}^{p} \frac{\omega_{\epsilon}^{n}}{n!} \\
&&+C(p) \int_{\pi^{-1}(U_{R}\setminus U_{\frac{R}{2}})}|\Psi(\overline{\textit{F}}, (\mu_{ 1} , \cdots , \mu_{l}) ,
\overline{H}_{R}) -\Lambda_{\omega_{\epsilon}}\Theta (\overline{H}_{1}))|_{\overline{H}_{R}}^{p} \frac{\omega_{\epsilon}^{n}}{n!}.\\
\end{array}
\end{eqnarray}

By the definition of $\Lambda_{\omega_{\epsilon }}$, we have
\begin{eqnarray}\label{ZZZ2}
\begin{array}{lll}
&&\int_{\pi^{-1}(U_{\frac{R}{2}})} |\Lambda_{\omega_{\epsilon}
}\Theta (\pi^{\ast}{H}_{D}) |_{\pi^{\ast}{H}_{D}}^{p} \frac{\omega_{\epsilon}^{n}}{n!}\\
&=&\int_{\pi^{-1}(U_{\frac{R}{2}})} |n\frac{\Theta (\pi^{\ast}{H}_{D})\wedge \omega_{\epsilon}^{n-1}}{\omega_{\epsilon}^{n}} |_{\pi^{\ast}{H}_{D}}^{p} \frac{\omega_{\epsilon}^{n}}{n!}\\
&\leq & C_{0}\int_{\pi^{-1}(U_{\frac{R}{2}})} |\Theta (\pi^{\ast}{H}_{D})|_{\pi^{\ast}{H}_{D}, \eta }^{p}(\frac{\eta^{n}}{\omega_{\epsilon}^{n}})^{p-1}\frac{\eta^{n}}{n!}\\
&\leq & C_{1} \int_{\pi^{-1}(U_{\frac{R}{2}})} (\frac{\eta^{n}}{\omega_{\epsilon}^{n}})^{p-1}\frac{\eta^{n}}{n!},
\end{array}
\end{eqnarray}
where $C_{1}$
is a constant independent of $\epsilon $ and $R$.

By directly calculation, we have
\begin{eqnarray}\label{TTT2}
\begin{array}{lll}
&&F_{\overline{H}_{R}}-F_{\overline{H}_{1}}= F_{\overline{H}_{R}}-F_{\overline{H}_{D}}+F_{\overline{H}_{D}}-F_{\overline{H}_{1}}\\
&=&\overline{\partial }_{\overline{E}} (h_{R}^{-1}\partial_{\overline{H}_{D}}h_{R})-\overline{\partial }_{\overline{E}} (h_{1}^{-1}\partial_{\overline{H}_{D}}h_{1})\\
&=& \overline{\partial }_{\overline{E}} (-h_{R}^{-1}h_{1}\partial \psi_{R} +(1-\psi_{R}) h_{R}^{-1}\partial_{\overline{H}_{D}}h_{1} +\partial \psi_{R} h_{R}^{-1})\\
&& -\overline{\partial }_{\overline{E}} (h_{1}^{-1}\partial_{\overline{H}_{D}}h_{1})\\
&=&\partial\overline{\partial } \psi_{R}( h_{R}^{-1}h_{1}-h_{R}^{-1}) -\overline{\partial }_{\overline{E}} (h_{R}^{-1}h_{1})\wedge \partial \psi_{R} -\partial \psi_{R}\wedge \overline{\partial }_{\overline{E}} h_{R}^{-1}\\
&& -\overline{\partial }\psi_{R} \wedge h_{R}^{-1}\partial_{H_{D}}h_{1}
 +(1-\psi_{R})\overline{\partial }_{\overline{E}} (h_{R}^{-1}h_{1})\wedge h_{1}^{-1}\partial_{\overline{H}_{D}}h_{1} \\ &&+((1-\psi_{R})h_{R}^{-1}h_{1}-Id)\overline{\partial }_{\overline{E}} (h_{1}^{-1}\partial_{\overline{H}_{D}}h_{1}),\\
&=& \partial\overline{\partial } \psi_{R}( h_{R}^{-1}h_{1}-h_{R}^{-1})
 -\overline{\partial }\psi_{R} \wedge h_{R}^{-1}\partial_{H_{D}}h_{1}+\partial \psi_{R}\wedge h_{R}^{-1}\overline{\partial }_{\overline{E}} h_{1}\\
 && +\partial \psi_{R} \wedge \overline{\partial}\psi_{R}(h_{R}^{-1}h_{1}h_{R}^{-1}-Id)(h_{1}-Id)\\
&& +(1-\psi_{R})(\overline{\partial }_{\overline{E}} h_{1}\wedge h_{1}^{-1}\partial_{\overline{H}_{D}}h_{1}-\partial \psi_{R} \wedge h_{R}^{-1}\overline{\partial }_{\overline{E}} h_{1} h_{R}^{-1}(h_{1}-Id)) \\
&&+ (1-\psi_{R})(h_{R}^{-1}h_{1}h_{R}^{-1}-Id)\overline{\partial }\psi_{R}\wedge \partial_{\overline{H}_{D}}h_{1} \\
&& -(1-\psi_{R})^{2}h_{R}^{-1}\overline{\partial }_{\overline{E}} h_{1}\wedge h_{R}^{-1}\partial_{\overline{H}_{D}}h_{1}\\
&&+((1-\psi_{R})h_{R}^{-1}h_{1}-Id)\overline{\partial }_{\overline{E}} (h_{1}^{-1}\partial_{\overline{H}_{D}}h_{1}),\\
\end{array}
\end{eqnarray}
where $\overline{H}_{D}=\pi^{\ast }H_{D}$, $h_{R}=\overline{H}_{D}^{-1}\overline{H}_{R}$ and $h_{1}=\overline{H}_{D}^{-1}\overline{H}_{1}$. Since the metrics $\overline{H}_{1}$ and $\overline{H}_{D}$ are fixed, so we have $\overline{C} ^{-1}Id \leq h_{R}\leq \overline{C} Id$, where constant $\overline{C}$ depends only on $\overline{H}_{1}$ and $\overline{H}_{D}$.  From the equality (\ref{TTT2}), we have
\begin{eqnarray}
|\Lambda_{\omega_{\epsilon}} (F_{\overline{H}_{R}}-F_{\overline{H}_{1}})|_{\overline{H}_{R}}\leq C_{2}(|\partial \overline{\partial }\psi_{R}|_{\omega_{\epsilon}}+|\partial \psi_{R}|_{\omega_{\epsilon}}^{2})+C_{3}\frac{\eta ^{n}}{\omega_{\epsilon}^{n}},
\end{eqnarray}
where $C_{2}$ and $C_{3}$ are constants independent of $\epsilon$ and $R$. On the other hand , we have
\begin{eqnarray}
|\Lambda_{\omega_{\epsilon }}([\overline{\phi} , \overline{\phi}^{\ast \overline{H}_{R}
}]-[\overline{\phi} , \overline{\phi}^{\ast \overline{H}_{1}
}])|_{\overline{H}_{R}}\leq C_{4} \frac{\eta ^{n}}{\omega_{\epsilon}^{n}},
\end{eqnarray}
where $C_{4}$ is a constant which  may depend on $\phi $, $\eta $, $\overline{H}_{D}$ and $\overline{H}_{1}$, but it is independent of $\epsilon$ and $R$.
Thus
\begin{eqnarray}\label{ZZZ3}
\begin{array}{lll}
&&\int_{\pi^{-1}(U_{R}\setminus U_{\frac{R}{2}})} |\Lambda_{\omega_{\epsilon}
}(\Theta (\overline{H}_{R}) -\Theta (\overline{H}_{1}))|_{\overline{H}_{R}}^{p} \frac{\omega_{\epsilon}^{n}}{n!}\\
&\leq & C_{5}R^{-2p} \int_{\pi^{-1}(U_{R}\setminus U_{\frac{R}{2}})} 1 \frac{\omega_{\epsilon}^{n}}{n!} +C_{6} \int_{\pi^{-1}(U_{R}\setminus U_{\frac{R}{2}})} ( \frac{\eta^{n}}{\omega_{\epsilon}^{n}})^{p-1}\frac{\eta ^{n} }{n!},\\
\end{array}
\end{eqnarray}
where $C_{4}$ and $C_{5}$ are constants independent of $\epsilon$ and $R$. Similarly, we have
\begin{eqnarray}\label{TTT1}
\begin{array}{lll}
&& \int_{\pi^{-1}(U_{R}\setminus U_{\frac{R}{2}})}|\Lambda_{\omega_{\epsilon}}\Theta (\overline{H}_{1}))|_{\overline{H}_{R}}^{p} \frac{\omega_{\epsilon}^{n}}{n!}\\
&\leq &C(n, p)\int_{\pi^{-1}(U_{R}\setminus U_{\frac{R}{2}})}|\Theta (\overline{H}_{1})|_{\overline{H}_{R}, \eta }^{p}(\frac{\eta^{n}}{\omega_{\epsilon}^{n}})^{p-1}\frac{\eta^{n}}{n!}\\
&\leq &C_{7}\int_{\pi^{-1}(U_{R}\setminus U_{\frac{R}{2}})}(\frac{\eta^{n}}{\omega_{\epsilon}^{n}})^{p-1}\frac{\eta^{n}}{n!},\\
\end{array}
\end{eqnarray}
where $C_{7}$ is a constant also independent of $\epsilon$ and $R$.

By (\ref{Int1}), we have $\int_{\overline{M}\setminus \pi^{-1}(\Sigma ) }( \frac{\eta ^{n}}{ (\pi^{\ast}\omega )^{n}})^{\gamma } \frac{\eta ^{n}}{n!}\leq C^{\ast}$. By the relation $\pi^{\ast}\omega < \omega_{\epsilon }<\omega_{1}$, it is easy to see that
\begin{eqnarray}\label{ZZZ4}
\begin{array}{lll}
\int_{\pi^{-1}(U_{R})} ( \frac{\eta^{n}}{\omega_{\epsilon}^{n}})^{p-1}\frac{\eta ^{n} }{n!} &\rightarrow &0,\\
\int_{\pi^{-1}(U_{R})}|\Psi(\overline{\textit{F}}, (\mu_{ 1} , \cdots , \mu_{l}) ,
\overline{H}_{D})|_{\overline{H}_{R}}^{p}  \frac{\omega_{\epsilon}^{n}}{n!} &\rightarrow &0 ,\\
\int_{\pi^{-1}(U_{R})}|\Psi(\overline{\textit{F}}, (\mu_{ 1} , \cdots , \mu_{l}) ,
\overline{H}_{R})|_{\overline{H}_{R}}^{p} \frac{\omega_{\epsilon}^{n}}{n!} &\rightarrow &0
\end{array}
\end{eqnarray}
as $R\rightarrow 0$, uniformly in $\epsilon $.

By (\ref{ZZZ1}), (\ref{ZZZ2}), (\ref{ZZZ3}), (\ref{TTT1}) and (\ref{ZZZ4}), and choosing  $R_{0}$ sufficiently small, we have
\begin{eqnarray}
\begin{array}{lll}
&&\int_{\overline{M}} |\Lambda_{\omega_{\epsilon}
}\Theta (\overline{H}_{R}) -\Psi(\overline{\textit{F}}, (\mu_{ 1} , \cdots , \mu_{l}) ,
\overline{H}_{R})|_{\overline{H}_{R}}^{p} \frac{\omega_{\epsilon}^{n}}{n!} \\
&\leq & \frac{\delta '}{2} + C_{5}R^{-2p} \int_{\pi^{-1}(U_{R}-U_{\frac{R}{2}})} 1 \frac{\omega_{\epsilon}^{n}}{n!},
\end{array}
\end{eqnarray}
for all $0<\epsilon \leq \epsilon_{1}$ and $0< R\leq R_{0}$.
Let $\epsilon \rightarrow 0$, we have
\begin{eqnarray}
\begin{array}{lll}
&&\int_{M} |\Lambda_{\omega
}\Theta (H_{R}) -\Psi(\textit{F}, (\mu_{ 1} , \cdots , \mu_{l}) ,
H_{R})|_{H_{R}}^{p} \frac{\omega ^{n}}{n!} \\
&\leq & \frac{\delta '}{2} + C_{5}R^{-2p} \int_{U_{R}} 1 \frac{\omega^{n}}{n!}.
\end{array}
\end{eqnarray}
Since $\Sigma $ has Hausdorff dimension at most $2n-2k$, it is easy to see that $Vol (U_{R}, \omega )\leq C R^{2k}$ for some uniform constant $C$. By the assumption of $p$, we know that $2k -2p>0$, choosing $R$ small enough, then (\ref{a3}) follows.

\hfill $\Box$ \\

\medskip

{\bf Theorem 5.12. } {\it Let $(E, A_{0} , \phi_{0} )$ be a
Higgs bundle on a smooth K\"ahler manifold $(M , \omega )$, and $(A_{t} , \phi_{t})$ be  the smooth solution of the Yang-Mills-Higgs flow (\ref{YMHH}) on the Hermitian vector bundle $(E, H_{0})$ with initial data $(A_{0}, \phi_{0})\in \textbf{B}_{(E, H_{0})}$.
Suppose that for any $\delta ' >0$ and any $1\leq p < p_{0}$ there is a smooth metric $H$ on $E$  such that (\ref{a3}) holds, where $\vec{\mu }_{0}$ is the Harder-Narasimhan type of $(E, A_{0} , \phi_{0} )$. Let $(A_{\infty }, \phi_{\infty})$ be an Uhlenbeck limit of $(A_{t}, \phi_{t})$, and $(E_{\infty}, H_{\infty})$ be the corresponding Hermitian vector bundle defined away from $\Sigma_{an}$.  Then
\begin{eqnarray}\label{CC1}
HYM_{\alpha , N} (A_{\infty }, \phi_{\infty})=\lim_{t\rightarrow \infty}HYM_{\alpha , N} (A_{t }, \phi_{t })=HYM_{\alpha , N} (\vec{\mu}_{0})
\end{eqnarray}
for all $1\leq \alpha < p_{0}$ and all $N\in R$; the HN type of the Higgs sheaf $(E_{\infty}, A_{\infty}, \phi_{\infty})$ is the same as $(E, A_{0} , \phi_{0})$.}

\medskip

{\bf Proof } Firstly, since the norm $(\int_{M}
\varphi_{\alpha }(\verb"a") dvol )^{\frac{1}{\alpha }}$ is equivalent to the
$L^{\alpha }$ norm on $\verb"u"(E)$,  we have,
\begin{eqnarray*}
\begin{array}{lll}
&& |(HYM_{\alpha , N}((\overline{\partial }_{A_{0}}, H), \phi_{0}))^{\frac{1}{\alpha }}-(HYM_{\alpha , N}(\vec{\mu }_{0}))^{^{\frac{1}{\alpha }}}|\\
&\leq &(\int_{M} |(\varphi_{\alpha }(\sqrt{-1}(\Theta ((\overline{\partial }_{A_{0}}, H), \phi_{0})+N Id_{E})))^{\frac{1}{\alpha}}-(\varphi_{\alpha }(\vec{\mu}+N))^{\frac{1}{\alpha}}|^{\alpha }\frac{\omega ^{n}}{n!})^{\frac{1}{\alpha }}\\
&\leq & (\int_{M} \varphi_{\alpha }(\sqrt{-1}(\Theta ((\overline{\partial }_{A_{0}}, H), \phi_{0}))- \Psi(\textit{F}, (\mu_{ 1} , \cdots , \mu_{l}) ,
H))\frac{\omega ^{n}}{n!} )^{\frac{1}{\alpha}}\\
&\leq &
C(\alpha )\|\Theta ((\overline{\partial }_{A_{0}}, H), \phi_{0})-\Psi(\textit{F}, (\mu_{ 1} , \cdots , \mu_{l}) ,
H)\|_{L^{\alpha }(M, \omega)}.
\end{array}
\end{eqnarray*}
By the above inequality and the condition  (\ref{a3}), we see for any $\delta >0$ and any $1\leq \alpha  < p_{0}$ there is $H$ such that
\begin{eqnarray}\label{BB2}
HYM_{\alpha , N}((\overline{\partial}_{A_{0}}, H), \phi_{0} )\leq  HYM_{\alpha , N}(\vec{\mu }_{0}) +\delta .
\end{eqnarray}

\medskip

 For fixed $\alpha $, $1\leq \alpha \leq \alpha_{0}$,
and fixed $N$, since the image of the degree map on line bundles is discrete, we can  define $\delta_{0}>0$ such that
\begin{eqnarray}
\begin{array}{lll}
&&2\delta_{0}+HYM_{\alpha , N}(\vec{\mu}_{0})\\ &=&min\{HYM_{\alpha ,
N}(\vec{\mu}): HYM_{\alpha , N}(\vec{\mu})>HYM_{\alpha ,
N}(\vec{\mu}_{0})\},\\
\end{array}
\end{eqnarray}
where $\vec{\mu }$ runs over all possible HN types of Higgs sheaves
on $M$ with the rank of $E$.

\medskip

Let $H$ be a Hermitian metric on the complex bundle $E$, and
$(A_{t}^{H} , \phi_{t}^{H})$ be the solution to the Yang-Mills-Higgs
flow (\ref{YMHH}) on the Hermitian vector bundle $(E , H)$ with initial
pair $(A_{0}^{H} , \phi_{0})\in \textbf{B}_{(E, H)}$ where
$A_{0}^{H}=(\overline{\partial }_{A_{0}}, H)$. Let $(A_{\infty }^{H} , \phi_{\infty })$ be an Uhlenbeck limit along
the flow (\ref{YMHH}).

Assume  that the $H$ satisfies:
\begin{eqnarray}\label{5.12}
HYM_{\alpha , N}((\overline{\partial}_{A_{0}}, H), \phi_{0} ))\leq  HYM_{\alpha , N}(\vec{\mu }_{0}) +\delta_{0} .
\end{eqnarray}
 By Prop. 5.3 and
Prop. 5.5, we obtain:
\begin{eqnarray*}
HYM_{\alpha , N}(\vec{\mu}_{0})\leq HYM_{\alpha , N}(A_{\infty }^{H} ,
\phi_{\infty})\leq
HYM_{\alpha , N}(\vec{\mu}_{0})+\delta_{0} .
\end{eqnarray*}
Hence, we must have $HYM_{\alpha , N}(A_{\infty}^{H} , \phi_{\infty}) =
HYM_{\alpha , N}(\vec{\mu}_{0})$. This shows that the result holds
if the metric $H_{0}$ satisfies (\ref{5.12}).

We are going to prove that for any metric $H$, for any fixed $\delta $ there is
$T\geq 0$ such that:
\begin{eqnarray}\label{ineq1}
HYM_{\alpha , N}(A_{t}^{H} , \phi_{t}^{H}) < HYM_{\alpha ,
N}(\vec{\mu}_{0})+\delta ,
\end{eqnarray}
for all $t\geq T$. Without loss of generality, we can assume
$0<\delta \leq \frac{\delta_{0}}{2}$.

Let us denote by $\textbf{H}_{\delta }$  the set of smooth Hermitian
metrics on $E$ with the property that the above inequality (\ref{ineq1}) holds for
some $T$. From (\ref{BB2}) and the discussion above, we know
$\textbf{H}_{\delta }$ is nonempty. In \cite{LZ1} (Proposition 2.1'), we have proved that the continuous dependence of the Donaldson's flow (\ref{D1}) on initial conditions. Following the argument in \cite{DW1} (Lemma 4.3), see also Theorem 5.13 in \cite{LZ1}, we can show that $\textbf{H}_{\delta }$ is closed and also open. The proof is exactly the same as that in \cite{LZ1} (Theorem 5.13), we omit it.  Since the
space of smooth metrics is connected, we conclude that every
metric is in $\textbf{H}_{\delta }$. Then, we have
$lim_{t\rightarrow \infty }HYM_{\alpha , N}(A_{t}^{H} ,
\phi_{t}^{H}) =HYM_{\alpha , N}(\vec{\mu}_{0})$ for any metric
$H$.

\medskip

Let $\vec{\lambda}_{\infty}$ be the HN type of $(E_{\infty}, A_{\infty}, \phi_{\infty})$, by (\ref{CC1}) and Proposition 5.5, we have $\varphi_{\alpha }(\vec{\mu}_{0}+N)=\varphi_{\alpha }(\vec{\lambda }_{\infty }+N )$ for all $1\leq \alpha < p_{0}$ and all $N$. We may choose $N$ sufficiently large so that every elements in the vectors $\vec{\mu}_{0}+N$ and $\vec{\lambda }_{\infty }+N$ are positive. Then, by Lemma 5.2, we have $\vec{\mu}_{0}+N = \vec{\lambda }_{\infty }+N$, and so $\vec{\mu}_{0} = \vec{\lambda }_{\infty }$.

\hfill $\Box$ \\

\medskip

Let the following filtrations of saturated sheaves
\begin{eqnarray}
0=E_{0}\subset E_{1}\subset \cdots \subset E_{l}=(E, \overline{\partial}_{A_{0}})
\end{eqnarray}
be a HN filtration of the Higgs bundle $(E, \overline{\partial}_{A_{0}}, \phi_{0})$. The action of $g_{j}$ produces a sequence of HN filtration
\begin{eqnarray}
0=E_{0}^{(j)}\subset E_{1}^{(j)}\subset \cdots \subset E_{l}^{(j)}=(E, \overline{\partial}_{A_{j}}),
\end{eqnarray}
where $E_{\alpha}^{(j)}=g_{j}(E_{\alpha})$, $\alpha =1, \cdots , l$. Let $\pi _{\alpha}^{(j)}$ be the orthogonal projection onto $E_{\alpha}^{(j)}$, then we have
 $\pi_{\alpha }^{(j)} \in L_{1}^{2}$ and satisfies the conditions in (\ref{WHC}). Using the above Theorem 5.12, by the same argument in \cite{DW1} (Proposition 4.5), we have the following lemma.

\medskip

{\bf Lemma 5.13. } {\it Let $(E, A_{0} , \phi_{0} )$ be a
Higgs bundle on a smooth K\"ahler manifold $(M , \omega )$, and satisfy the same assumptions as that in Theorem 5.12.

(1) Let $\{\pi_{\alpha}^{\infty}\}$ be the HN filtration of the reflexive Higgs sheaf $(E_{\infty}, A_{\infty}, \phi_{\infty})$, then there is a subsequence of HN filtration $\{ \pi_{\alpha }^{(j)} \}$ converges to a filtration $\{ \pi_{\alpha }^{\infty } \}$  strongly in $L^{p}\cap L_{1, loc}^{2}$ off $\Sigma_{an}$.

(2) Assume the Higgs bundle $(E, A_{0}, \phi_{0})$ is semi-stable and $\{E_{\alpha }\}$ is the Seshadri filtration of $(E, A_{0}, \phi_{0})$, then, after passing to a subsequence, $\{ \pi_{\alpha }^{(j)} \}$ converges to a filtration $\{ \pi_{\alpha }^{\infty } \}$  strongly in $L^{p}\cap L_{1, loc}^{2}$ off $\Sigma_{an}$
, the rank and degree of $ \pi_{\alpha }^{\infty } $ is equal to the rank and degree of $ \pi_{\alpha }^{j } $ for all $\alpha $ and $j$.

}

\medskip

{\bf Proof } By the formula (\ref{5.4}) in Lemma 5.4, we have
\begin{eqnarray}\label{GC2}
\begin{array}{lll}
&&  deg (E_{\alpha }^{(j)})+\frac{1}{2\pi }\int_{M} |\overline{\partial}_{A_{j}} \pi_{\alpha }^{(j)}|_{H_{0}}^{2}+|[\phi_{j}, \pi_{\alpha }^{(j)}]|^{2} \frac{\omega^{n}}{n!}\\
&\leq & \sum_{i \leq rank E_{\alpha}}\lambda_{i} +\|\theta (A_{j}, \phi_{j})-\theta (A_{\infty}, \phi_{\infty})\|_{L^{1}}.
\end{array}
\end{eqnarray}
By theorem 5.12, we have $\vec{\mu}_{0} = \vec{\lambda }_{\infty }$. Since $E_{\alpha }^{(j)}$ is a Higgs sheaf of the Higgs bundle $(E_{\infty}, \overline{\partial }_{A_{\infty}}, \phi_{\infty})$ and $\mu (E_{\alpha }^{(j)})=\mu (E_{\alpha })=\mu (E)$, we have
\begin{eqnarray}\label{pi1}
\begin{array}{lll}
&& \int_{M}|\overline{\partial}_{A_{j}} \pi_{\alpha }^{(j)}|_{H_{0}}^{2} +|[\phi_{j}, \pi_{\alpha }^{(j)}]|^{2}\frac{\omega^{n}}{n!} \\
&\leq & \|\theta (A_{j}, \phi_{j})-\theta (A_{\infty}, \phi_{\infty})\|_{L^{1}} \rightarrow  0 \\
\end{array}
\end{eqnarray}
as $j \rightarrow \infty $, where we have the property that $\theta (A_{j}, \phi_{j})\rightarrow \theta (A_{\infty}, \phi_{\infty})$ strongly in $L^{p}$ for all $p$. After perhaps passing to a subsequence, we have $ \pi_{\alpha }^{(j)}\rightarrow  \tilde{\pi }_{\alpha }^{\infty}$ weakly in $L_{1}^{2}$, for some $L_{1}^{2}$ projection $\tilde{\pi}_{\alpha }^{\infty}$. Since $\pi_{\alpha}^{(j)}$ is uniformly bounded, we see $ \pi_{\alpha }^{(j)}\rightarrow  \tilde{\pi}_{\alpha }^{\infty}$ strongly in $L^{p}$ for all $p$. Recall that $(A_{j}, \phi_{j}) \rightarrow (A_{\infty}, \phi_{\infty})$ in $C^{\infty}_{loc}$ topology on $M\setminus \Sigma_{an}$ and write
\begin{eqnarray}
\overline{\partial}_{A_{\infty}} \pi_{\alpha }^{(j)}=\overline{\partial}_{A_{j}} \pi_{\alpha }^{(j)}+(A_{\infty}^{0,1}-A_{j}^{0,1 })\circ \pi_{\alpha }^{(j)}- \pi_{\alpha }^{(j)}\circ (A_{\infty}^{0,1}-A_{j}^{0,1 }),
\end{eqnarray}
 then as in the proof of Lemma 4.5 in \cite{DW1}, we conclude from (\ref{pi1}) that
$\overline{\partial}_{A_{\infty}}\tilde{\pi}_{\alpha }^{\infty}=0$, and $ \pi_{\alpha }^{(j)}\rightarrow  \tilde{\pi}_{\alpha }^{\infty}$ strongly in $L^{p}\cap L_{1, loc}^{2}$ off $\Sigma_{an}$. On the other hand, it is easy to check that $[\phi_{\infty}, \tilde{\pi}_{\alpha }^{\infty}]=0$ and $(\tilde{\pi}_{\alpha }^{\infty})^{2}=\tilde{\pi}_{\alpha }^{\infty}=(\tilde{\pi}_{\alpha }^{\infty})^{\ast}$. By Proposition 2.12., we know that $\tilde{\pi}_{\alpha }^{\infty}$ determines a $\phi_{\infty}$-invariant coherent subsheaf  $\tilde{E}_{\alpha }^{\infty}$ of $(E_{\infty}, \overline{\partial}_{A_{\infty}})$. Furthermore, it is clear $rank (\tilde{E}_{\alpha }^{\infty})=rank (E_{\alpha})=rank (\pi_{\alpha }^{\infty})$. Using (\ref{pi1}), we have
\begin{eqnarray}
\begin{array}{lll}
deg (\tilde{E}_{\alpha }^{\infty})&=& \int_{M}tr (\sqrt{-1}\theta (A_{\infty}, \phi_{\infty})\tilde{\pi}_{\alpha }^{\infty})\frac{\omega^{n}}{n!}\\
&=& \lim_{j\rightarrow \infty } \int_{M}tr (\sqrt{-1}\theta (A_{j}, \phi_{j})\tilde{\pi}_{\alpha }^{\infty})\frac{\omega^{n}}{n!}\\
&=& deg (E_{\alpha }^{(j)}) +\lim_{j\rightarrow \infty} \int_{M}|\overline{\partial}_{A_{j}} \pi_{\alpha }^{(j)}|_{H_{0}}^{2} +|[\phi_{j}, \pi_{\alpha }^{(j)}]|^{2}\frac{\omega^{n}}{n!}\\
&=& deg (E_{\alpha })=deg (\pi_{\alpha}^{\infty}).\\
\end{array}
\end{eqnarray}
So, the rank and degree of $ \tilde{\pi}_{\alpha }^{\infty } $ is equal to the rank and degree of
$\pi_{\alpha }^{\infty }$ for all $\alpha$. By the uniqueness of the maximal destabilizing subsheaf $\pi_{1}^{\infty}$ in the HN filtration of Higgs sheaf $(E_{\infty}, \overline{\partial}_{A_{\infty}}, \phi_{\infty})$, then we have $\tilde{\pi }_{1 }^{\infty}=\pi_{1 }^{\infty}$. Proceeding by induction, it is easy to conclude that $ \tilde{\pi}_{\alpha }^{\infty }=\pi_{\alpha }^{\infty }$
for all $\alpha $. This completes the proof of part (1) of the lemma.

\medskip

For part (2), notice that the argument given above applies to Seshadri
filtration as well, where because of the lack of uniqueness of Seshadri filtration we may
conclude only that the ranks and degrees of limiting filtration are same with that of the original filtration.

\hfill $\Box$ \\

\medskip

\medskip

{\bf Proposition 5.14. } {\it Let $(E, A_{0} , \phi_{0} )$ be a
Higgs bundle on a smooth K\"ahler manifold $(M , \omega )$, and satisfy the same assumptions as that in Theorem 5.12. Then given $\delta >0 $ and $1\leq p <\infty $,
$(E, A_{0} , \phi_{0} )$ has a $L^{p}$ $\delta$-approximate Hermitian structure.}

\medskip

{\bf Proof. } Let $(A_{t}, \phi_{t})$ be the solution of the Yang-Mills-Higgs flow (\ref{YMHH}) with initial data $(A_{0} , \phi_{0})\in \textbf{B}_{(E, H_{0})}$, and $H_{t}$ be the solution of the heat flow (\ref{D1}) on Higgs bundle $(E, A_{0}, \phi_{0})$ with initial data $H_{0}$. Let $(A_{\infty}, \phi_{\infty})$ be the Uhlenbeck limiting of some sequence $(A_{t_{j}}, \phi_{t_{j}})$. Applying the previous lemma, we have $\Psi ^{HN}_{\omega }((A_{t_{j}}, \phi_{t_{j}}), H_{0})\rightarrow  \Psi ^{HN}_{\omega }((A_{\infty}, \phi_{\infty }), H_{\infty })$ strongly in $L^{p}$ for all $1
\leq p < \infty $. By Corollary 2.8., we have
\begin{eqnarray}\label{5.56}
\begin{array}{lll}
&&\|\frac{\sqrt{-1}}{2\pi}\Lambda_{\omega }(F_{A_{H(t_{j})}}+[\phi_{0} , \phi_{0}^{\ast H(t_{j})
}])-\Psi_{\omega }^{HN}((E, A_{0}, \phi_{0}) , H(t_{j}))\|_{L^{p}(\omega )}\\
&&=\|\theta (A_{t_{j}}, \phi_{t_{j}})-\Psi_{\omega }^{HN}((A_{t_{j}}, \phi_{t_{j}}) , H_{0})\|_{L^{p}(\omega)}\\
&&=\|\theta (A_{t_{j}}, \phi_{t_{j}})-\theta (A_{\infty }, \phi_{\infty }))\|_{L^{p}(\omega)}\\ &&+\|\Psi ^{HN}_{\omega }((A_{\infty}, \phi_{\infty }), H_{\infty })-\Psi_{\omega }^{HN}(( A_{t_{j}}, \phi_{t_{j}}) , H_{0})\|_{L^{p}(\omega)} \rightarrow 0.\\
\end{array}
\end{eqnarray}

\hfill $\Box$ \\

\medskip

{\bf Theorem 5.15. } {\it Let $(E, A_{0} , \phi_{0} )$ be a
Higgs bundle on a smooth K\"ahler manifold $(M , \omega )$. Then given $\delta >0$ and $1\leq p <\infty $,
$(E, A_{0} , \phi_{0} )$ has a $L^{p}$ approximate Hermitian structure.}

\medskip

{\bf Proof. }
By Proposition 3.7, we can resolve the singularity set $\Sigma_{al}$ by blowing up finitely many times, i.e. we have a sequence of blow-ups:
\begin{eqnarray}
\pi_{i}: \overline{M}_{i}\rightarrow \overline{M}_{i-1}, \quad i=1, \cdots , r
\end{eqnarray}
where $\overline{M}_{0}=M$, such that every $\pi_{i}$ is blow up along a smooth complex submanifold, every $\overline{E}_{i}=\pi^{\ast} (E_{i-1})$ is bundle, and the pull back filtration $(\pi_{r}\circ \cdots \circ \pi_{1})^{\ast}(\textit{F}^{HNS}(E, A_{0}, \phi_{0}))$ of $\overline{E}_{r}$ is given by $\overline{\phi}_{r}$-invarient sub-bundles, where $\overline{\phi}_{i}=(\pi_{i}\circ \cdots \circ \pi_{1})^{\ast}(\phi_{0})$. On each blow-up $\overline{M}_{i}$, we have a family of K\"ahler metrics defined iteratively by $\omega_{\epsilon_{1}\cdots \epsilon_{i}}=\pi_{i}^{\ast}\omega_{\epsilon_{1}\cdots \epsilon_{i-1}}+\epsilon_{i}\eta_{i}$, where $\eta_{i}$ is a K\"ahler metric on $\overline{M}_{i}$. By  proposition 5.14., for any fixed small $\epsilon_{1}, \cdots , \epsilon_{r-1} $, and any $\delta '$, we have a metric $H$ on  bundle $E_{r-1}$ such that
\begin{eqnarray*}
\|\Theta_{\omega_{\epsilon_{1}, \cdots , \epsilon_{r-1}}} ((\overline{\partial }_{\overline{E}_{r-1}}, H), \overline{\phi}_{r-1})-\Psi ^{HN}_{\omega_{\epsilon_{1}, \cdots , \epsilon_{r-1}} }((\overline{\partial }_{\overline{E}_{r-1}}, \overline{\phi}_{r-1 }), H)\|_{L^{2}(\omega_{\epsilon_{1}, \cdots , \epsilon_{r-1}})}\leq \delta '.
\end{eqnarray*}

By induction, we can assume that, for any fixed small $\epsilon_{1} $ and any $\delta '$, we have a metric $H$ on  bundle $E_{1}$ such that
\begin{eqnarray}
\|\Theta_{\omega_{\epsilon_{1}}} ((\overline{\partial }_{\overline{E}_{1}}, H), \overline{\phi}_{1})-\Psi ^{HN}_{\omega_{\epsilon_{1}} }((\overline{\partial }_{\overline{E}_{1}}, \overline{\phi}_{1 }), H)\|_{L^{2}(\omega_{\epsilon_{1}})}\leq \delta '.
\end{eqnarray}
Since $\pi_{1}: M_{1}\rightarrow M$ is the blow-up along a smooth complex submanifold,
by Proposition 5.10, then for any $\delta >0$ and any $1\leq p < 1+\frac{1}{2k-1}$ there are $\epsilon_{1}>0$ and a smooth Hermitian metric $\overline{H}_{1}$ on $E_{1}$ such that
\begin{eqnarray}
\|\Theta_{\omega_{\epsilon}} ((\overline{\partial }_{\overline{E}_{1}}, H_{1}), \overline{\phi}_{1})-\Psi ^{HN}_{\omega_{\epsilon} }((\overline{\partial }_{\overline{E}_{1}}, \overline{\phi}_{1 }), H_{1})\|_{L^{p}(\omega_{\epsilon})}\leq \delta '.
\end{eqnarray}
for all $0<\epsilon \leq \epsilon_{1}$.  By Proposition 5.11, Theorem 5.12 and Proposition 5.14,  we see that for any given $\delta $ and $1\leq p <\infty $,
$(E, A_{0} , \phi_{0} )$ has a $L^{p}$ $\delta $-approximate Hermitian structure.

\hfill $\Box$ \\

\medskip

Then repeating the argument in Theorem 5.12, we have:

\hspace{0.3cm}

{\bf Theorem 5.16. } {\it Let $(A_{t} , \phi_{t})$ be a smooth
solution of the gradient flow (\ref{YMHH}) on the Hermitian vector bundle
$(E , H_{0})$ with initial condition $(A_{0} , \phi_{0})\in \textbf{B}_{(E, H_{0})}$, and $(A_{\infty } , \phi_{\infty })$ be a
Uhlenbeck limit. Let $E_{\infty }$ denote the  vector bundle
obtained from $(A_{\infty } , \phi_{\infty })$ as that in Proposition 2.10.
Then the Harder-Narasimhan type of the extended reflexive Higgs sheaf $(E_{\infty } , A_{\infty } ,
\phi_{\infty })$ is same as that of the original Higgs bundle $(E_{0 } , A_{0 } , \phi_{0
})$.}

\section{Proof of theorem 1.1. }
\setcounter{equation}{0}

Let $\{E_{\alpha , \beta}\}$ be the
HNS-filtration  of the Higgs bundle $(E ,
\overline{\partial}_{A_{0}}, \phi_{0})$,   the associated graded object $Gr^{HNS}(E, A_{0},
\phi_{0})=\oplus_{\alpha =1}^{l}\oplus_{\beta=1}^{r_{\alpha}}Q_{\alpha , \beta}$ be uniquely determined by the
isomorphism class of $(A_{0} , \phi_{0})$, where $Q_{\alpha , \beta }=E_{\alpha , \beta}/E_{\alpha , \beta -1}$. We refer to $\Sigma_{al}$ as the singular set of the double filtration $\{E_{\alpha , \beta }\}$, it is a complex analytic subset of $M$ of complex codimensional at least $2$. We will prove the result
inductively on the length of the HNS filtration. The inductive
hypotheses on a  sheaf $Q $ are following:

\hspace{0.3cm}

{\bf Inductive hypotheses:} {\it  There is a sequence of Higgs structures $(A_{ j}^{Q} ,
\phi_{ j}^{Q})$ on $Q$ such that:

(1) $(A_{ j}^{Q} ,
\phi_{ j}^{Q})\rightarrow (A_{ \infty}^{Q_{\infty }} , \phi_{ \infty}^{Q_{\infty }})$ in
$C^{\infty }_{loc}$ off $\Sigma_{al}\cup \Sigma_{an}$;

(2) $(A_{ j}^{Q} ,
\phi_{ j}^{Q})=g_{j}(A_{ 0}^{Q} ,
\phi_{ 0}^{Q})$
 for some
$g_{j}\in  \textbf{G}^{C}(Q)$;

(3) $(Q, \overline{\partial }_{A_{0}^{Q}}, \phi_{0}^{Q})$ and $(Q_{\infty }, \overline{\partial }_{A_{\infty}^{Q_{\infty}}}, \phi_{\infty }^{Q_{\infty}})$ extended to $M$ as reflexive Higgs sheaves with the same HN type;

 (4) $\| \phi_{j}^{Q}\|_{C^{0}}$ and $\|\sqrt{-1}\Lambda _{\omega }(F_{A_{ j}^{Q}})\|_{L^{1}(\omega)}$ is uniformly bounded  in $j$.}

\hspace{0.3cm}

Let $S=E_{1, 1}$ be the first stable Higgs sub-sheaf corresponding to the HNS-filtration of the Higgs bundle $(E, \overline{\partial }_{A_{0}}, \phi_{0})$,
 $\pi : \tilde{M}\rightarrow M$ be the resolution of singularities $\Sigma_{al}$ then the filtration of $\tilde{E}=\pi^{\ast }E $ is given by subbundles $\{\tilde{E}_{\alpha ,\beta }\}$, isomorphic to $\{E_{\alpha , \beta }\}$ off the exception divisor $\tilde{\Sigma}=\pi^{-1}(\Sigma_{al})$.
Setting $(\tilde{A}_{j}, \tilde{\phi}_{i})=\pi ^{\ast }(A_{j}, \phi_{j})$ and $\tilde{g}_{j}=\pi^{\ast}g_{j} $, then we have $(\tilde{A}_{j}, \tilde{\phi}_{i})=\tilde{g}_{j}(\tilde{A}_{0}, \tilde{\phi}_{0})$. By Theorem 2.7., we know that $(\tilde{A}_{j}, \tilde{\phi }_{j}
)$ converges to $(\tilde{A}_{\infty}, \tilde{\phi} _{\infty } )$   in $C^{\infty}_{loc}$ topology outside
$\pi^{-1}(\Sigma_{al}\cup \Sigma_{an})$. By Corollary 2.8. and the uniform $C^{0}$ bound on $\phi (t)$ (Lemma 2.3), we have
$|\sqrt{-1}\Lambda_{\omega_{0}}(F_{\tilde{A}_{j}})|_{L^{\infty}}$, specially  $\|\sqrt{-1}\Lambda _{\omega_{0}}(F_{\tilde{A}_{j}})\|_{L^{1}}$ is uniformly bounded in $j$, where $\omega_{0}=\pi^{\ast }\omega $.

 Using Proposition 4.1, we have a   subsequence of $\tilde{g}_{j} \circ i_{0}$, up to rescale,  converges to a nonzero smooth $\tilde{\phi }$-invariant holomorphic map
$\tilde{f}_{\infty }: \tilde{S}\rightarrow \tilde{E}_{\infty }$ off $\pi^{-1}(\Sigma_{al}\cup \Sigma_{an})$. Since $\tilde{S}$ is isomorphic to $S$ off the exception divisor, then we obtain a subsequence of $f_{j}=g_{j}\circ i_{0}$ up to rescale, which converges to a  nonzero smooth $\phi $-invariant holomorphic map
$f_{\infty }: S \rightarrow (E_{\infty}, \overline{\partial }_{A_{\infty}})$ in $C^{\infty}_{loc}$ on $M\setminus \Sigma_{an}\cup \Sigma_{al}$, where $i_{0}: S \rightarrow (E, \overline{\partial}_{A_{0}})$ is the holomorphic inclusion. By Hartog's
theorem, $f_{\infty}$ extends to a Higgs sheaf homomorphism  $f_{\infty }: (S, \phi_{0})\rightarrow (E_{\infty}, \overline{\partial }_{A_{\infty}}, \phi_{\infty})$ on $M$ (where $(E_{\infty}, \overline{\partial }_{A_{\infty}}, \phi_{\infty})$ is the extended reflexive Higgs sheaf).

As above,  $\pi_{1}^{(j)}$ denotes the projection
to $g_{j}(S)$. Since $\pi_{1}^{(j)}\circ f_{j}=f_{j}$, we see that in the limit $\pi_{1}^{\infty}\circ f_{\infty}=f_{\infty}$. By Lemma 5.13, we know that
$\pi_{1}^{\infty}$ determines a Higgs subsheaf $E_{1, 1}^{\infty}$ of $(E_{\infty}, \overline{\partial}_{A_{\infty}}, \phi_{\infty})$, with $rank (E_{1, 1}^{\infty})=rank (S)$ and $\mu (E_{1, 1}^{\infty})=\mu(S)$. Since  $(E_{\infty}, \overline{\partial}_{A_{\infty}}, \phi_{\infty})$ and  $(E_{0}, \overline{\partial}_{A_{0}}, \phi_{0})$ have the same HN type, thus we have the Higgs subsheaf $(E_{1, 1}^{\infty}, \phi_{\infty})$ is semistable and
\begin{eqnarray}f_{\infty}: S\rightarrow E_{1,1}^{\infty}.\end{eqnarray} Recall that  $S=E_{1, 1}$  is Higgs stable.
    By Lemma 3.4., we see that the non-zero holomorphic map
$f_{\infty }$ must be injective, then \begin{eqnarray}S\simeq
E_{1, 1}^{\infty}=f_{\infty }(S)\end{eqnarray}
 on $M\setminus (\Sigma_{al}\cup \Sigma_{an})$. It is easy to see that $E_{1, 1}^{\infty }$ is a stable Higgs subsheaf of $(E_{\infty}, \overline{\partial}_{A_{\infty}}, \phi_{\infty})$.

Let $\{e_{\alpha }\}$ be a local frame of $S$, and $H_{j, \alpha \bar{\beta}}=<f_{j}(e_{\alpha }), f_{j}(e_{\beta})>_{H_{0}}$. We can write the orthogonal projection $\pi_{1}^{(j)}$ as
\begin{eqnarray}
\pi_{1}^{(j)} (X)=<X, f_{j}(e_{\beta })>H_{j}^{\alpha , \bar{\beta} }f_{j}(e_{\alpha})
\end{eqnarray}
for any $X \in E$, where $(H_{j}^{\alpha , \bar{\beta} })$ is the inverse of the matrix $(H_{j, \alpha \bar{\beta}})$.
Because $f_{j}\rightarrow
f_{\infty}$ in $C^{\infty }(\Omega )$, and $f_{\infty}$ is injective,  then we can prove that $\pi_{1}^{(j)}\rightarrow
\pi_{1}^{\infty}$ in $C^{\infty }_{loc }$ off $\Sigma_{an}\cup \Sigma_{al}$.

\medskip

Let $Q=E/S$, then we have $Gr^{HNS}(E, \overline{\partial}_{A_{0}}, \phi_{0})=S\oplus Gr^{HNS}(Q, \overline{\partial}_{A_{0}^{Q}}, \phi_{0}^{Q} )$. Write the orthogonal holomorphic decomposition $(E_{\infty },
\overline{\partial }_{A_{\infty }}, \phi_{\infty})=E_{1}^{\infty }\oplus Q_{\infty}$, where $Q_{\infty}=(E_{1}^{\infty})^{\bot }$ because $H_{\infty}$ is admissible Hermitian-Einstein metric.
Using
 Lemma 5.12 in \cite{Da}, we can choose a sequence of unitary gauge
transformation $u_{j}$ such that
$\pi_{1}^{(j)}=u_{j}\tilde{\pi}_{j}u_{j}^{-1}$ where
$\tilde{\pi}_{j}(E)=\pi_{1}^{\infty }(E)=E_{1}^{\infty}$ and $u_{j}\rightarrow Id_{E}$
in $C^{\infty}(loc )$ on $M\setminus (\Sigma_{al}\cup \Sigma_{an})$. It is easy to check that $u_{j}(Q_{\infty})=(\pi_{1}^{(j)} (E))^{\bot}$. Noting the bundle isomorphisms $p^{\ast}: Q \rightarrow S^{\bot}$ and the unitary gauge
transformation $u_{0}: Q_{\infty}\rightarrow S^{\bot }$, and considering the induced
connections on $Q$, we have \begin{eqnarray}D_{A_{j}^{Q}}=u_{0}\circ
u_{j}^{-1}\circ \pi_{j}^{\bot}\circ D_{A_{j}} \circ \pi_{j}^{\bot} \circ
u_{j}\circ u_{0}^{-1},\end{eqnarray}
\begin{eqnarray}\phi_{j}^{Q}=u_{0}\circ u_{j}^{-1}\circ \pi_{j}^{\bot}\circ
\phi_{j} \circ \pi_{j}^{\bot} \circ u_{j}\circ u_{0}^{-1}\in
\Omega^{ 1, 0}(End (Q)),\end{eqnarray}
\begin{eqnarray}h_{j}=u_{0}\circ u_{j}^{-1}\circ \pi_{j}^{\bot}\circ
g_{j}\in  \textbf{G}^{C}(Q).\end{eqnarray}
Then, we have
\begin{eqnarray}
\begin{array}{lll}
\overline{\partial}_{A_{j}^{Q}}&=&u_{0}\circ u_{j}^{-1}\circ
\pi_{j}^{\bot}\circ \overline{\partial}_{j} \circ \pi_{j}^{\bot}
\circ u_{j}\circ u_{0}^{-1}\\
&=&u_{0}\circ u_{j}^{-1}\circ \pi_{j}^{\bot}\circ g_{j}\circ
\pi_{0}^{\bot}\circ \overline{\partial}_{0} \circ
\pi_{0}^{\bot}\circ g_{j}^{-1}
\circ u_{j}\circ u_{0}^{-1}\\
&=& h_{j}\circ \overline{\partial}_{A_{0}^{Q}} \circ h_{j}^{-1},
\end{array}
\end{eqnarray}
\begin{eqnarray}
\partial_{A_{j}^{Q}}=(h_{j}^{\ast})^{-1}\circ \partial_{A_{0}^{Q}}\circ h_{j}^{\ast},
\end{eqnarray}

\begin{eqnarray}
\begin{array}{lll}
\phi_{j}^{Q}&=&u_{0}\circ u_{j}^{-1}\circ \pi_{j}^{\bot}\circ
g_{j}\circ \phi_{0} \circ g_{j}^{-1}\circ \pi_{j}^{\bot} \circ
u_{j}\circ u_{0}^{-1}\\
&=&u_{0}\circ u_{j}^{-1}\circ \pi_{j}^{\bot}\circ g_{j}\circ
\pi_{0}^{\bot}\circ \phi_{0}\circ \pi_{0}^{\bot} \circ
g_{j}^{-1}\circ  u_{j}\circ u_{0}^{-1}\\
&=& h_{j}\circ \phi_{0}^{Q} \circ h_{j}^{-1},\\
\end{array}
\end{eqnarray}
and
\begin{eqnarray}
\overline{\partial}_{A_{j}^{Q}}\phi_{j}^{Q}=\pi_{0}^{\bot}\circ
(\overline{\partial}_{0}\circ \phi_{0}+\phi_{0}\circ
\overline{\partial}_{0})\pi_{0}^{\bot}=0,
\end{eqnarray}
 where we have used $h_{j}^{-1}=\pi_{0}^{\bot} \circ
g_{j}^{-1}\circ u_{j}\circ u_{0}^{-1}$. On the other hand, by the
definition, it is easy to check that $u_{0}^{\ast}(A_{j}^{Q}, \phi_{j}^{Q})\rightarrow (A_{\infty
}^{Q_{\infty}}, \phi_{\infty}^{Q_{\infty}})$ in
$C^{\infty}_{loc }$.
Now we check the third statement in the inductive hypotheses. Let's consider the Gauss-Codazzi equation on $(\pi_{1}^{(j)}(E))^{\bot}=Q^{j}$
\begin{eqnarray}
F_{A_{Q^{j}}}=(\pi_{1}^{(j)})^{\bot}\circ F_{A_{j}}+\partial_{A_{j}}\pi_{1}^{(j)}\wedge \overline{\partial}_{A_{j}}\pi_{1}^{(j)},
\end{eqnarray}
where $D_{A_{Q^{j}}}=(\pi_{1}^{(j)})^{\bot}\circ D_{A_{j}}$. Setting the Higgs field $\phi_{Q^{j}}=(\pi_{1}^{(j)})^{\bot}\circ \phi_{j}$, by (\ref{pi1}) and (\ref{5.56}), we have
\begin{eqnarray}
\begin{array}{lll}
&&\int_{M}|\sqrt{-1}\Lambda_{\omega}(F_{A_{j}^{Q}}+[\phi_{j}^{Q} , (\phi_{j}^{Q})^{\ast }
])-\Psi^{hn} ((A_{j}^{Q}, \phi_{j}^{Q}), H_{0})|\frac{\omega^{n}}{n!}\\
&=&\int_{M}|\sqrt{-1}\Lambda_{\omega}(F_{A_{Q^{j}}}+[\phi_{Q^{j}} , (\phi_{Q^{j}})^{\ast }
])-\Psi^{hn} ((A_{Q^{j}}, \phi_{Q^{j}}), H_{0})|\frac{\omega^{n}}{n!}\\
&=&\int_{M}|(\pi_{1}^{(j)})^{\bot}\{\sqrt{-1}\Lambda_{\omega}(F_{A_{j}}+[\phi_{j} , (\phi_{j})^{\ast }
])-\Psi^{hn} ((A_{j}, \phi_{j}), H_{0})\}(\pi_{1}^{(j)})^{\bot}\\ && +\sqrt{-1}\Lambda_{\omega }(\partial_{A_{j}}\pi_{1}^{(j)}\wedge \overline{\partial}_{A_{j}}\pi_{1}^{(j)})\\
& &- (\pi_{1}^{(j)})^{\bot}([\phi_{j}, \pi_{1}^{(j)}]\wedge \phi_{j}^{\ast}+([\pi_{1}^{(j)}, \phi_{j}])^{\ast}\wedge \phi_{j})
|\frac{\omega^{n}}{n!}\\
&\leq &  \int_{M}|\sqrt{-1}\Lambda_{\omega}(F_{A_{j}}+[\phi_{j} , (\phi_{j})^{\ast }
])-\Psi^{hn} ((A_{j}, \phi_{j}), H_{0})| +| \overline{\partial}_{A_{j}}\pi_{1}^{(j)}|^{2}\\
& &+|\phi_{j}||[\phi_{j}, \pi_{1}^{(j)}]|\frac{\omega^{n}}{n!}\\
&& \rightarrow 0.
\end{array}
\end{eqnarray}
Since  $C^{0}$ norm of $\phi_{j}$ is uniformly bounded, then $\| \phi_{j}^{Q}\|_{C^{0}}$ and $\|\sqrt{-1}\Lambda _{\omega }(F_{A_{ j}^{Q}})\|_{L^{1}(\omega)}$ is uniformly bounded  in $j$.  So, $(Q , A_{j}^{Q} , \phi_{j}^{Q})$ satisfy the
inductive hypotheses. Since we can resolve the singularity set $\Sigma_{al}$ by blowing up finitely many times with non-singular center, and the pulling back of the HNS filtration is given by sub-bundles. The sheaf $Q$  and every geometric objects which we considered  are induced by the HNS filtration, so their pulling back are all smooth. Using Proposition 4.1 again,  by induction  we have
\begin{eqnarray}
E_{\infty}\simeq Gr^{HNS
}(E , \overline{\partial }_{A_{0}} , \phi_{0})= \oplus_{i=1}^{l} \oplus_{j=1}^{r_{i}}Q_{i, j}
\end{eqnarray}
on $M\setminus (\Sigma_{al}\cup \Sigma_{an})$. By Proposition 2.7, we know that  $(E_{\infty}, \overline{\partial }_{A_{\infty}}, \phi_{\infty})$ can be extended to the whole $M$ as a reflexive Higgs sheaf. By the uniqueness of reflexive extension in \cite{Siu}, we know that there exists a sheaf isomorphism \begin{eqnarray}f: (E_{\infty}, \overline{\partial }_{A_{\infty}}, \phi_{\infty}) \rightarrow Gr^{HNS
}(E , \overline{\partial }_{A_{0}} , \phi_{0})^{\ast \ast }\end{eqnarray} on $M$.
 This completes the
proof of Theorem 1.1.

\hfill $\Box$ \\

\hspace{0.4cm}

\hspace{0.3cm}


\begin{thebibliography}{99}


\bibitem{AB} M.Atiyah and R.Bott, {\em The Yang-Mills equations over
Riemann surfaces},  Phil. Trans. Roy. Soc. London A {\bf
308}(1982), 524-615.

\bibitem{AG1} L.Alvarez-Consul and O. Garcis-Prada, {\em Dimensional
reduction, SL(2, C)-equivariant bundles and stable holomorphic
chains},  Int. J. Math., {\bf 2}(2001), 159-201.



\bibitem{Bi} O.Biquard, {\em On parabolic bundles over a complex surface}, J. London. Math. Soc., {\bf 53}(1996),
no.2, 302-316.

\bibitem{Br1} S.B.Bradlow, {\em Vortices in holomorphic line bundles over
closed K\"ahler manifolds},  Commun.Math.Phys. {\bf 135}(1990),
1-17.



\bibitem{BG} S.B.Bradlow and O. Garcia-Prada, {\em Stable triples,
equivariant bundles and dimensional reduction},  Math. Ann. {\bf
304 }(1996), 225-252.



\bibitem{BS} S.Bando and Y.T.Siu. {\em Stable sheaves and Einstein-Hermitian metrics}, in {\it Gemetry and Analysis on Complex Manifolds}, World Scientific, 1994, 39-50.

\bibitem{BT} P.D.Bartolomeis and G.Tian, {\em Stability of complex vector
bundles},  J.Differential Geometry, {\bf 43} (1996),
232-275.


\bibitem{Da} G.Daskalopoulos, {\em The topology of the space of stable bundles on a Riemann surface}, J.Differential Geom. {\bf 36}(1992), 699-746.

\bibitem{DW1} G.Daskalopoulos and R.Wentworth, {\em Convergence
properties of the Yang-Mills flow on K\"ahler surfaces}, J.
Reine Angew. Math. {\bf 575}(2004), 69--99.

\bibitem{D0} S.K.Donaldson, {\em A new proof of a theorem of Narasimhan and Seshadri}, J.Differential Geom., {\bf 18}(198), 279-315.

\bibitem{D1} S.K.Donaldson, {\em Anti-self-dual Yang-Mills connections
over complex algebraic surfaces and stable vector bundles},
Proc.London Math.Soc. {\bf 50}(1985),
1-26.


\bibitem{GH} P.Griffiths and J.Harris, {\em Principles of algebraic geometry}, Pure and Applied Mathematics, Wiley-Interscience, New York, 1978.


\bibitem{GP} O.Garcia-Prada, {\em Dimensional reduction of stable bundles,
vortices and stable pairs},  Int.J.Math. {\bf 5 }(1994), 1-52.

\bibitem{HT} Tamas Hausel and Michael Thaddeus, {\em Generators for the cohomology ring
of the moduli space of rank 2 Higgs bundles}, Proc. London Math. Soc. (3),
{\bf 88}(3):632¨C658, 2004.

\bibitem{Hio1} H.Hironaka, {\em Resolution of singularities of an algebraic variety over a field of characteristic zero}, Ann. of Math. (2) {\bf 79} no.1, 109-203

\bibitem{Hi} N.J.Hitchin, {\em The self-duality equations on a Riemann
surface},  Proc.London Math.Soc. {\bf 55}(1987), 59-126.

\bibitem{HT} M.C.Hong and G.Tian, {\em
Asymptotical behaviour of the Yang-Mills flow and singular
Yang-Mills connections}, Math. Ann. {\bf 330}(2004), no. 3,
441--472.



\bibitem{Ja3} A.Jacob, {\em The Yang-Mills flow and the Atiyah-Bott formula on compact K\"ahler manifolds}, arXiv:1109.1550.



\bibitem{Ko2} S.Kobayashi, {\em Differential geometry of complex vector bundles}, Publications
of the Mathematical Society of Japan, {\bf 15}. Kano Memorial Lectures, 5. Princeton
University Press, Princeton, NJ (1987).

\bibitem{LN1} J.Y.Li, {\em Hermitian-Einstein metrics and Chern number inequalities on parabolic stable bundles over K\"hler manifolds},
 Comm. Anal. Geom. {\bf 8}(2000), no. 3, 445--475.

\bibitem{LN2} J.Y.Li and M.S.Narasimhan,
{\em Hermitian-Einstein metrics on parabolic stable bundles},
 Acta Math. Sin. (Engl. Ser.) {\bf 15} (1999), no. 1, 93--114.


\bibitem{LY} J.Li and S.T.Yau, {\em Hermitian-Yang-Mills connection on non-K\"ahler manifolds},
 Mathematical aspects of string theory (San Diego, Calif., 1986), 560--573, {\it Adv. Ser. Math. Phys.}, 1, World Sci. Publishing, Singapore, 1987.



\bibitem{LZ1} J.Y.Li and X.Zhang, {\em The gradient flow of Higgs pairs}, J. Eur. Math. Soc., {\bf 13}(2011), 1373-1422.

\bibitem{LZ2} J.Y.Li and X.Zhang, {\em Existence of approximate Hermitian-Einstein structures on semi-stable Higgs bundles},  Calc. Var. Partial Differential Equations, {\bf 52}(2015),  no. 3-4, 783¨C795.

\bibitem{PL} Peter. Li,
{\em Geometry analysis},
 {\em Cambridge Studies in Advanced Mathematics},	(No. 134).

\bibitem{M} I. Mundet i Riera, {\em A Hitchin-Kobayashi correspondence for
K\"ahler fibrations},  J.reine angew. Math. {\bf 528}(2000),
41-80.

\bibitem{NS} M.S.Narasimhan and C.S.Seshadri, {\em Stable and unitary
vector bundles on compact Riemann surfaces},  Ann. of Math., {\bf
82} (1965) 540-567.


\bibitem{Sib} B.Sibley, {\em Asymptotics of the Yang-Mills flow for holomorphic vector bundles over K\"ahler manifolds: the canonical structure of the limit}, arXiv:1206.5491v1.

\bibitem{Si} C.T.Simpson, {\em Constructing variations of Hodge structures
using Yang-Mills connections and applications to uniformization},
 J.Amer.Math.Soc., {\bf 1,} (1988)867-918.

\bibitem{Siu} Y.T.Siu, {\em A Hartogs type extension theorem for coherent analytic sheaves}, Ann. of Math. (2) {\bf 93}(1971), no.1,
166-188.

\bibitem{Ta} T.Mochizuki, {\em Kobayashi-Hitchin correspondence for tame harmonic bundles and an application},
Ast\'erisque {\bf 309} (2006), ISBN: 978-2-85629-226-6, +117pp.



\bibitem{UY} K.K.Uhlenbeck and S.T.Yau, {\em On existence of
Hermitian-Yang-Mills connection in stable vector bundles},
Comm.Pure Appl.Math., {\bf 39S}(1986), 257-293.

\bibitem{Wi} G.Wilkin, {\em Morse
theory for the space of Higgs bundles}, Comm. Anal. Geom.,
{\bf 16(2)}(2008), 283-332.



\end{thebibliography}
\end{document}